\documentclass[11pt]{paper}
\usepackage{a4wide, amsthm, amssymb, amscd, mathptmx}

\DeclareSymbolFont{rsfs}{OMS}{rsfs}{m}{n}
\DeclareSymbolFontAlphabet{\mathscr}{rsfs}
\DeclareSymbolFont{bbold}{U}{dsrom}{m}{n}
\DeclareSymbolFontAlphabet{\mathbb}{bbold}

\renewcommand{\Bbb}{\mathbb}
\renewcommand{\frak}{\mathfrak}
\newcommand{\catqot}{/\hskip-3pt/}
\newcommand{\C}{{\Bbb C}}

\newcommand{\End}{\mathop{\rm End}}
\newcommand{\ev}{\mathop{\rm ev}}

\newcommand{\G}{{\cal G}}
\newcommand{\PGL}{\mathop{\rm PGL}}
\newcommand{\GL}{\mathop{\rm GL}}

\newcommand{\gr}{\mathop{\rm gr}}
\newcommand{\Hom}{\mathop{\rm Hom}}

\newcommand{\Imm}{\mathop{\rm Im}}

\newcommand{\id}{\mathop{\rm id}}
\newcommand{\Jac}{\mathop{\rm Jac}}
\newcommand{\EL}{{\cal L}}

\newcommand{\n}{{\cal N}}
\newcommand{\Oh}{{\cal O}}

\newcommand{\Pe}{{\Bbb P}}

\newcommand{\Q}{{\Bbb Q}}
\newcommand{\R}{{\Bbb R}}

\newcommand{\SL}{\mathop{\rm SL}}
\newcommand{\Spec}{\mathop{\rm Spec}}

\newcommand{\Z}{{\Bbb Z}}
\newcommand{\N}{{\Bbb N}}
\newcommand{\la}{\lambda}
\newcommand{\lra}{\longrightarrow}

\newcommand{\lma}{\longmapsto}
\newcommand{\p}{\prime}
\newcommand{\q}{\quad}

\renewcommand{\phi}{\varphi}
\newcommand{\rk}{\mathop{\rm rk}}
\newcommand{\eps}{\varepsilon}

\newcommand{\ul}{\underline}
\newcommand{\ol}{\overline}

\theoremstyle{plain}
\newtheorem{Thm}{Theorem}[section]
\newtheorem{Cor}[Thm]{Corollary}
\newtheorem{Prop}[Thm]{Proposition}
\newtheorem{Lem}[Thm]{Lemma}

\theoremstyle{remark}
\newtheorem{Rem}[Thm]{Remark}

\newtheorem{Ex}[Thm]{\it Example}

\begin{document}
\pagestyle{myheadings}
\markboth{\hfill\rm Alexander H.W.\ Schmitt}{\rm Decorated Vector Bundles}

\title{\sf A universal construction for moduli spaces of decorated
vector bundles over curves}

\author{Alexander H.W.\ Schmitt
\institution{\rm \small Universit\"at Duisburg-Essen\\ \rm \small FB6 Mathematik \& Informatik\\ \small\rm  D-45117 Essen\\ \rm \small Germany
\\ \small\tt alexander.schmitt@uni-essen.de\rm}}
\date{}
\maketitle
\begin{abstract}
Let $X$ be a smooth projective curve over the field of complex numbers, and fix a homogeneous
representation
$\rho\colon \GL(r)\lra \GL(V)$. Then, one can associate to every vector bundle $E$ of rank $r$
over $X$ a vector bundle $E_\rho$ with fibre $V$. We would like to study triples $(E,L,\phi)$
where $E$ is a vector bundle of rank $r$ over $X$, $L$ is a line bundle over $X$, and $\phi\colon
E_\rho\lra L$ is a non-trivial homomorphism. This set-up comprises well-known objects such as
framed vector bundles, Higgs bundles, and conic bundles. In this paper, we will formulate a general
(parameter dependent)
semistability concept for such triples, which generalizes the classical Hilbert-Mumford criterion,
and establish the existence of moduli spaces for the semistable objects. In the examples which have been
studied so far, our semistability concept reproduces the known ones. Therefore, our results give in
particular a unified construction for many moduli spaces considered in the literature.
\end{abstract}

\tableofcontents
\section*{Introduction}
The present paper is devoted to the study of vector bundles
with an additional structure from a unified point of view.
We have picked the name ``decorated vector bundles"
suggested in \cite{HL}.
\par
Before we outline our paper, let us give some background. The
first problem to treat is the problem of classifying vector
bundles over an algebraic curve $X$, assumed here to be smooth,
projective and defined over $\C$. From the point of view of
projective geometry, this is important because it is closely
related to classifying projective bundles over $X$, so-called \it
ruled manifolds\rm. The basic invariants of a vector bundle $E$
are its rank and its degree. They determine $E$ as topological
$\C$-vector bundle. The problem of classifying all vector bundles
of fixed degree $d$ and rank $r$ is generally accessible only in
a few cases:
\begin{itemize}
\item The case $r=1$, i.e., the case of line bundles which is
      covered by the theory of Jacobian varieties.
\item The case $X=\Pe_1$ where Grothendieck's splitting theorem \cite{Gr}
      provides the classification.
\item The case $g(X)=1$. In this case, the classification has been
      worked out by Atiyah \cite{At}.
\end{itemize}
As is clear from the theory of line bundles, over a curve of
genus $g\ge 1$, vector bundles of degree $d$ and rank $r$ cannot
be parameterized by discrete data. Therefore, one seeks a variety
parameterizing all vector bundles of given degree $d$ and rank $r$
characterized by a universal property like the Jacobian. Such a
universal property was formulated by Mumford in his definition of
a \it coarse moduli space \rm \cite{GIT}. However, one checks
that the family of all vector bundles of degree $d$ and rank $r$
is not bounded which implies that a coarse moduli space cannot
exist. For this reason, one has to restrict one's attention to
suitable bounded subfamilies of the family of all vector bundles
of degree $d$ and rank $r$. Motivated by his general procedure to
construct moduli spaces via his Geometric Invariant Theory
\cite{GIT}, Mumford suggested that these classes should be the
classes of stable and semistable vector bundles. His definition,
given in \cite{Mu2}, is the following: A vector bundle $E$ is
called \it (semi)stable\rm, if for every non-trivial, proper
subbundle $F\subset E$
$$
\mu(F):={\deg F\over\rk F}\q (\le)\q \mu(E).
$$
Here, ``($\le$)" means that ``$\le$" is to be used for defining
``semistable" and ``$<$" for stable. Seshadri then succeeded to
give a construction of the coarse moduli space of stable vector
bundles, making use of Geometric Invariant Theory \cite{Sesh}.
This moduli space is only a quasi-projective manifold. To
compactify it, one has also to look at semistable vector bundles.
Seshadri formulated the notion of \it S-equivalence \rm of
semistable bundles which agrees with isomorphy for stable bundles
but is coarser for properly semistable ones. The moduli space of
S-equivalence classes exists by the same construction and is a
normal projective variety compactifying the moduli space of
stable bundles. Later Gieseker, Maruyama, and Simpson generalized
the results to higher dimensions \cite{Gies}, \cite{Ma},
\cite{Si}. Their constructions also apply to curves and replace
Seshadri's (see \cite{LP}). Narasimhan and Seshadri related
stable bundles to unitary representations of fundamental groups,
a framework in which vector bundles had been formerly studied
\cite{NS1}, \cite{NS2}.
\par
The next step is to consider vector bundles with extra structures.
Let us mention a few sources for this kind of problems:
\begin{itemize}
\item \bf Classification of algebraic varieties\rm.
     We have already mentioned that the classification of vector bundles
     is related to the classification of projective bundles via
     the assignment $E\lma \Pe(E)$. Suppose, for example, that we want
     to study divisors in projective bundles. For this, let
     $E$ be a vector bundle, $\Pe(E)$ its associated projective bundle,
     $k$ a positive integer, and $M$ a line bundle on $X$.
     To give a divisor $D$ in the linear system $|\Oh _{\Pe(E)}(k)\otimes\pi^* M|$
     we have to give a section $\sigma\colon \Oh _{\Pe(E)}\lra
     \Oh _{\Pe(E)}(k)\otimes\pi^* M$ which is the same as giving a non-zero
     homomorphism $\Oh _X\lra S^k E\otimes M$, or $S^k E^\vee\lra M$.
     Thus, we are led to classify triples $(E,M,\tau)$ where
     $E$ is a vector bundle over $X$, $M$ a line bundle, and $\tau\colon
     S^k E\lra M$ a non-trivial homomorphism.
     In case the rank of $E$ is three and $k$ is two, this is the theory
     of \it conic bundles\rm, recently studied by G\'omez and Sols \cite{GS}.
\item \bf Dimensional reduction\rm.
     Here, one looks at vector bundles $\G$ on $X\times \Pe_1$ which can
     be written as extensions
     $$
     0\lra \pi_X^*F\lra \G\lra \pi_X^*E\otimes\pi_{\Pe_1}^*\Oh _{\Pe_1}(2)
     \lra 0
     $$
     where $E$ and $F$ are vector bundles on $X$. These extensions are
     parameterized by $H^0(E^\vee\otimes F)=\Hom(E,F)$.
     The study of such vector bundles is thus related to the study
     of triples $(E,F,\phi)$ where $E$ and $F$ are vector bundles
     on $X$ and $\phi\colon E\lra F$ is a non-zero homomorphism.
     These are the \it holomorphic triples \rm of Bradlow and Garc\'\i a-Prada
     \cite{HT1} and \cite{HT2}. They were also studied from the algebraic
     point of view by the author \cite{Sch5}.
     For the special case $E=\Oh _X$, we find the
     problem of vector bundles with a section, so-called \it Bradlow pairs
     \rm \cite{Br}. An important application of Bradlow pairs was given by
     Thaddeus in his proof of the Verlinde formula \cite{Th}.
\item \bf Representations of fundamental groups\rm.
     \it Higgs bundles \rm are pairs $(E,\phi)$, consisting of a vector bundle
     $E$ and a twisted endomorphism $\phi\colon E\lra E\otimes\omega_X$.
     Simpson used in \cite{Si} the higher dimensional analogues of these
     objects to study representations of fundamental groups of projective
     manifolds. This ties up nicely with the work of Narasimhan and
     Seshadri.
\item \bf Gauge theory\rm.
     Here, one starts with differentiable vector bundles together with
     an additional structure and considers certain differential equations
     associated to these data. The solutions of the equations
     then have --- via a Kobayashi-Hitchin correspondence --- interpretations
     as holomorphic decorated vector bundles over $X$, satisfying
     certain  stability conditions. Again, the first case where this
     arose was the theory of Hermite-Einstein equations and stable vector
     bundles (see \cite{LT}) and was later studied in more complicated
     situations like the above examples. Recently, Banfield \cite{Ba}
     and Mundet i Riera \cite{MR} investigated this in a broad context.
     We will come back to this again.
\end{itemize}
Now, for all of these problems and many more, there exist notions
of semistability, depending on a rational parameter. The task of projective
geometry is then to generalize the construction of Seshadri and the
successors to obtain moduli spaces for the respective semistable and
stable objects. These constructions, where existent, were done
case by case and follow a certain
pattern inspired by Gieseker's, Maruyama's, and Simpson's constructions.
One is therefore led to ask for a single unifying construction
incorporating the known examples. This would complete the algebraic
counterpart to the work of Banfield and Mundet i Riera.
\par
We will consider this problem in the present article. Our framework
is as follows: We fix a representation $\rho\colon \GL(r)\lra
\GL(V)$, such that the restriction to the centre
$\C^*\subset\GL(r)$ is $z\lma z^\alpha\cdot\id_V$ for some integer
$\alpha$. Then, to any vector bundle $E$, we can associate a
vector bundle $E_\rho$ of rank $\dim V$. The objects we will
treat are triples $(E,M,\tau)$ where $E$ is a vector bundle of
rank $r$, $M$ is a line bundle, and $\tau\colon E_\rho\lra M$ is
a non-zero homomorphism. E.g., for $\rho\colon \GL(3)\lra
\GL(S^2\C^3)$, we recover conic bundles. The list of problems we
then have to solve is
\begin{itemize}
\item Formulate an appropriate notion of semistability for the above
      objects!
\item Prove boundedness of the semistable triples $(E,M,\tau)$ where
      $\deg E$ and $\deg M$ are fixed!
\item Construct a parameter space ${\frak P}$ for the semistable objects together
      with an action of a general linear group $G$, such that the equivalence
      relation induced by this action is the natural equivalence relation
      on those triples!
\item Show that the categorical quotient ${\frak P}\catqot G$ exists!
\end{itemize}
The latter space will then be the moduli space. As one sees from this list,
especially in view of the existing constructions,
Geometric Invariant Theory will play a central r\^ole.
Let us explain how one can find the semistability concept.
First, assume that we are given a bounded family of triples $(E,M,\tau)$.
Using the theory of quot-schemes it is by now not too hard a task
to construct a parameter space ${\frak P}$ for the members of the family in such a way
that we have a group action  as required together with a family of linearizations
--- depending on a rational parameter --- in line bundles over ${\frak P}$.
Therefore, we have realized the input for the GIT process.
The Hilbert-Mumford criterion now tells us how to find the semistable
points.
Thus, it is clear that our notion of semistability should mimic the
Hilbert-Mumford criterion as closely as possible.
Such an approach was also taken in gauge theory \cite{Ba} and \cite{MR}.
The structure of one parameter subgroups of the special linear group
suggests that one parameter subgroups should be replaced by \it weighted
filtrations of vector bundles\rm.
For weighted filtrations, one then defines the necessary numerical
quantities resembling Mumford's ``$\mu$" and arrives at the desired
semistability concept.
\par
Our paper is organized as follows: In the first section, we collect the
necessary background material from representation theory and GIT.
Then, we come to the definition
of semistability for the triples $(E,M,\tau)$ which depends on a
positive rational parameter and describe the associated
moduli functors. We state the main result, namely the existence of
moduli spaces, and proceed to the proofs along the lines outlined before.
The paper concludes with a long discussion of examples in order to
show that the known problems in that context can be recovered from
our results and that, in some cases, additional light is shed
on them.
The reader will notice that our general
semistability concept is in the known
cases more complicated than the existing ones and has to be simplified
to recover the known ones.
This is one of the key points of the paper:
The notion of semistability should be simplified after doing the GIT
construction and not before. This is why a unifying construction is
feasible. However, we will present a general method to simplify the
semistability concept in terms of the representation $\rho$. This method
enables us to write down in every concrete situation the semistability
concept in a more classical form. Applying this procedure, e.g., to
framed bundles or conic bundles immediately reproduces the known
semistability concepts. This provides us with a mechanism for finding the
correct notion of semistability without guessing or referring to gauge
theory.
\par
Finally, we remark that we have confined ourselves to the case of
curves in order to have a nice moduli functor associated to every
representation of the general linear group. However, if one
restricts to direct sums of tensor powers, the construction can
also be performed over higher dimensional manifolds \cite{GS2}.
These higher dimensional versions have, in turn, important
applications in the problem of compactifying moduli spaces of
principal bundles with \sl singular objects \rm (\cite{SchPrin},
\cite{GS3}). Finally, there is now also a version for product
groups $\GL(r_1)\times\cdots\times\GL(r_s)$ over base manifolds
of arbitrary dimension \cite{SchQuiv} the construction of which
is based on the results of the present paper.
\subsection*{Acknowledgments}
In November 1999, Professor Okonek organised a workshop on the
gauge theoretic aspects of decorated vector bundles. I would like to thank
him for the invitation to that workshop and financial support through his
Swiss National Funds project and him and the other
speakers, Professor L\"ubke and Professor Andrei Teleman, for
their enlightening talks which provided crucial impulses
for the present work
\par
This paper is an English version of Chapter I-III of the author's
work ``Universelle Konstruktionen f\"ur Modulr\"aume dekorierter
Vektorb\"undel \"uber Kurven" which was accepted as habilitation
thesis by the University of Essen.
The author acknowledges support by the DFG through the ``Schwerpunkt" program
``Globale Methoden in der komplexen Geometrie --- Global Methods in Complex Geometry''.
\subsection*{Notations and conventions}
\begin{itemize}
\item
All schemes will be defined over the field of complex numbers,
$X$ will be a smooth projective curve of genus $g\ge 2$.
We denote by $\ul{\mathop{\rm Sch}}_\C$ the category of separated schemes
of finite type over $\C$. A \it point \rm will be a closed point unless
otherwise mentioned.
\item
For a vector bundle $E$ over a scheme $S$,
we denote by $\Pe(E)$ the projective bundle of
hyperplanes
in the fibres of $E$.
\item
Given a product $X\times Y$ of schemes, $\pi_X$ and $\pi_Y$
stand for the projections from $X\times Y$ onto the respective factors.
\item
Let $V$ be a finite dimensional $\C$-vector space and $\rho\colon G\lra
\GL(V)$ a representation of the algebraic group $G$.
This yields an action of $G$ on $\Pe(V)$ and a linearization
$G\times\Oh _{\Pe(V)}(1)\lra \Oh _{\Pe(V)}(1)$. We will denote this
linearization again by $\rho$.
\item
Let $E$ be a vector bundle of rank $r$. Then, the associated
$\GL(r)$-principal bundle is given as ${\frak P}(E) =\bigcup_{x\in
X}\mathop{\rm Isom}(\C^r,E_x)\subset \ul{\Hom}(\Oh _X^{\oplus r},
E)$. If we are furthermore given an action $\Gamma\colon
\GL(r)\times F\lra F$ of $\GL(r)$ on a quasi-projective manifold
$F$, we set ${\frak P}(E)\times^{\GL(r)} F:= ({\frak P}(E)\times
F)/\GL(r)$. Here, $\GL(r)$ acts on ${\frak P}(E)\times F$ by
$(x,y)\cdot g= (x\cdot g, g^{-1}\cdot y)$. If $F$ is a vector
space and the action $\Gamma$ comes from a representation
$\rho\colon \GL(r)\lra \GL(F)$, we write $E_\rho$ for the vector
bundle ${\frak P}(E)\times^{\GL(r)} F$.
\item
For any $x\in\R$, we set $[x]_+:=\max\{\,0,x\,\}$.
\end{itemize}

\section{Preliminaries}
\subsection{Representations of the general linear group}
\label{Rep}
First, let $\rho\colon \GL(r)\lra \GL(V)$ be an irreducible representation
on the finite dimensional $\C$-vector space $V$.
\begin{Thm}
There are integers $a_1,...,a_r$ with $a_i\ge 0$ for
$i=1,...,r-1$, such that $\rho$ is a direct summand of the natural
representation of $\GL(r)$ on $$
S^{a_1}\bigl(\C^r\bigr)\otimes\cdots\otimes
S^{a_{r-1}}\bigl(\bigwedge^{r-1}
\C^r\bigr)\otimes\bigl(\bigwedge^r\C^r\bigr)^{\otimes a_n}.
$$
\end{Thm}
\begin{proof} 
See \cite{FH}, Proposition 15.47.
\end{proof}
\par
For any vector space $W$, the representations of $\GL(W)$ on
$S^i(W)$ and $\bigwedge^i W$ are direct summands of the representation
of $\GL(W)$ on $W^{\otimes i}$. Setting $a:=a_1+\cdots +a_{r-1}(r-1)$
and $b:=a_n$, we see that $\rho$ is a direct summand of the
representation $\rho_{a,b}$ of $\GL(r)$ on $(\C^r)^{\otimes a}\otimes
(\bigwedge^r\C^r)^{\otimes b}$.
\begin{Cor}
Let $\rho\colon \GL(r)\lra \GL(V)$ be a (not necessarily irreducible)
representation of $\GL(r)$ on the finite dimensional $\C$-vector space
$V$, such that the centre $\C^*\subset \GL(r)$ acts by $z\lma z^\alpha\cdot
\id_V$ for some $\alpha\in\Z$.
Then, there exist $a,b,c\in\Z_{\ge 0}$, $c>0$, such that $\rho$
is a direct summand of the natural representation $\rho_{a,b,c}$
of $\GL(r)$ on
$$
V_{a,b,c}\q:=\q \Bigl(\bigl(\C^r)^{\otimes a}\otimes
\bigl(\bigwedge^r\C^r\bigr)^{\otimes -b}\Bigr)^{\oplus c}.
$$
\end{Cor}
\begin{proof}
We can decompose $\rho=\rho_1\oplus\cdots\oplus\rho_c$ where
the $\rho_i$'s are irreducible representations.
By what we have said before, there are integers $a_i,b_i$, $i=1,...,c$,
with $a_i\ge 0$, $i=1,...,c$, such that $\rho$ is a direct summand
of $\rho_{a_1,b_1}\oplus\cdots\oplus\rho_{a_c,b_c}$.
Our assumption on the action of $\C^*$ implies that
$a_1+rb_1=\cdots=a_c+rb_c$.
Let $b$ be a positive integer which is so large that $b_i+b>0$ for
$i=1,...,c$. Then, $\rho_{a_i,b_i}$ is the natural representation
of $\GL(r)$ on
$$
\bigl(\C^r\bigr)^{\otimes a_i}\otimes \bigl(\bigwedge^r\C^r\bigr)^{\otimes
b_i+b}\otimes \bigl(\bigwedge^r\C^r\bigr)^{\otimes-b},\q i=1,...,c.
$$
Now, the $\GL(r)$-module
$$
\bigl(\C^r\bigr)^{\otimes a_i}\otimes \bigl(\bigwedge^r\C^r\bigr)^{\otimes
b_i+b}
$$
is a direct summand
of
$
\bigl(\C^r\bigr)^{\otimes a},\ a:=a_1+r(b_1+b)=\cdots=a_c+r(b_c+b),
$
and we are done.
\end{proof}
\subsection{Basic concepts from GIT}
We briefly summarize the main steps in Geometric Invariant Theory
to fix the notation.
References are \cite{GIT} and \cite{Ne}.
\subsubsection{The GIT-process}
\label{git}
Let $G$ be a reductive algebraic group and $G\times F\lra F$ an
action of $G$ on the projective scheme $F$. Let $L$ be an ample line
bundle on $F$. A \it linearization \rm of the given action in $L$ is a
lifting of that action to an action $\rho\colon G\times L\lra L$,
such that for every $g\in G$ and $x\in F$ the induced map
$L_x\lra L_{g\cdot x}$ is a linear isomorphism.
Taking tensor powers, $\rho$ provides us with linearizations
of the action in any power $L^{\otimes k}$, $k>0$, and actions of $G$
on $H^0(F,L^{\otimes k})$ for any $k>0$.
A point $x_0\in F$ is called \it semistable\rm, if there exist an
integer $k>0$ and a $G$-invariant section $\sigma\in H^0(F, L^{\otimes k})$
not vanishing in $x_0$. If, moreover, the action of $G$ on the
set $\{\, x\in F\,|\, \sigma(x)\neq 0\,\}$ is closed
and $\dim G\cdot x_0=\dim G$, $x_0$ is called
\it stable\rm.
The sets $F^{(s)s}$ of (semi)stable points are open $G$-invariant
subsets of $F$.
Finally, a point $x\in F$ is called \it polystable\rm, if it is semistable
and its $G$-orbit is closed in $F^{ss}$.
Using this definition, the stable points are precisely the polystable
points with finite stabilizer.
The core of Mumford's Geometric Invariant Theory is that the
categorical quotients $F^{ss}\catqot G$ and $F^s\catqot G$ do exist
and that $F^{ss}\catqot G$ is a projective scheme whose closed
points are in one to one correspondence to the orbits
of polystable points, so that $F^s\catqot G$ is in particular an
orbit space.
\par
A finite dimensional representation $\rho\colon G\lra \GL(V)$
provides an action of $G$ on $\Pe(V)$ and a linearization of this
action in $\Oh _{\Pe(V)}(1)$, called again $\rho$. A point
$[v]\in\Pe(V)$ represented by $v\in V^\vee$ is then semistable if
and only if the closure of the orbit of $v$ in $V^\vee$ does not
contain $0$, stable if, furthermore, its orbit is closed and the
dimension of this orbit equals the dimension of $G$, and
polystable if the orbit of $v$ in $V^\vee$ is closed.
\subsubsection{Around the Hilbert-Mumford criterion}
\label{HM}
Let $F$ be a projective variety on which the reductive group $G$ acts.
Suppose this action is linearized in the line bundle $L$.
Call the linearization $\rho$.
Then, given a one parameter subgroup $\la$ of $G$ and $x\in F$, we can form
$$
x_\infty\q:=\q \lim_{z\lra\infty} \la(z)\cdot x.
$$
The point $x_\infty$ is clearly a fix point for the $\C^*$-action on $F$
induced by $\la$. Thus, $\C^*$ acts on the fibre of $L$ over
$x_\infty$, say, with weight $\gamma$.
One defines
$$
\mu_\rho(\la, x)\q:=\q-\gamma.
$$
\begin{Thm}[Hilbert-Mumford criterion~\cite{GIT}]
A point $x\in F$ is (semi)\allowbreak sta\-ble, if and only if for
every non-trivial one parameter subgroup $\la\colon \C^*\lra G$
$$
\mu_\rho\bigl(\la, x\bigr)\q (\ge)\q 0.
$$
Moreover, a point $x\in F$ is polystable if and only if it is semistable and, for
every one parameter subgroup $\la$ of $G$ with $\mu_\rho(\la,x)=0$, there is a $g\in G$
with $x_\infty=g\cdot x$.
\end{Thm}
As we have explained in the introduction, our concept of stability
for decorated vector bundles is basically  a Hilbert-Mumford criterion.
To define the necessary numerical invariants, we need the following
preparatory
\begin{Lem}
\label{Def}
Let $S$ be a scheme and $\sigma\colon S\lra F$ a morphism.
Suppose the $G$-action on $F$ is linearized in the ample line
bundle $L$.
Then
$$
\mu_{\rho}(\la, \sigma)\q :=\q \max\bigl\{\,\mu_{\rho}(\la,\sigma(s))\,|\,
s\in S\,\bigr\}
\hbox{\q exists}.
$$
\end{Lem}
\begin{proof}
We may assume that $L$ is a very ample
line bundle.
Set $V:=H^0(F,L)$. The linearization $\rho$ provides us with
a representation $\rho\colon G\lra \GL(V)$ and a $G$-equivariant
embedding $\iota\colon F\hookrightarrow \Pe(V)$.
Since obviously $\mu_{\rho}(\la,x)=\mu_{{\id}_{\GL(V)}}(\la,\iota(x))$
for all points $x\in F$ and all one parameter subgroups $\la$
of $G$, we can assume $F=\Pe(V)$.
Now, there are a basis $v_1,...,v_n$ of $V$ and integers
$\gamma_1\le\cdots\le \gamma_n$ with
$$
\la(z)\cdot \sum_{i=1}^n c_i v_i\q=\q \sum_{i=1}^n z^{\gamma_i}c_iv_i.
$$
A point $[l]\in \Pe(V)$ can be thought of as the equivalence
class of a linear form
$$
l\colon V\lra \C.
$$
Then,
$$
\mu_\rho\bigl(\la, [l]\bigr)\q=\q -\min\bigl\{\, \gamma_i\,|\, l(v_i)\neq 0\,\bigr\}.
$$
Therefore, $\mu_\rho(\la,\sigma(s))\in\{\,-\gamma_1,...,-\gamma_n\,\}$,
and this
implies the assertion.

\end{proof}
\begin{Rem}
\label{SemCont}
Let $F\subset\Pe(V)$ and $\la$ a one parameter subgroup of $G$.
Choose a basis $v_1,...,v_n$ of $V$ and $\gamma_1\le\cdots\le \gamma_n$
as before. Suppose $\mu_{\rho}(\la, \sigma)=-\gamma_{i_0}$ and
let $V^0\subset V$ be the eigenspace for the weight $\gamma_{i_0}$.
Let $U\subset S$ be the open set where the rational
map $S\stackrel{\sigma}{\lra} F\hookrightarrow \Pe(V)\dasharrow \Pe(V^0)$
is defined. Then $\mu_{\rho}(\la, \sigma(s))=-\gamma_{i_0}$
for all $s\in U$.
In other words, if $S$ is irreducible,
$\mu_{\rho}(\la,\sigma)$ is just the generic weight
occurring for a point $\sigma(s)$, $s\in S$.
\end{Rem}
\subsubsection{Semistability for actions coming from direct sums
of representations}
Let $G$ be a reductive algebraic group and $V_1$,...,$V_s$ finite
dimensional vector spaces.
Suppose we are given representations $\rho_i\colon G\lra \GL(V_i)$,
$i=1,...,s$. The direct sum $\rho_1\oplus\cdots\oplus\rho_s$
provides us with a linear action of $G$ on $\Pe(V)$,
$V:=V_1\oplus\cdots
\oplus V_s$.
Furthermore, for any $\underline{\iota}=(\iota_1,...,\iota_t)$
with $0<t\le s$, $\iota_1,...,\iota_t\in\{\,1,...,s\,\}$,
and $\iota_1<\cdots<\iota_t$, the $\rho_i$'s yield an action
$\sigma_{\underline{\iota}}$ of $G$ on
$\Pe_{\underline{\iota}}:=
\Pe(V_{\iota_1})\times\cdots\times \Pe(V_{\iota_t})$,
and, for any sequence of positive integers $k_1,...,k_t$, a linearization
of $\sigma_{\underline{\iota}}$ in the very ample line bundle
$\Oh (k_1,...,k_t)$.
The computation of the semistable points in $\Pe(V)$ can be
reduced to the computation of the semistable points in the
$\Pe_{\underline{\iota}}$'s by means of the following
\begin{Thm}
\label{ch2} Let $w^\p:=([w_{\iota_1}, w_{\iota_2}],
[w_{\iota_3}],..., [w_{\iota_t}])$ be a point in the space
$\Pe(V_{\iota_1}\oplus V_{\iota_2}) \times\allowbreak
\Pe_{(\iota_3,...,\iota_t)}$. Then $w^\p$ is semistable
(polystable) w.r.t.\ the given linearization in the line bundle
$\Oh (k, k_3,\allowbreak ...,k_t)$, if and only if either
$([w_{\iota_i}], [w_{\iota_3}],..., [w_{\iota_t}])$ is semistable
(polystable) in $\Pe_{(\iota_i,\iota_3,...,\iota_t)}$ w.r.t.\ the
linearization in $\Oh (k,k_3,\allowbreak ...,k_t)$ for either $i=1$
(and $w_{\iota_2}=0$) or $i=2$ (and $w_{\iota_1}=0$), or there
are positive natural numbers $n$, $k_1$, and $k_2$, such that
$k_1+k_2=nk$ and the point $([w_{\iota_1}], [w_{\iota_2}],
[w_{\iota_3}],..., [w_{\iota_t}])$ is semistable (polystable) in
$\Pe_{(\iota_1,\iota_2,\iota_3,...,\iota_t)}$ w.r.t.\ the
linearization in $\Oh (k_1,k_2,nk_3,...,nk_t)$.
\end{Thm}
\begin{Rem}
As one easily checks, for stable points only the ``if''-direction
remains true.
\end{Rem}
\begin{proof} 
This theorem can be proved with the methods
developed in \cite{OST} for $s=2$. A more elementary approach is
contained in the note~\cite{Sch4}.
\end{proof}
\subsection{One parameter subgroups of $\SL(r)$}
\label{1PSG}
Let $\GL(r)\times F\lra F$ be an action of the general linear
group on the projective manifold $F$.
For our definition of semistability, only
the induced action of $\SL(r)\times F\lra F$ will matter.
Since the Hilbert-Mumford criterion will play a central r\^ole
throughout our considerations, we will have
to describe the one parameter subgroups of $\SL(r)$.
\par
Given a one parameter subgroup $\la\colon\C^*\lra \SL(r)$,
we can find a basis $\ul{w}=(w_1,...,w_r)$ of $\C^r$ and
a weight vector $\ul{\gamma}=(\gamma_1,...,\gamma_r)$ with integral entries,
such that
\begin{itemize}
\item $\gamma_1\le\cdots\le \gamma_r$ and $\sum_{i=1}^r\gamma_i=0$, and
\item $\la(z)\cdot \sum_{i=1}^rc_iw_i=\sum_{i=1}^r z^{\gamma_i}c_iw_i$.
\end{itemize}
Conversely, a basis $\ul{w}$ of $\C^r$ and a weight vector $\ul{\gamma}$
with the above properties define a one parameter subgroup
$\la(\ul{w},\ul{\gamma})$ of $\SL(r)$.
\par
To conclude, we remark that, for any vector $\ul{\gamma}=(\gamma_1,...,\gamma_r)$
of integers with $\gamma_1\le\cdots\le\gamma_r$ and $\sum\gamma_i=0$,
there is a decomposition
$$
\ul{\gamma}\q=\q \sum_{i=1}^{r-1} {\gamma_{i+1}-\gamma_i\over r} \gamma^{(i)}
$$
with
$$
\gamma^{(i)}\q:=\q \bigl(\, {\underbrace{i-r,...,i-r}_{i\times}},
\underbrace{i,...,i}_{(r-i)\times}\,\bigr),\q i=1,...,r-1.
$$
\subsection{Estimates for the weights of some special representations}
In the following, $\rho_{a,b,c}$ will stand for the induced
representation
of $\GL(r)$ on the vector space
$V_{a,b,c}:=\bigl((\C^r)^{\otimes a}\otimes (\bigwedge^r
\C^r)^{\otimes -b}\bigr)^{\oplus c}$
where $a,b\in \Z_{\ge 0}$, $c\in\Z_{>0}$.
Then, $\Pe(V_{a,b,c})=\Pe(V_{a,0,c})$ and $V_{a,b,c}\cong V_{a,0,c}$
as $\SL(r)$-modules.
\par
Let $\ul{w}=(w_1,...,w_r)$ be a basis for $\C^r$ and $\ul{\gamma}=
\sum_{i=1}^{r-1}\alpha_i \gamma^{(i)}$, $\alpha_i\in \Q_{\ge 0}$,
an integral weight vector.
Let $I^a$ be the set of all $a$-tuples $\ul{\iota}=(\iota_1,...,\iota_a)$
with $\iota_j\in\{\, 1,...,r\,\}$, $j=1,...,a$.
For $\ul{\iota}\in I^a$ and $k\in \{\, 1,...,c\,\}$,
we define $w_{\ul{\iota}}:=w_{\iota_1}\otimes\cdots
\otimes w_{\iota_a}$, and $w_{\ul{\iota}}^k:=(0,...,0,w_{\ul{\iota}},0,...,0)$,
$w_{\ul{\iota}}$ occupying the $k$-th entry.
The elements $w^k_{\ul{\iota}}$ with $\ul{\iota}\in I^a$ and
$k\in\{\, 1,...,c\,\}$ form a basis for $V_{a,0,c}$.
We let ${w^k_{\ul{\iota}}}^\vee$, $\ul{\iota}\in I^a$,
$k\in\{\, 1,...,c\,\}$ be the dual basis of $V^\vee_{a,0,c}$.
Now, let $[l]\in \Pe(V_{a,0,c})$
where $l=\sum a^k_{\ul{\iota}}{w^k_{\ul{\iota}}}^\vee$.
Then, there exist  $k_0$ and $\ul{\iota}_0$ with
$a^{k_0}_{\ul{\iota}_0}\neq 0$ and
$$
\mu_{\rho_{a,b,c}}\bigl(\la(\ul{w},\ul{\gamma}), [l]\bigr)
=
\mu_{\rho_{a,0,c}}\bigl(\la(\ul{w},\ul{\gamma}),
[l]\bigr)
=
\mu_{\rho_{a,0,c}}\bigl(\la(\ul{w},\ul{\gamma}),
[{w^{k_0}_{\ul{\iota}_0}}^\vee
]\bigr),
$$
and for any other $k$ and $\ul{\iota}$ with $a^k_{\ul{\iota}}\neq 0$
$$
\mu_{\rho_{a,b,c}}\bigl(\la(\ul{w},\ul{\gamma}),[l]\bigr)\q
\ge\q
\mu_{\rho_{a,0,c}}\bigl(\la(\ul{w},\ul{\gamma}),
[{w^{k}_{\ul{\iota}}}^\vee
]\bigr).
$$
We also find that for $i\in\{\, 1,...,r-1\,\}$
$$
\mu_{\rho_{a,0,c}}\bigl(\la(\ul{w},\gamma^{(i)}),
[{w^{k_0}_{\ul{\iota}_0}}^\vee
]\bigr)= \nu\cdot r- a\cdot i,\
\nu=\#\bigl\{\, \iota_j\le i\,|\, \ul{\iota}_0=
(\iota_1,...,\iota_a),\ j=1,...,a\,
\bigr\}.
$$
One concludes
\begin{Lem}
\label{Add}
{\rm i)}
For every basis $\ul{w}=(w_1,...,w_r)$ of $\C^r$, every integral
weight vector $\ul{\gamma}=
\sum_{i=1}^{r-1}\alpha_i \gamma^{(i)}$, $\alpha_i\in \Q_{\ge 0}$,
and every point $[l]\in\Pe(V_{a,b,c})$
$$
\left(\sum_{i=1}^{r-1}\alpha_i\right)a(r-1)\q\ge\q
\mu_{\rho_{a,b,c}}\bigl(\la(\ul{w},\ul{\gamma}),[l]\bigr)\q\ge\q
-\left(\sum_{i=1}^{r-1}\alpha_i\right)a(r-1).
$$
\par
{\rm ii)}
For every basis $\ul{w}=(w_1,...,w_r)$ of $\C^r$, every two
integral
weight vectors $\ul{\gamma}_1=
\sum_{i=1}^{r-1}\alpha_i \gamma^{(i)}$, $\alpha_i\in \Q_{\ge 0}$,
$\ul{\gamma}_2=
\sum_{i=1}^{r-1}\beta_i \gamma^{(i)}$, $\beta_i\in \Q_{\ge 0}$,
and every point $[l]\in\Pe(V_{a,b,c})$
$$
\mu_{\rho_{a,b,c}}\bigl(\la(\ul{w},\ul{\gamma}_1+\ul{\gamma}_2),[l]\bigr)
\q\ge\q
\mu_{\rho_{a,b,c}}\bigl(\la(\ul{w},\ul{\gamma}_1),[l]\bigr)
-\left(\sum_{i=1}^{r-1}\beta_i\right)a(r-1).
$$
\end{Lem}
\section{Decorated vector bundles}
\subsection{The moduli functors}
In this section, we will introduce the vector bundle problems we would
like to treat. The main topic will be the definition of the
semistability concept.
Having done this, we describe the relevant moduli functors to be studied
throughout the rest of this chapter.
\subsubsection{Semistable objects}
The input data for our construction are:
\begin{itemize}
\item a positive integer $r$,
\item an action of the general linear group $\GL(r)$ on the projective
      manifold $F$, such that the centre $\C^*\subset\GL(r)$ acts trivially.
\end{itemize}
The objects we want to classify are pairs $(E,\sigma)$
where
\begin{itemize}
\item $E$ is a vector bundle of rank $r$, and
\item $\sigma\colon X\lra {\frak F}(E):={\frak P}(E)\times^{\GL(r)} F$
is a section.
\end{itemize}
Here, ${\frak P}(E)$ is the principal $\GL(r)$-bundle associated with $E$.
Uninspired as we are, we call $(E,\sigma)$ an \it $F$-pair\rm.
Two $F$-pairs $(E^1,\sigma^1)$ and $(E^2,\sigma^2)$ are called
\it equivalent\rm, if there exists an isomorphism $\psi\colon E^1\lra E^2$
such that $\sigma^1=\sigma^2\circ \widehat{\psi}$,
$\widehat{\psi}\colon {\frak F}(E^1)\lra {\frak F}(E^2)$ being the
induced isomorphism.
\par
It will be our task to formulate a suitable semistability concept
for these objects and to perform a construction of the moduli spaces.
Let $E$ be a vector bundle over $X$. A \it weighted filtration of $E$ \rm
is a pair $(E^\bullet, \ul{\alpha})$ consisting of a filtration
$E^\bullet: 0\subset E_1\subset\cdots\subset E_s\subset E$ of
$E$ by non-trivial proper subbundles and a vector
$\ul{\alpha}=(\alpha_1,...,\alpha_s)$ of positive rational numbers.
Given such a weighted filtration, we set
$$
M\bigl(E^\bullet,\ul{\alpha}\bigr)\q :=\q
\sum_{j=1}^s \alpha_{j}\bigl(\deg(E)\rk E_j-\deg(E_j)\rk E\bigr).
$$
Suppose we are also given
a linearization $\rho$ of the $\GL(r)$-action on $F$ in an ample line
bundle $L$.
Let $(E,\sigma)$ be as above and $(E^\bullet, \ul{\alpha})$ be a weighted
filtration of $E$.
We define
$\mu_\rho(E^\bullet,\ul{\alpha};\sigma)$ as
follows:
Let $\ul{w}=(w_1,...,w_r)$ be an arbitrary basis of $W:=\C^r$.
For every $i\in\{\,1,...,r-1\,\}$, we set $W_{\ul{w}}^{(i)}
:=\langle\, w_1,...,w_i\,
\rangle$. Define $i_j:=\rk E_j$, $j=1,...,s$. This provides
a flag
$$
W^\bullet: 0\subset W_{\ul{w}}^{(i_1)}
\subset\cdots\subset W_{\ul{w}}^{(i_s)}
\subset W
$$
and thus a parabolic subgroup $P\subset \SL(r)$, namely the stabilizer
of the flag $W^\bullet$. Finally, set $\ul{\gamma}=\sum_{j=1}^s
\alpha_j \gamma^{(i_j)}$.
Next, let $U$ be an open subset of $X$ over which there is an
isomorphism $\psi\colon E_{|U}\lra W\otimes\Oh _U$
with $\psi(E^\bullet_{|U})=W^\bullet\otimes\Oh _U$.
Then, $\psi$ gives us an isomorphism ${\frak F}(E_{|U})\lra U\times F$
and $\sigma$ a morphism $\widetilde{\sigma}\colon U\lra U\times F\lra F$.
If $\ul{\gamma}$ is a vector of integers, we set
$\mu_{\rho}(E^\bullet,\ul{\alpha};\sigma):=
\mu_{\rho}(\la(\ul{w},\ul{\gamma}),\widetilde{\sigma})$
as in Lemma~\ref{Def}. Otherwise, we choose $k>0$ such that
$k\cdot\ul{\gamma}$ is a vector of integers and set
$\mu_{\rho}(E^\bullet,\ul{\alpha};\sigma):=
(1/k)\mu_{\rho}(\la(\ul{w}, k\ul{\gamma}),\widetilde{\sigma})$.
Since for an integral weight vector $\ul{\gamma}^\p$ and
a positive integer $k^\p$ one has $\mu_\rho(\la(\ul{w},\allowbreak k^\p\ul{\gamma}^\p),\sigma)=
k^\p\mu_\rho(\la(\ul{w},\ul{\gamma}^\p),\sigma)$,
this is well-defined. Note that the weight vector $\ul{\gamma}$ is canonically defined by $(E^\bullet,\ul{\alpha})$,
but that we have to verify that definition does not depend on the basis $\ul{w}$
and the trivialization
$\psi$.
First, let $\ul{w}^\p=(w_1^\p,...,w_r^\p)$ be a different basis.
Let $g\in \GL(r)$ be the element which maps $w_i$ to $w^\p_i$, $i=1,...,r$,
and set $\psi^\p:=(g\otimes\id_{\Oh _U})\circ\psi$.
This defines the morphism $\widetilde{\sigma}^\p\colon U\lra F$.
Then, $\la(\ul{w}^\p,\ul{\gamma})=g\cdot\la(\ul{w},\ul{\gamma})\cdot g^{-1}$
and $\widetilde{\sigma}^\p(x)=g\cdot\widetilde{\sigma}(x)$
for every $x\in U$.
Since $\mu_{\rho}(\la(\ul{w}^\p,\ul{\gamma}), \widetilde{\sigma}^\p(x))=
\mu_{\rho}(g\cdot\la(\ul{w},\ul{\gamma})\cdot g^{-1},
g\cdot\widetilde{\sigma}(x))
=\mu_{\rho}(\la(\ul{w},\ul{\gamma}), \widetilde{\sigma}(x))$,
we may fix the
basis $\ul{w}$.
Any other trivialization $\widetilde{\psi}$ defined w.r.t.\ $\ul{w}$
differs from $\psi$ by a map $U\lra P$.
Now, for every $g\in P$ and every point $x\in U$,
$\mu_{\rho}(\la,\widetilde{\sigma}(x))=\mu_{\rho}(g\la g^{-1}, g\cdot
\widetilde{\sigma}(x))
=\mu_{\rho}(\la, g\cdot \widetilde{\sigma}(x))$.
The last equality results from \cite{GIT}, Prop.~2.7, p.~57.
This shows our assertion.
To conclude, Remark~\ref{SemCont} shows that the definition is also
independent
of the choice of the open subset $U$.
\par
Fix also a number $\delta\in \Q_{>0}$.
With these conventions, we call $(E,\sigma)$ \it
$\delta$-$\rho$-(semi)stable\rm,
if for every weighted filtration $(E^\bullet,\ul{\alpha})$ of $E$
$$
M\bigl(E^\bullet,\ul{\alpha}\bigr)\q+\q \delta\cdot
\mu_{\rho}\bigl(E^\bullet, \ul{\alpha};\sigma\bigr)\q (\ge) \q 0.
$$
\par
Next, we remark that we should naturally fix
the degree of $E$. Then, the topological
fibre space $\pi\colon {\frak F}^{d,r}\lra X$ underlying ${\frak F}(E)$
will be independent of
$E$, so that it makes sense to fix the homology class $[\sigma(X)]\in
H_2({\frak F}^{d,r},\Z)$.
Given $d\in\Z$, $r\in\Z_{>0}$, and $h\in H_2({\frak F}^{d,r},\Z)$, we say
that $(E,\sigma)$ is \it of type $(d,r,h)$\rm, if $E$ is a vector
bundle of degree $d$ and rank $r$, and $[\sigma(X)]=h$.
Before we define the moduli functor, we enlarge our scope.
\par
For a given linearization of the $\GL(r)$-action on $F$ in the
line bundle $L$, we can choose a positive
integer $k$ such that $L^{\otimes k}$ is very ample.
Therefore, we obtain a $\GL(r)$-equivariant
embedding $F\hookrightarrow \Pe(V)$, $V:=H^0(F,L^{\otimes k})$.
Note that $\C^*$ acts trivially on $\Pe(V)$.
Therefore, we formulate the following classification problem:
The input now consists of
\begin{itemize}
\item a positive integer $r$, a finite dimensional vector space $V$, and
\item a representation $\rho\colon \GL(r)\lra \GL(V)$ whose restriction
      to the centre $\C^*$ is of the form $z\lma z^\alpha\cdot \id_V$
      for some integer $\alpha$,
\end{itemize}
and the objects we want to classify are pairs $(E,\sigma)$
where
\begin{itemize}
\item $E$ is a vector bundle of rank $r$, and
\item $\sigma\colon X\lra \Pe(E_\rho)$ is a section. Here, $E_\rho$ is
      the vector bundle of rank $\dim V$ associated to $E$ via the
      representation $\rho$.
\end{itemize}
The equivalence relation is the same as before.
Now, giving a section $\sigma\colon X\lra \Pe(E_\rho)$ is the same
as giving a line bundle $M$ on $X$ and a surjection $\tau\colon E_\rho\lra M$.
Remember that $(M,\tau)$ and $(M^\p,\tau^\p)$ give the same section
if and only if there exists an isomorphism $M\lra M^\p$ which carries
$\tau$ into $\tau^\p$. Moreover, fixing the homology class $[\sigma(X)]$
amounts to the same as fixing the degree of $M$.
Since the condition that $\tau$ be surjective will be an open condition
in a suitable parameter space, we formulate
the following classification problem:
The input data are
\begin{itemize}
\item a tuple $(d,r,m)$ called the \it type\rm, where $d$, $r$, and $m$
      are integers, $r>0$,
\item a representation $\rho\colon \GL(r)\lra \GL(V)$,
\end{itemize}
and the objects to classify are triples $(E,M,\tau)$
where
\begin{itemize}
\item $E$ is a vector bundle of rank $r$ and degree $d$,
\item $M$ is a line
      bundle of degree $m$, and
\item $\tau\colon E_\rho\lra M$ is a non-zero
      homomorphism.
\end{itemize}
Then, $(E,M,\tau)$ is called a \it $\rho$-pair of type $(d,r,m)$\rm,
and $(E^1,M^1,\tau^1)$ and $(E^2,M^2,\tau^2)$
are said to be \it equivalent\rm,
if there exist isomorphisms $\psi\colon E^1\lra E^2$ and $\chi\colon
M^1\lra M^2$ with $\tau^1=\chi^{-1}\circ \tau^2\circ \psi_\rho$,
where $\psi_\rho\colon E^1_\rho\lra E_\rho^2$ is the induced isomorphism.
Let $(E,M,\tau)$ be a $\rho$-pair of type $(d,r,m)$. A \it weak automorphism \rm
of $(E,M,\tau)$ is the class $[\psi]\in \Pe(\End(E))$ of an automorphism
$\psi\colon E\lra E$ with $\tau=\tau\circ \psi_\rho$. We call $(E,M,\tau)$
\it simple\rm, if there are only finitely many weak automorphisms.
\begin{Rem}
i)
A representation $\rho\colon \GL(r)\lra \GL(V)$ of the general linear group
with $\rho(z\cdot E_n)=z^\alpha\cdot \id_V$ is called \it homogeneous of
degree $\alpha$\rm. Every representation of $\GL(r)$ obviously splits into
a direct sum of homogeneous representations.
Some cases of inhomogeneous representations $\rho$ can be treated within our
framework. Indeed, if $\rho$ is a representation, such that its homogeneous components $\rho_1,...,\rho_n$
have \sl positive \rm degrees $\alpha_1,...,\alpha_n$, let
$\kappa$ be a common multiple of the $\alpha_i$. Then, we pass to the homogeneous
representation
$$
\rho^\p\q:=\q  \bigoplus_{\nu_1\alpha_1+\cdots+\nu_n\alpha_n=\kappa} S^{\nu_1}\rho_1
\otimes\cdots\otimes S^{\nu_n}\rho_n.
$$
The solution of the moduli problem associated with $\rho^\p$ can be used to solve
the moduli problem associated with $\rho$. This trick was used in \cite{OST} and
will be recalled in the section on examples.
\par
ii)
The identification of $\tau$ and $\la\cdot \tau$, or equivalently,
considering sections in $\Pe(E_\rho)$ rather than in $E_\rho$ seems a little
artificial. First of all, this identification is mandatory to get projective
moduli spaces. Second, for homogeneous representations of degree $\alpha\neq 0$,
this is naturally forced upon us. Third, if we are given a homogeneous representation
$\rho$
of degree zero and are interested in the moduli problem without the identification
of $\tau$ and $\la\tau$, we may pass to the representation $\rho^\p$, obtained
from $\rho$ by adding the trivial one dimensional representation.
Then, one gets from the solution of the moduli problem associated with $\rho^\p$
a compactification of the moduli problem associated with $\rho$.
This will be explained within the context of Hitchin pairs in the examples.
\end{Rem}
In order to define a functor, we first fix a Poincar\'e line
bundle ${\cal L}$ on $\mathop{\rm Jac}^m\times X$.
For every scheme $S$ and every morphism $\kappa\colon S\lra \Jac^m$,
we define $\EL[\kappa]:=(\kappa\times\id_X)^*\EL$.
Now, let $S$ be a scheme of finite type over $\C$.
Then, a \it family of $\rho$-pairs of type $(d,r,m)$ parameterized by $S$
\rm
is a tuple $(E_S,\kappa_S,{\frak N}_S,\tau_S)$ with
\begin{itemize}
\item $E_S$ a vector bundle of rank $r$ having degree $d$ on $\{s\}\times X$
      for all $s\in S$,
\item $\kappa_S\colon S\lra \Jac^m$ a morphism,
\item ${\frak N}_S$ a line bundle on $S$,
\item $\tau_S\colon E_{S,\rho}\lra \EL[\kappa_S]\otimes\pi_S^*{\frak N}_S$
      a homomorphism whose
      restriction to $\{s\}\times X$ is non-zero for every closed
      point $s\in S$.
\end{itemize}
Two families $(E^1_S,\kappa_S^1,{\frak N}^1_S,\tau_S^1)$
and $(E^2_S,\kappa_S^2,{\frak N}^2_S,\tau_S^2)$
are called \it equivalent\rm, if $\kappa_S^1=\kappa_S^2=:\kappa_S$
and there exist
isomorphisms $\psi_S\colon E^1_S\lra E^2_S$ and $\chi_S\colon{\frak N}_S^1
\lra {\frak N}_S^2$ with
$$
\tau^1_S= ({\id}_{\EL[\kappa_S]}\otimes \pi_S^*\chi_S)^{-1}\circ
\tau^2_S\circ \psi_{S,\rho}.
$$
To define the semistability concept for $\rho$-pairs, observe that
for given $(E,M,\tau)$, the homomorphism $\tau\colon E\lra M$ will
be generically surjective, therefore we get a rational section
$\sigma^\p\colon X\dasharrow \Pe(E_\rho)$ which can, of course, be
prolonged to a section $\sigma\colon X\lra \Pe(E_\rho)$, so that we can
define for every weighted filtration $(E^\bullet,\ul{\alpha})$
of $E$
$$
\mu_\rho\bigl(E^\bullet, \ul{\alpha};\tau\bigr)\q :=\q
\mu_\rho\bigl(E^\bullet,\ul{\alpha};\sigma\bigr).
$$
We will occasionally use the following short hand notation:
If $E^\p$ is a non-zero, proper subbundle of $E$, we
set
$$
\mu_{\rho}(E^\p,\tau)\q:=\q\mu_{\rho}(0\subset E^\p\subset E, (1);\tau).
$$
\par
Now, for fixed $\delta\in\Q_{>0}$, call a $\rho$-pair $(E,M,\tau)$
\it $\delta$-(semi)stable\rm,
if for every weighted filtration $(E^\bullet,\ul{\alpha})$
$$
M\bigl(E^\bullet,\ul{\alpha}\bigr)\q+\q\delta\cdot
\mu_{\rho}\bigl(E^\bullet,\ul{\alpha};\tau\bigr)\q (\ge) \q 0.
$$
\begin{Rem}
For the $F$-pairs, one can formulate the semistability concept in a more
intrinsic way. For this, one just has to choose a linearization $\rho$
of the given action in an ample $\Q$-line bundle. Then,
$\mu_\rho(E^\bullet, \ul{\alpha}; \Phi)$ can still be defined, and an
$F$-pair $(E,\Phi)$ will be called \it $\rho$-(semi)stable\rm, if
$$
M\bigl(E^\bullet,\ul{\alpha}\bigr)+\mu_\rho\bigl(E^\bullet,\ul{\alpha};
\Phi\bigr)\q (\ge)\q 0.
$$
In gauge theory, one would say that the notion of semistability depends
only on the metric chosen on the fibre $F$. If $\rho$ is a linearization
in an ample line bundle $L$ and $\delta\in\Q$, we can pass to the induced
linearization ``$\rho^{\otimes \delta}$" in the $\Q$-line bundle $\delta L$
to recover $\delta$-$\rho$-semistability.
For the moduli problems associated with a representation $\rho$,
the formulation with the parameter $\delta$ seems more appropriate and
practical and, since we treat $F$-pairs only as special cases of $\rho$-pairs,
we have decided for the given definition of $\delta$-$\rho$-semistability.
\end{Rem}
We define the functors
$$
\begin{array}{rrcl}
\ul{\hbox{M}}(\rho)_{d/r/m}^{\delta-(s)s}\colon & \ul{\hbox{Sch}}_\C &\lra & \ul{\hbox{Set}}\\
& S & \lma & \left\{{\hbox{Equivalence classes of families of $\delta$-(semi)stable}
\atop\hbox{$\rho$-pairs of type $(d,r,m)$ parameterized by $S$}}\right\}.
\end{array}
$$
\begin{Rem}
\label{Indep}
The definition of the moduli functor involves the choice of
the Poin\-car\'e sheaf $\EL$.
Nevertheless, the above moduli functor
is independent of that choice. Indeed, choosing another Poincar\'e line
bundle $\EL^\p$ on $\Jac^m\times X$, there is a line bundle
${\frak N}_{\Jac^m}$ on $\Jac^m$ with $\EL\cong\EL^\p\otimes\pi_{\Jac^m}^*
{\frak N}_{\Jac^m}$. Therefore, assigning to a family
$(E_S,\kappa_S,{\frak N}_S,\tau_S)$ defined via $\EL$ the
family $(E_S,\kappa_S,{\frak N}_S\otimes\kappa_S^*{\frak N}_{\Jac^m},\tau_S)$
defined via $\EL^\p$ identifies the functor which is defined w.r.t.\ $\EL$
with the one defined w.r.t.\ $\EL^\p$.
\end{Rem}
We also define the open subfunctors
$\ul{\hbox{M}}(\rho)_{d/r/m/\mathop{\rm surj}}^{\delta-(s)s}$
of equivalence classes of families $(E_S, \kappa_S, {\frak N}_S,\tau_S)$
where $\tau_{S|\{s\}\times X}$ is surjective for all $s\in S$.
\par
Next,
let $(E,M,\tau)$ be a $\rho$-pair where $\tau$ is surjective, and
let ${\Bbb P}^{d,r}$ be the oriented topological projective
bundle underlying $\Pe(E_\rho)$. This is independent of $E$,
and as explained before, the degree of $M$ determines
the cohomology class $h_m:=[\sigma(X)]\in H_2({\Bbb P}^{d,r},\Z)$
where $\sigma$ is the section associated with $\tau$.
Set $h:=h_m\cap [{\frak F}^{d,r}]\in H_2({\frak F}^{d,r},\Z)$.
We can now define $\ul{\hbox{M}}(F,\rho)_{d/r/h}^{\delta-(s)s}$ as the
closed subfunctor of $\ul{\hbox{M}}(\rho)_{d/r/m/\mathop{\rm surj}}^{\delta-(s)s}$
of equivalence classes of families $(E_S,\kappa_S, {\frak N}_S,\tau_S)$
for which the section $S\times X\lra \Pe(E_{S,\rho})$ factorizes
over ${\frak P}(E_{S})\times^{\GL(r)} F$.
\subsubsection{Polystable pairs}
Fix a basis $\ul{w}=(w_1,...,w_r)$ of $\C^r$.
Let $(E,M,\tau)$ be a $\delta$-semistable $\rho$-pair of type $(d,r,m)$.
We call $(E,M,\tau)$ \it $\delta$-polystable\rm,
if for every weighted filtration $(E^\bullet,\ul{\alpha})$,
$E^\bullet: 0=:E_0\subset E_1\subset\cdots \subset  E_s\subset E_{s+1}:=E$,
with
$$
M\bigl(E^\bullet,\ul{\alpha}\bigr)\q+\q \delta\cdot
\mu_{\rho}\bigl(E^\bullet,\ul{\alpha};\tau\bigr)\q = \q 0
$$
the following holds true
\begin{itemize}
\item $E\cong \bigoplus_{j=1}^{s+1} E_j/E_{j-1}$
\item the $\rho$-pair $(E,M,\tau)$ is equivalent to the $\rho$-pair $(E,M,\tau|_{E^\gamma})$,
	  $\gamma:=-\mu_\rho(E^\bullet,\ul{\alpha};\tau)$.
\end{itemize}
Here, one uses the fact that giving an isomorphism $E\lra \bigoplus^{s+1}_{j=1}
E_j/E_{j-1}$ is the same as giving a cocycle for $E$ in the group
$$
Z\bigl(\la(\ul{w},\ul{\gamma})\bigr):=\bigl\{\,g\in \GL(r)\,|\, g\cdot \la(\ul{w},\ul{\gamma})(z)=
\la(\ul{w},\ul{\gamma})(z)
\cdot g\ \forall z\in\C^*\,\bigr\}.
$$ 
It follows that $E_\rho\cong E^{\gamma_1}\oplus\cdots\oplus E^{\gamma_{t}}$ where $\gamma_1$,...,$\gamma_t$ 
are the weights of
$\la(\ul{w},\ul{\gamma})$ on $V$ and $E^{\gamma_i}$ is the ``eigenbundle" for the weight $\gamma_i$, $i=1,...,t$.
As before, $W^\bullet: 0\subset W_{\ul{w}}^{(\rk E_1)}\subset\cdots
\subset W_{\ul{w}}^{(\rk E_s)}\subset W$ and $\ul{\gamma}:=\sum_{j=1}^s
\alpha_j\gamma^{(\rk E_j)}$. The stated condition is again independent of
the involved choices.
\begin{Rem}
\label{st=pst+s}
i)
If $(E,M,\tau)$ is $\delta$-stable, the stated condition is void, so that
$(E,M,\allowbreak\tau)$ is also $\delta$-polystable.
\par
ii) It will follow from our GIT construction that $(E,M,\tau)$ is
$\delta$-stable, if and only if it is $\delta$-polystable and has
only finitely many weak automorphisms.
\par
iii) For the description of S-equivalence in the case of $\rho=\rho_{a,b,c}$
for some $a,b,c\in\Z_{\ge 0}$, the reader may consult \cite{GS2}.
\end{Rem}
\subsection{The main result}
\begin{Thm}
\label{Main}
{\rm i)} There exist a projective scheme
        ${\cal M}(\rho)^{\delta-ss}_{d/r/m}$
        and an open subscheme
        ${\cal M}(\rho)^{\delta-s}_{d/r/m}\subset
        {\cal M}(\rho)^{\delta-ss}_{d/r/m}$ together with
        natural transformations
        $$\vartheta^{(s)s}\colon\ul{\hbox{\rm M}}(\rho)^{\delta-(s)s}_{d/r/m}
        \lra h_{{\cal M}(\rho)^{\delta-(s)s}_{d/r/m}}
        $$ with the following properties:
\begin{enumerate}
\item[\rm 1.] For every scheme ${\cal N}$ and
       every natural transformation $\vartheta^\p\colon
       \ul{\hbox{\rm M}}(\rho)^{\delta-ss}_{d/r/m}\allowbreak
       \lra\allowbreak h_{\cal N}$, there
       exists a unique morphism $\phi\colon
       {\cal M}(\rho)^{\delta-ss}_{d/r/m}\lra {\cal N}$ with
       $\vartheta^\p=h(\phi)\circ\vartheta^{ss}$.
\item[\rm 2.] ${\cal M}(\rho)^{\delta-s}_{d/r/m}$ is a coarse moduli space
      for the functor $\ul{\hbox{\rm M}}(\rho)^{\delta-s}_{d/r/m}$.
\item[\rm 3.] $\vartheta^{ss}(\Spec\C)$
      induces a bijection between the set of equivalence
      classes of $\delta$-polystable $\rho$-pairs of type $(d,r,m)$
      and the set of closed points of ${\cal M}(\rho)^{\delta-ss}_{d/r/m}$.
\end{enumerate}
\par
{\rm ii)} There exist a locally closed subscheme
          ${\cal M}(F,\rho)^{\delta-s}_{d/r/h}$
          of ${\cal M}(\rho)^{\delta-s}_{d/r/m}$
          and a natural transformation
          $$
          \vartheta_F\colon\ul{\hbox{\rm M}}(F,\rho)^{\delta-s}_{d/r/h}
          \lra h_{{\cal M}(F,\rho)^{\delta-s}_{d/r/h}}
          $$
          which turns ${\cal M}(F,\rho)^{\delta-s}_{d/r/h}$
          into the coarse moduli space for
          $\ul{\hbox{\rm M}}(F,\rho)^{\delta-s}_{d/r/h}$.
\end{Thm}
\subsection{The proof of the main result}
Given any homogeneous representation $\rho\colon \GL(r)\lra \GL(V)$,
we have seen in Section~\ref{Rep} that we can find
integers $a,b\ge 0$ and $c>0$, such that $\rho$ is a direct summand
of the representation $\rho_{a,b,c}$. Write $\rho_{a,b,c}=\rho\oplus
\ol{\rho}$. For every vector bundle $E$ of rank $r$, we find
$E_{\rho_{a,b,c}}\cong E_\rho\oplus E_{\ol{\rho}}$.
Every $\rho$-pair $(E,M,\tau)$ can therefore also be viewed
as a $\rho_{a,b,c}$-pair. Since $\mu_\rho(E^\bullet,\ul{\alpha};\tau)=
\mu_{\rho_{a,b,c}}(E^\bullet,\ul{\alpha}; \tau)$ for every weighted filtration
$(E^\bullet,\ul{\alpha})$, the triple $(E,M,\tau)$
is $\delta$-(semi)stable
as $\rho$-pair, if and only if it is $\delta$-(semi)stable as
$\rho_{a,b,c}$-pair.
More precisely, we can recover
$\ul{\hbox{M}}(\rho)^{\delta-(s)s}_{d/r/m}$ as closed
subfunctor of $\ul{\hbox{M}}(\rho_{a,b,c})^{\delta-(s)s}_{d/r/m}$.
Indeed, for every scheme of finite type over $\C$,
\begin{eqnarray*}
\ul{\hbox{M}}(\rho)^{\delta-(s)s}_{d/r/m}(S)
&=&\Bigl\{\, [E_S,\kappa_S,{\frak N}_S,
\tau_S]\in \ul{\hbox{M}}(\rho_{a,b,c})^{\delta-(s)s}_{d/r/m}(S)\,|\,\\
&&\qquad\tau_S\colon E_{S,\rho_{a,b,c}}\lra
\EL[\kappa_S]\otimes\pi_S^*{\frak N}_S
\hbox{ vanishes on } E_{S,\ol{\rho}}\,\Bigr\}.
\end{eqnarray*}
Therefore, we will assume from now on that $\rho=\rho_{a,b,c}$
for some $a,b,c$.
\subsubsection{Boundedness}
\begin{Thm}
\label{bound2}
There is a non-negative
constant $C_1$, depending only on $r$, $a$, and $\delta$, such that
for every $\delta$-semistable $\rho_{a,b,c}$-pair $(E,M,\tau)$
of type $(d,r,m)$ and every non-trivial proper subbundle $E^\p$ of $E$
$$
\mu(E^\p)\q \le\q {d\over r}+C_1.
$$
\end{Thm}
\begin{proof}
Let $0 \subsetneq E^\p\subsetneq E$ be any
subbundle. By Lemma~\ref{Add}~i), $\mu_{\rho_{a,b,c}}
(E^\p,\tau)\le a(r-1)$,
so that $\delta$-semistability
gives
\begin{eqnarray*}
&& d\rk E^\p-\deg(E^\p)r+\delta\cdot a\cdot (r-1)\\
&\ge&d\rk E^\p -\deg(E^\p)r+\delta\cdot \mu_{\rho_{a,b,c}}\bigl(E^\p,
\tau\bigr)\q\ge\q 0,
\end{eqnarray*}
i.e.,
$$
\mu(E^\p)\q\le\q {d\over r}+{\delta\cdot a\cdot(r-1)\over r\cdot\rk E^\p}
\q\le\q  {d\over r}+{\delta\cdot a\cdot(r-1)\over r},
$$
so that the theorem holds for $C_1:=\delta\cdot a\cdot (r-1)/r$.

\end{proof}
\subsubsection{Construction of the parameter space}
Recall that,
for a scheme $S$ of finite type over $\C$, a family of $\rho_{a,b,c}$-pairs
parameterized by $S$ is a quadruple $(E_S,\kappa_S,{\frak N}_S,\tau_S)$
where $E_S$ is a vector bundle of rank $r$ on $S\times X$
with $\deg(E_{S|\{s\}\times X})=d$ for all $s\in S$,
$\kappa_S\colon S\lra \Jac^m$ is a morphism, ${\frak N}_S$ is a line
bundle on $S$, and
$\tau_S\colon {E_S^{\otimes a}}^{\oplus c}\lra
\det(E_S)^{\otimes b}\otimes\EL[\kappa_S]\otimes\pi_S^*{\frak N}_S$
is a homomorphism which is non zero on every fibre $\{s\}\times X$.
\par
Pick a point $x_0\in X$, and write $\Oh _X(1)$ for $\Oh _X(x_0)$.
According to~\ref{bound2}, we can choose
an integer $n_0$, such that for every $n\ge n_0$
and every $\delta$-semistable $\rho_{a,b,c}$-pair $(E,M,\tau)$
of type $(d,r,m)$
\begin{itemize}
\item $H^1(E(n))=0$ and $E(n)$ is globally generated,
\item $H^1(\det(E)(rn))=0$ and $\det(E)(rn)$ is globally
      generated,
\item $H^1(\det(E)^{\otimes b}\otimes M\otimes\Oh _X(na))=0$ and
      $\det(E)^{\otimes b}\otimes M\otimes\Oh _X(na)$ is globally generated.
\end{itemize}
Choose some $n\ge n_0$ and set $p:=d+rn+r(1-g)$.
Let $U$ be a complex vector space of dimension $p$.
We define ${\frak Q}^0$ as the quasi-projective scheme parameterizing
equivalence classes of quotients
$q\colon U\otimes \Oh _X(-n)\lra E$
where $E$ is a vector bundle of rank $r$ and degree $d$ on $X$
and $H^0(q(n))$ is an isomorphism.
Then there exists a universal quotient
$$
q_{{\frak Q}^0}\colon U\otimes\pi_X^*\Oh _X(-n)\lra E_{{\frak Q}^0}
$$
on ${\frak Q}^0\times X$.
Let
$$
q_{{\frak Q}^0\times \Jac^m}\colon U\otimes\pi_X^*\Oh _X(-n)\lra
E_{{\frak Q}^0\times\Jac^m}
$$
be the pullback of $q_{{\frak Q}^0}$ to ${\frak Q}^0\times\Jac^m\times X$.
Set $U_{a,c}:={U^{\otimes a}}^{\oplus c}$.
By our assumption, the sheaf
$$
\ul{\Hom}\Bigl(U_{a,c}\otimes\Oh _{{\frak Q}^0\times
\Jac^m},\pi_{{\frak Q}^0\times\Jac^m *}
\bigl(\det(E_{{\frak Q}^0\times\Jac^m})^{\otimes b}\otimes\EL[\pi_{\Jac^m}]
\otimes \pi_X^*\Oh _X(na)\bigr)\Bigr)
$$
is locally free, call it ${\cal H}$, and set ${\frak H}:=\Pe({\cal H}^\vee)$.
We let
$$
q_{\frak H}\colon U\otimes\pi_X^*\Oh _X(-n)\lra
E_{\frak H}
$$
be the pullback of $q_{{\frak Q}^0\times \Jac^m}$ to ${\frak H}\times X$.
Now, on ${\frak H}\times X$, there is the tautological homomorphism
$$
s_{\frak H}\colon
U_{a,c}\otimes\Oh _{\frak H}\lra \det(E_{\frak H})^{\otimes b}\otimes
\EL[\kappa_{\frak H}]
\otimes\pi_X^*\Oh _X(na)\otimes\pi_{\frak H}^*\Oh _{\frak H}(1).
$$
Here, $\kappa_{\frak H}
\colon {\frak H}\lra {\frak Q}^0\times \Jac^m\lra \Jac^m$
is the natural morphism.
Let ${\frak T}$ be the closed subscheme defined by the condition that
$s_{\frak H}\otimes
\pi_X^*\id_{\Oh _X(-na)}$ vanish on
$$
\ker\bigl(
U_{a,c}\otimes\pi_X^*\Oh _X(-na)\lra {{E}_{\frak H}^{\otimes a}}^{\oplus c}
\bigr).
$$
Let
$$
q_{\frak T}\colon U\otimes\pi_X^*\Oh _X(-n)\lra E_{\frak T}
$$
be the restriction of $q_{\frak H}$ to ${\frak T}\times X$.
By definition, there is a universal homomorphism
$$
\tau_{\frak T}\colon {{E}_{\frak T}^{\otimes a}}^{\oplus c}
\lra \det(E_{\frak T})^{\otimes b}\otimes\EL[\kappa_{\frak T}]
\otimes \pi_{\frak T}^* {\frak N}_{\frak T}.
$$
Here, ${\frak N}_{\frak T}$ and $\kappa_{\frak T}$ are the restrictions
of $\Oh _{\frak H}(1)$ and $\kappa_{\frak H}$ to ${\frak T}$.
Note, that the parameter space ${\frak T}$ is equipped
with a universal family $(E_{\frak T}, \kappa_{\frak T}, {\frak N}_{\frak T},
\tau_{\frak T})$.
\begin{Rem}
\label{rep}
Let $S$ be a scheme of finite type over $\C$.
Call a tuple $(q_S, \kappa_S, {\frak N}_S, \tau_S)$
where
\begin{itemize}
\item $q_S\colon U\otimes \pi_X^*\Oh _X(-n)\lra E_S$ is a family
      of quotients, such that its restriction to $\{s\}\times X$
      lies in ${\frak Q}^0$ for every $s\in S$,
\item $\kappa_S\colon S\lra \Jac^m$ is a morphism,
\item ${\frak N}_S$ is a line bundle on $S$, and
\item $\tau_S\colon {E_{S}^{\otimes a}}^{\oplus c}
      \lra \det(E_S)^{\otimes b}\otimes
      \EL[\kappa_S]\otimes \pi_S^*{\frak N}_S$ is a homomorphism
      which is non trivial on all fibres $\{s\}\times X$, $s\in S$,
\end{itemize}
a \it quotient family of $\rho_{a,b,c}$-pairs of type $(d,r,m)$
parameterized by $S$\rm. We say that the families $(q^1_S,
\kappa^1_S, {\frak N}^1_S,\allowbreak \tau^1_S)$ and $(q^2_S, \kappa^2_S,
{\frak N}^2_S, \tau^2_S)$ are \it equivalent\rm, if
$\kappa_S^1=\kappa_S^2=:\kappa_S$ and there are isomorphisms
$\psi_S\colon E^1_S\lra E_S^2$ and $\chi_S\colon {\frak
N}_S^1\lra {\frak N}_S^2$ with $q_S^2=\psi_S\circ q_S^1$ and
$$
\tau_S^1=
\bigl(
{\id}_{\EL[\kappa_S]}\otimes\pi_S^*(\chi_S)
\bigr)^{-1}\circ \tau_S^2\circ
\Bigl({\psi_S^{\otimes a}}^{\oplus c}\Bigr).
$$
It can be easily inferred from the construction of ${\frak T}$
and the base change theorem that ${\frak T}$ represents the
functor which assigns to a scheme $S$ of finite type over $\C$
the set of equivalence classes of quotient families of
$\rho_{a,b,c}$-pairs of type $(d,r,m)$ parameterized by $S$.
\end{Rem}
\begin{Prop}[Local universal property]
\label{LUP}
Let $S$ be a scheme of finite type over $\C$, and $(E_S,\kappa_S,{\frak N}_S,\allowbreak
\tau_S)$
a family of $\delta$-semistable $\rho_{a,b,c}$-pairs parameterized by $S$.
Then, there exist an open covering $S_i$, $i\in I$, of $S$,
and morphisms $\beta_i\colon S_i\lra {\frak T}$, $i\in I$,
such that the restriction of the family $(E_S,\kappa_S,{\frak N}_S,\tau_S)$
to $S_i\times X$ is equivalent to the pullback of
$(E_{\frak T}, \kappa_{\frak T},\allowbreak
{\frak N}_{\frak T}, \tau_{\frak T})$ via $\beta_i\times\id_X$,
for all $i\in I$.
\end{Prop}
\begin{proof}
By our assumptions, the sheaf $\pi_{S*}(E_S\otimes\pi_X^*\Oh _X(n))$
is locally free of rank $p$.
Therefore, we can choose a covering $S_i$, $i\in I$, of $S$, such that
it is free over $S_i$ for all $i\in I$.
For each $i$, we can choose a trivialization
$$
U\otimes\Oh _{S_i}\q\cong\q \pi_{S*}\bigl(E_S\otimes\pi_X^*\Oh _X(n)_{|S_i}\bigr),
$$
so that we obtain a surjection
$$
q_{S_i}\colon U\otimes\pi_X^* \Oh _X(-n)\lra E_{S|S_i\times X}
$$
on $S_i\times X$.
Therefore, $(q_{S_i}, \kappa_{S|S_i}, {\frak N}_{S|S_i}, \tau_{S|S_i\times X})$
is a quotient family of $\rho_{a,b,c}$-pairs of type $(d,r,m)$
parameterized by $S_i$, and
we can conclude by Remark~\ref{rep}.

\end{proof}
\subsubsection{The group action}
Let $m\colon U\otimes\Oh _{\SL(U)}\lra U\otimes\Oh _{\SL(U)}$
be the universal automorphism over $\SL(U)$.
Let $({E}_{\SL(U)\times {\frak T}}, \kappa_{\SL(U)\times {\frak T}},
\allowbreak
{\frak N}_{\SL(U)\times {\frak T}}, \tau_{\SL(U)\times {\frak T}})$
be the pullback of the universal family on ${\frak T}\times X$
to $\SL(U)\times {\frak T}\times X$.
Define
$$
q_{\SL(U)\times{\frak T}}\colon
U\otimes\pi_X^*\Oh _X(-n)\stackrel{\pi_{\SL(U)}^*(m^{-1})\otimes
{\id}_{\pi_X^*\Oh _X(-n)}}{\lra}U\otimes\pi_X^*\Oh _X(-n)\lra
{E}_{\SL(U)\times {\frak T}}.
$$
Thus, $(q_{\SL(U)\times{\frak T}},\kappa_{\SL(U)\times {\frak T}},
{\frak N}_{\SL(U)\times {\frak T}}, \tau_{\SL(U)\times {\frak T}})$
is a quotient family of $\rho_{a,b,c}$-pairs
para\-meterized by $\SL(U)\times{\frak T}$, and hence, by~\ref{rep},
defines a morphism
$$
\Gamma\colon \SL(U)\times {\frak T}\lra {\frak T}.
$$
It is not hard to see that $\Gamma$ is indeed a group action. Note that this action descends
to a $\PGL(U)$-action!
\begin{Rem}
\label{lin}
By construction, the universal family
$(E_{\frak T}, \kappa_{\frak T}, {\frak N}_{\frak T}, \tau_{\frak T})$
comes with a linearization, i.e., with an isomorphism
$$
\bigl(\Gamma\times{\id}_X\bigr)^*
\bigl(E_{\frak T}, \kappa_{\frak T}, {\frak N}_{\frak T}, \tau_{\frak T}\bigr)
\lra
\bigl(\pi_{\frak T}\times{\id}_X\bigr)^*
\bigl(E_{\frak T}, \kappa_{\frak T}, {\frak N}_{\frak T},
\tau_{\frak T}\bigr).
$$
Therefore, elements of the $\PGL(U)$-stabilizer of a point $t\in {\frak T}$
correspond to weak automorphisms of the $\rho_{a,b,c}$-pair
$(E_t, M_t, \tau_t):= (E_{\frak T}, \kappa_{\frak T}, {\frak
N}_{\frak T}, \tau_{\frak T})_{|\{t\}\times X}$. In particular,
the $\SL(U)$-stabilizer of $t$ is finite if and only if $(E_t,M_t,\tau_t)$
has only finitely many weak automorphisms.
\end{Rem}
\begin{Prop}
\label{Glue}
Let $S$ be a scheme of finite type over $\C$ and $\beta_{1,2}\colon S
\lra {\frak T}$ two morphisms, such that the pullbacks of
$(E_{\frak T},\kappa_{\frak T}, {\frak N}_{\frak T}, \tau_{\frak T})$
via $\beta_1\times \id_X$ and $\beta_2\times \id_X$ are equivalent.
Then, there exist an \'etale covering $\eta\colon T\lra S$
and a morphism $\Xi\colon T\lra \SL(U)$, such that
the morphism $\beta_2\circ \eta\colon T\lra {\frak T}$ equals
the morphism
$$
T\stackrel{\Xi\times (\beta_1\circ\eta)}{\lra} \SL(U)\times {\frak T}
\stackrel{\Gamma}{\lra}{\frak T}.
$$
\end{Prop}
\begin{proof}
The two morphisms $\beta_1$ and $\beta_2$ provide us with quotient
families $(q_S^1,\kappa_S^1,\allowbreak {\frak N}_S^1, \tau_S^1)$
and $(q_S^2,\kappa_S^2,\allowbreak {\frak N}_S^2, \tau_S^2)$
of $\rho_{a,b,c}$-pairs parameterized by $S$.
By hypothesis, $\kappa_S^1=\kappa_S^2=:\kappa_S$, and we have
isomorphisms $\psi_S\colon E_S^1\lra E_S^2$ and $\chi_S\colon
{\frak N}_S^1\lra {\frak N}_S^2$
with
$$
\tau_S^1=({\id}_{\EL[\kappa_S]}\otimes\pi_S^*\chi_S)^{-1}\circ\tau_S^2
\circ ({\psi_S^{\otimes a}}^{\oplus c}).
$$
In particular, there is an isomorphism
\begin{eqnarray*}
U\otimes\Oh _S&\stackrel{\pi_{S*}(q^1_S\otimes{\id}_{\pi_X^*\Oh _X(n)})}{\lra}&
\pi_{S_*} \bigl(E^1_S\otimes\pi_X^*\Oh _X(n)\bigr)\q\lra\\
&\stackrel{\pi_{S*}(\psi_S\otimes{\id}_{\pi_X^*\Oh _X(n)})}{\lra}&
\pi_{S_*} \bigl(E^2_S\otimes\pi_X^*\Oh _X(n)\bigr)
\stackrel{\pi_{S*}(q^2_S\otimes{\id}_{\pi_X^*\Oh _X(n)})^{-1}}{\lra}
U\otimes\Oh _S.
\end{eqnarray*}
This yields a morphism $\Xi_S\colon S\lra \GL(U)$ and
$\Delta_S:=(\det)\circ \Xi_S\colon S\lra \C^*$.
Let $T:=S\times_{\C^*}\C^*$ be the fibre product taken w.r.t.\ $\Delta_S$
and $\C^*\lra\C^*$, $z\lma z^p$. The morphism $\eta\colon T\lra S$ is then
a $p$-sheeted \'etale covering coming with the projection
map $\widetilde{\Delta}\colon T\lra \C^*$.  In the following, we set
$\widetilde{\Delta}^e:=(z\lma z^e)\circ \widetilde{\Delta}$, $e\in\Z$.
One has
$\widetilde{\Delta}^p=\Delta_S\circ \eta$.
By construction, the morphism
$$
T\stackrel{\widetilde{\Delta}^{-1}\times (\Xi_S\circ\eta)}
{\lra}\C^*\times\GL(U)\stackrel{\mathop{\rm mult}}{\lra}
\GL(U)
$$
factorizes over a morphism $\Xi\colon T\lra \SL(U)$.
The quotient family defined by the morphism
$$
T\stackrel{\Xi\times (\beta_1\circ\eta)}{\lra} \SL(U)\times {\frak T}
\stackrel{\Gamma}{\lra}{\frak T}
$$
is just $(\widetilde{q}^1_S, \kappa_S\circ\eta, \eta^*{\frak
N}^1_S, (\eta\times {\id}_X)^*\tau^1_S)$ with
$$
\widetilde{q}^1_S\colon U\otimes\pi_X^*\Oh _X(-n)
\stackrel{\Xi^*(m^{-1})\otimes{\id}_{\pi_X^*\Oh _X(-n)}} {\lra}
U\otimes\pi_X^*\Oh _X(-n) \stackrel{(\eta\times
{\id}_X)^*q_S^1}{\lra} (\eta\times {\id}_X)^* E_S^1.
$$
The assertion of the proposition is that this family is equivalent
to the quotient family $((\eta\times \id_X)^* q_S^2,
\kappa_S\circ\eta, \eta^*{\frak N}^2_S,
(\eta\times \id_X)^*\tau^2_S)$.
But this is easily seen, using 
$$
\widetilde{\psi}_T:=\widetilde{\Delta}\cdot
((\eta\times {\id}_X)^*\psi_S^{-1})\colon \allowbreak
(\eta\times {\id}_X)^* E_S^2\lra (\eta\times {\id}_X)^* E_S^1
$$
and
$\widetilde{\chi}_T:=\widetilde{\Delta}^{a-rb}\cdot (\eta^*\chi_S^{-1})\colon
\eta^*{\frak N}_S^2\lra \eta^*{\frak N}_S^1$.
\end{proof}
\subsubsection{The Gieseker space and map}
Choose a Poincar\'e sheaf ${\cal P}$ on $\Jac^d\times X$. By our assumptions
on $n$, the sheaf
$$
\G_1:=\ul{\Hom}\Bigl(\bigwedge^rU\otimes\Oh _{\Jac^d}, \pi_{\Jac^d*} \bigl(
{\cal P}\otimes\pi_X^*\Oh _X(rn)\bigr)\Bigr)
$$
is locally free. We set ${\Bbb G}_1:=\Pe(\G_1^\vee)$.
By replacing ${\cal P}$ with ${\cal P}\otimes\pi_{\Jac^d}^*$
(sufficiently ample),
we may assume that $\Oh _{{\Bbb G}_1}(1)$ is very ample.
Let ${\frak d}\colon {\frak T}\lra \Jac^d$ be the morphism
associated with $\bigwedge^r E_{\frak T}$, and let ${\frak A}_{\frak T}$
be a line bundle on ${\frak T}$ with
$\bigwedge^r E_{\frak T}\cong ({\frak d}\times\id_X)^*{\cal P}\otimes
\pi_{\frak T}^*{\frak A}_{\frak T}$.
Then
$$
\bigwedge^ r \bigl(q_{\frak T}\otimes {\id}_{\pi_X^*\Oh _X(n)}\bigr)
\colon \bigwedge^r U\otimes\Oh _{\frak T}\lra
({\frak d}\times{\id}_X)^*{\cal P}
\otimes\pi_X^*\Oh _X(rn)
\otimes\pi_{\frak T}^*{\frak A}_{\frak T}
$$
defines a morphism $\iota_1\colon {\frak T}\lra {\Bbb G}_1$
with $\iota_1^*\Oh _{{\Bbb G}_1}(1)={\frak A}_{\frak T}$.
\par
Set $J^{d,m}:=\Jac^{d}\times \Jac^{m}$. The sheaf
$$
\G_2:=\ul{\Hom}\Bigl(U_{a,c}\otimes \Oh _{J^{d,m}},
\pi_{J^{d,m}*}\bigl(\pi_{{\rm Jac}^d\times X}^*({\cal P})^{\otimes
b} \otimes \pi_{{\rm Jac}^m\times X}^*({\cal
L})\otimes\pi_X^*\Oh _X(na)\bigr)\Bigr)
$$
on $J^{d,m}$ is also locally free. Set ${\Bbb
G}_2:=\Pe(\G_2^\vee)$. Making use of Remark~\ref{Indep}, it is
clear that we can assume $\Oh _{{\Bbb G}_2}(1)$ to be very ample.
The homomorphism
\begin{eqnarray*}
U_{a,c}\otimes \Oh _{\frak T}&\lra& {E_{\frak T}^{\otimes a}}^{\oplus c}
\otimes\pi_X^*\Oh _X(na)\q\lra\\
&\lra& ({\frak d}\times{\id}_X)^*{\cal P}^{\otimes b}
\otimes \EL[\kappa_{\frak T}]
\otimes\pi_X^*\Oh _X(na)
\otimes \pi_{\frak T}^*\bigl({\frak A}_{\frak T}^{\otimes b}
\otimes{\frak N}_{\frak T}\bigr)
\end{eqnarray*}
provides a morphism $\iota_2\colon {\frak T}\lra {\Bbb G}_2$
with
$\iota_2^*\Oh _{{\Bbb G}_2}(1)=
{\frak A}_{\frak T}^{\otimes b}\otimes{\frak N}_{\frak T}$.
Altogether, setting ${\Bbb G}:={\Bbb G}_1\times {\Bbb G}_2$
and $\iota:=\iota_1\times \iota_2$, we have an injective and
$\SL(U)$-equivariant morphism
$$
\iota\colon {\frak T}\lra {\Bbb G}.
$$
Linearize the $\SL(U)$-action on ${\Bbb G}$ in $\Oh _{\Bbb G}(\eps,1)$
with
$$
\eps\q:=\q {p-a\cdot\delta\over r\delta},
$$
and denote by ${\Bbb G}^{\eps-(s/p)s}$ the sets of points in ${\Bbb G}$
which are $\SL(U)$-(semi/poly)stable w.r.t.\ the given linearization.
\begin{Thm}
\label{Quot}
For $n$ large enough, the following two properties hold true:
\par
{\rm i)} The preimages $\iota^{-1}({\Bbb G}^{\eps-(s/p)s})$ consist
exactly of those points $t\in {\frak T}$ for which $(E_t, M_t,\allowbreak \tau_t)$
(notation as in Rem.~\ref{lin}) is a $\delta$-(semi/poly)stable
$\rho_{a,b,c}$-pair of type $(d,r,m)$.
\par
{\rm ii)} The restricted
morphism $\iota_{|\iota^{-1}\bigl({\Bbb G}^{\eps-ss}\bigr)}\colon
\iota^{-1}\bigl({\Bbb G}^{\eps-ss}\bigr)
\lra {\Bbb G}^{\eps-ss}$ is proper.
\end{Thm}
The proof of this theorem will be given in a later section.
\subsubsection{Proof of Theorem~\ref{Main}}
Set ${\frak T}^{\delta-(s)s}:=\iota^{-1}\bigl({\Bbb G}^{\eps-(s)s}\bigr)$.
Theorem~\ref{Quot} now shows that the categorical quotients
$$
{\cal M}(\rho_{a,b,c})_{d/r/m}^{\delta-(s)s}
\q:=\q
{\frak T}^{\delta-(s)s}
\catqot\SL(U)
$$
exist and that
${\cal M}(\rho_{a,b,c})_{d/r/m}^{\delta-s}$ is an orbit space.
Proposition~\ref{LUP} and Proposition~\ref{Glue} tell us that we have
a natural transformation of the functor
$\ul{\hbox{\rm M}}(\rho_{a,b,c})^{\delta-(s)s}_{d/r/m}$ into
the functor of points of ${\cal M}(\rho_{a,b,c})_{d/r/m}^{\delta-(s)s}$.
The asserted minimality property of
${\cal M}(\rho_{a,b,c})_{d/r/m}^{\delta-ss}$
and ${\cal M}(\rho_{a,b,c})_{d/r/m}^{\delta-s}$'s being a coarse
moduli space follow immediately from the universal property of the
categorical quotient.
Finally, the assertion about the closed points is a consequence
of the ``polystable" part of~\ref{Quot}.
Therefore, Theorem~\ref{Main} is settled for representations of
the form $\rho_{a,b,c}$.
\par
For an arbitrary representation $\rho$, we may find $a,b,c$ and
a decomposition $\rho_{a,b,c}=\rho\oplus\ol{\rho}$.
Define ${\frak T}(\rho)$ as the closed
subscheme of ${\frak T}$ where the homomorphism
$$
\widetilde{\tau}_{\frak T}\colon
E_{S,\rho_{a,b,c}}={E_S^{\otimes a}}^{\oplus c}\otimes\Bigl(\bigwedge^rE_S
\Bigr)^{\otimes -b}
\lra \EL[\kappa_{\frak T}]\otimes \pi_{\frak T}^*{\frak N}_{\frak T}
$$
vanishes on $E_{S,\ol{\rho}}$.
Set ${\frak T}(\rho)^{\delta-(s)s}:={\frak T}(\rho)\cap
{\frak T}^{\delta-(s)s}$. It follows that the categorical quotients
$$
{\cal M}(\rho)^{\delta-(s)s}_{d/r/m}\q:=\q {\frak T}(\rho)^{\delta-(s)s}
\catqot\SL(U)
$$
also exist.
By our characterization of $\ul{\hbox{\rm M}}(\rho)^{\delta-(s)s}_{d/r/m}$
as a closed subfunctor of the functor
$\ul{\hbox{\rm M}}(\rho_{a,b,c})_{d/r/m}^{\delta-(s)s}$,
the theorem follows likewise for $\rho$.
\par
Next, we let ${\frak T}_{\mathop{\rm surj}}$
be the open subscheme
of ${\frak T}$ consisting of those points $t$
for which $\widetilde{\tau}_{{\frak T}|\{t\}\times X}$ is surjective
and set ${\frak T}(\rho)_{\mathop{\rm surj}}^{\delta-s}:=
{\frak T}(\rho)^{\delta-s}\cap  {\frak T}_{\mathop{\rm surj}}$.
Thus, there is a section
$$
\sigma_{{\frak T}(\rho)_{\mathop{\rm surj}}^{\delta-s}}\colon
{\frak T}(\rho)_{\mathop{\rm surj}}^{\delta-s}\times X
\lra \Pe(E_{S,\rho}).
$$
Moreover, ${\frak P}(E_S)\times^{\GL(r)} F$ is a closed subscheme
of $\Pe(E_{S,\rho})$. Now, we define ${\frak
T}(F,\rho)^{\delta-\rho-s}$ as the closed subscheme of those
points $t\in {\frak T}(\rho)_{\mathop{\rm surj}}^{\delta-s}$ for
which the restricted morphism $\sigma_{{\frak
T}(\rho)_{\mathop{\rm surj}}^{\delta-s}|\{t\} \times X}$
factorizes over ${\frak P}(E_S)\times^{\GL(r)} F$. Since the
action of $\SL(U)$ on ${\frak T}(\rho)^{\delta-s}$ is closed, the
categorical quotient ${\frak T}(\rho)_{\mathop{\rm
surj}}^{\delta-s} \catqot\SL(U)$ exists as an open subscheme of
the moduli space ${\cal M}(\rho)^{\delta-s}_{d/r/m}$, whence
$$
{\cal M}(F,\rho)_{d/r/h}^{\delta-\rho-s}\q:=\q
{\frak T}(F,\rho)^{\delta-\rho-s}\catqot \SL(U)
$$
exists as a closed subscheme of ${\frak T}(\rho)_{\mathop{\rm surj}}^{\delta-s}
\catqot\SL(U)$ and hence as a locally closed subscheme of
${\cal M}(\rho)^{\delta-s}_{d/r/m}$ as asserted.
\subsubsection{Proof of Theorem~\ref{Quot}}
\subsubsection{Notation and Preliminaries}
The remarks about one parameter subgroups of $\SL(r)$ in Section~\ref{1PSG}
naturally apply to one parameter subgroups of $\SL(U)$.
We set
$$
\gamma_p^{(i)}\q:=\q \bigl(\, {\underbrace{i-p,...,i-p}_{i\times}},
\underbrace{i,...,i}_{(p-i)\times}\,\bigr),\q i=1,...,p-1.
$$
Given a basis $\ul{u}=(u_1,...,u_p)$ of $U$ and a weight vector
$\widetilde{\ul{\gamma}}=\sum_{i=1}^{p-1}\beta_i \gamma_p^{(i)}$,
we denote the corresponding one parameter subgroup of $\SL(U)$ by
$\la(\ul{u}, \widetilde{\ul{\gamma}})$. We hope that these
conventions will not give rise to too much confusion. Having
fixed a basis $\ul{u}= (u_1,...,u_p)$ of $U$ and an index
$l\in\{\,1,...,p\,\}$, we set
$U_{\ul{u}}^{(l)}:=\langle\,u_1,...,u_l\,\rangle$.
\par
Let $\rho_{{\Bbb G}_1}$ be the natural linearization of the $\SL(U)$-action
on ${\Bbb G}_1$ in $\Oh _{{\Bbb G}_1}(1)$.
Then, we write $\mu_{{\Bbb G}_1}(.,.)$ instead of
$\mu_{\rho_{{\Bbb G}_1}}(.,.)$. In the same way, $\mu_{{\Bbb G}_2}(.,.)$
is to be read.
Finally, $\mu^\eps_{\Bbb G}(.,.):=\eps\mu_{{\Bbb G}_1}(.,.)+\mu_{{\Bbb G}_2}
(.,.)$, i.e., $\mu^\eps_{\Bbb G}(.,.)=\mu_{\rho^\eps_{\Bbb G}}(.,.)$,
where $\rho_{\Bbb G}^\eps$ stands for the linearization of the $\SL(U)$-action
on ${\Bbb G}$ in $\Oh (\eps,1)$, $\eps\in \Q_{>0}$.
\par
Let $q\colon U\otimes\Oh _X(-n)\lra E$ be a generically surjective
homomorphism and $E$ a vector bundle of degree $d$ and rank $r$.
Set $Z:=H^0(\det(E)(rn))$. Then
$h:=\bigwedge^r(q\otimes\id_{\Oh _X(n)}) \in \Hom(\bigwedge^r U,
Z)$ is non trivial, and we can look at $[h]\in\Pe(\Hom(\bigwedge^r
U, Z)^\vee)$. On this space, there is a natural $\SL(U)$-action.
Then, it is well-known (e.g., \cite{HL2}) that for any basis
$\ul{u}=(u_1,...,u_p)$ and any two weight vectors
$\ul{\gamma}^i=(\gamma_1^i,...,\gamma_p^i)$ with
$\gamma^i_1\le\cdots\le\gamma^i_p$ and $\sum \gamma^i_j=0$,
$i=1,2$,
$$
\mu\bigl(\la(\ul{u},\ul{\gamma}^1),[h]\bigr)+
\mu\bigl(\la(\ul{u},\ul{\gamma}^2),[h]\bigr)\q=\q
\mu\bigl(\la(\ul{u},\ul{\gamma}^1+\ul{\gamma}^2),[h]\bigr)
$$
and for every $l\in\{\, 1,...,p-1\,\}$
\begin{equation}
\label{Weights}
\mu\bigl(\la(\ul{u},\ul{\gamma}^{(l)}),[h]\bigr)\q=\q
p\rk E_l- l r.
\end{equation}
Here, $E_l\subset E$ stands for the subbundle generated by
$q\bigl(U_{\ul{u}}^{(l)}\otimes\Oh _X(-n)\bigr)$.
\subsubsection{Sectional semistability}
\begin{Thm}
\label{Sect1} Fix the tuple $(d,r,m)$ and $a,b,c$ as before. Then,
there exists an $n_1$, such that for every $n\ge n_1$ and every
$\delta$-(semi)stable $\rho_{a,b,c}$-pair $(E,M,\tau)$, the
following holds true: For every weighted filtration
$(E^\bullet,\ul{\alpha})$, $E^\bullet: 0\subset E_1\subset\cdots
\subset  E_s\subset E$, of $E$
$$
\sum_{j=1}^s\alpha_{i_j}\bigl(\chi(E(n))\rk E_i-h^0(E_i(n))\rk E\bigr)
\q+\q \delta\cdot\mu_{\rho}\bigl(E^\bullet, \ul{\alpha};\tau\bigr)
\q (\ge) \q 0.
$$
\end{Thm}
\begin{proof}
First, suppose we are given a weighted filtration $(E^\bullet,\ul{\alpha})$,
$E^\bullet: 0\subset E_1\subset\cdots \subset  E_s\subset E$,
such that
$E_i(n)$ is globally generated and $H^1(E_i(n))=0$ for
$i=1,...,s$.
Then, for $i=1,...,s$,
\begin{eqnarray*}
&&\chi(E(n))\rk E_i-h^0(E_i(n))\rk E\\
&=&
\bigl(d+r(n+1-g)\bigr)\rk E_i - \bigl(\deg (E_i)+\rk E_i(n+1-g)\bigr)r
\\
&=& d\rk E_i-\deg(E_i)r,
\end{eqnarray*}
so that the claimed condition follows from $(E,M,\tau)$ being
$\delta$-(semi)stable.
\par
Next, recall that we have found a universal positive constant $C_1$
depending only on $r$, $a$, and $\delta$, such that for every $d$,
every semistable $\rho_{a,b,c}$-pair $(E,M,\tau)$,
and every non-trivial subbundle
$E^\p$ of $E$
$$
\mu(E^\p)\q\le\q {d\over r}+C_1.
$$
If we fix another positive constant $C_2$, then the set of
isomorphy classes of vector bundles $E^\p$ such that
$\mu(E^\p)\ge (d/r)-C_2$, $\mu_{\max}(E^\p)\le (d/r)+C_1$, and
$1\le\rk E^\p \le r-1$ is bounded. From this, we infer that there
is a natural number $n(C_2)$, such that for every $n\ge n(C_2)$,
every semistable $\rho_{a,b,c}$-pair $(E,M,\tau)$ of type
$(d,r,m)$, and every proper subbundle $E^\p$ of $E$
\begin{itemize}
\item either $\mu(E^\p)<(d/r)-C_2$
\item or $E^\p(n)$ is globally generated and $H^1(E^\p(n))=0$.
\end{itemize}
Moreover, the Le Potier-Simpson estimate (cf.\ \cite{LP},
Lemma~7.1.2 and proof of 7.1.1, p.~106)
gives in the first case
$$
h^0\bigl(E^\p(n)\bigr)\ \le\  \rk E^\p\cdot\left({\rk
E^\p-1\over\rk E^\p} \Bigl[{d\over r}+C_1+n+1\Bigr]_++ {1\over \rk
E^\p}\Bigl[{d\over r}-C_2+n+1\Bigr]_+ \right),
$$
i.e., for large $n$
$$
h^0\bigl(E^\p(n)\bigr)\q\le\q \rk E^\p\left({d\over r} + n+1+
(r-2) C_1 -{C_2\over r}\right),
$$
and
thus
\begin{eqnarray*}
\chi\bigl(E(n)\bigr)\rk E^\p-h^0\bigl(E^\p(n)\bigr) r&\ge& K(g,r,C_1,C_2)\\
&:=&
-r(r-1)g-r(r-1)(r-2)C_1+C_2.
\end{eqnarray*}
Our contention is now that for $C_2$
with
$$
K(g,r,C_1,C_2)\q >\q \delta\cdot a\cdot (r-1)
$$
and $n_1:=n(C_2)$, the theorem holds true.
\par
So, assume that we are given a weighted filtration
$(E^\bullet,\ul{\alpha})$ with $E^\bullet: 0\subset
E_1\subset\cdots \subset  E_s\subset E$ and
$\ul{\alpha}=(\alpha_1,...,\alpha_s)$. Let ${j_1},...,{j_t}$ be
the indices such that $\mu(E_{j_i})\ge d/r-C_2$, for $i=1,...,t$,
so that $E_{j_i}(n)$ is globally generated and
$H^1(E_{j_i}(n))=0$, $i=1,...,t$. We let
$\widetilde{\jmath}_1,...,\widetilde{\jmath}_{s-t}$ be the indices in $\{\,
1,...,s\,\}\setminus\{\, j_1,...,j_t\,\}$ in increasing order. We
introduce the weighted filtrations $(E_1^\bullet,\ul{\alpha}_1)$
and $(E^\bullet_2,\ul{\alpha}_2)$ with $E_1^\bullet: 0\subset
E_{j_1}\subset\cdots\subset E_{j_t}\subset E$,
$\ul{\alpha}_1:=(\alpha_{j_1},...,\alpha_{j_t})$ and
$E_2^\bullet: 0\subset E_{\widetilde{\jmath}_1}\subset\cdots \subset
E_{\widetilde{\jmath}_{s-t}}\subset E$, $\ul{\alpha}_2=
(\alpha_{\widetilde{\jmath}_1},...,\alpha_{\widetilde{\jmath}_{s-t}})$.
Lemma~\ref{Add}~ii) yields
$$
\mu_\rho\bigl(E^\bullet,\ul{\alpha};\tau\bigr)\q\ge\q
\mu_\rho\bigl(E_1^\bullet,\ul{\alpha}_1;\tau\bigr)-\left(\sum_{i=1}^{s-t}
\alpha_{\widetilde{\jmath}_{i}}\right)\cdot \delta\cdot a\cdot(r-1),
$$
whence
\begin{eqnarray*}
& &\sum_{j=1}^{s}\alpha_j\bigl(\chi(E(n))\rk E_j-h^0(E_j(n))\rk E\bigr)
\q+\q \delta\cdot\mu_{\rho}\bigl(E^\bullet,\ul{\alpha};\tau\bigr)
\\
&\ge&
\sum_{i=1}^t\alpha_{j_i}\bigl(\chi(E(n))\rk E_{j_i}-h^0(E_{j_i}(n))\rk E\bigr)
\q+\q \delta\cdot\mu_{\rho}\bigl(E^\bullet_1,\ul{\alpha}_1;\tau\bigr)\\
&&+ \sum_{i=1}^{s-t}\alpha_{\widetilde{\jmath}_i} \bigl(\chi(E(n))\rk
E_{\widetilde{\jmath}_i}-h^0(E_{\widetilde{\jmath}_i}(n))\rk
E\bigr)-\left(\sum_{i=1}^{s-t}
\alpha_{\widetilde{\jmath}_i}\right)\cdot \delta\cdot a\cdot (r-1)\\
&\ge &
\sum_{i=1}^t\alpha_{j_i}\bigl(\chi(E(n))\rk E_{j_i}-h^0(E_{j_i}(n))\rk E\bigr)
\q+\q \delta\cdot\mu_{\rho}\bigl(E^\bullet_1,\ul{\alpha}_1;\tau\bigr)\\
&&+ \left(\sum_{i=1}^{s-t}
\alpha_{\widetilde{\jmath}_i}\right)K(g,r,C_1,C_2)
-\left(\sum_{i=1}^{s-t} \alpha_{\widetilde{\jmath}_i}\right)\cdot
\delta\cdot  a\cdot (r-1).
\end{eqnarray*}
Since this last expression is positive by assumption,
we are done.
\end{proof}
\subsubsection*{The implication
{$t\in\iota^{-1}({\Bbb G}^{\eps-(s)s})
\Rightarrow (E_t,M_t,\tau_t)$} is $\delta$-(semi)stable}
To begin with, we fix a constant $K$ with the property that
$$
rK\q>\q\max\bigl\{\, d(s-r)+\delta\cdot a\cdot (r-1)\,|\, s=1,...,r-1\,
\bigr\}.
$$
Now, let $t=[q\colon U\otimes\Oh _X(-n)\lra E_t, M_t,\tau_t]$ be
a point with $\iota(t)\in {\Bbb G}^{\eps-(s)s}$.
\par
We first claim that there can be no subbundle $E^\p\subset E_t$
with $\deg(E^\p)\ge d+K$. Let $E^\p$ be such a subbundle.
Then, for every natural number $n$,
$$
h^0\bigl(E^\p(n)\bigr)\q\ge\q d+K+\rk E^\p\bigl(n+1-g).
$$
Let $\widetilde{E}$ be the subbundle of $E_t$ which is generated
by $\Imm(\ev\colon H^0(E^\p(n) )\otimes\Oh _X(-n)\lra E_t)$. Thus,
$H^0(\widetilde{E}(n))=H^0(E^\p(n))$ and $\widetilde{E}$ is
generically generated by its global sections. Now, choose a basis
$u_1,...,u_i$ for $H^0(E^\p(n))$, complete it to a basis
$\ul{u}:=(u_1,...,u_p)$ of $U$, and set
$\la:=\la(\ul{u},\gamma_p^{(i)})$. Then, we have seen that
\begin{eqnarray*}
\mu_{{\Bbb G}_1}\bigl(\la, \iota_1(t)\bigr)&=& p\cdot\rk \widetilde{E}-h^0(\widetilde{E}(n))\cdot r\\
          &\le & p\cdot\rk E^\p-h^0(E^\p(n))\cdot r.
\end{eqnarray*}
Our discussion preceding Lemma~\ref{Add} applies to $\SL(U)$ as well, whence
$$
\mu_{{\Bbb G}_2}\bigl(\la,\iota_2(t)\bigr)\q\le\q a\cdot (p-i).
$$
Therefore,
\begin{eqnarray*}
\mu^\eps_{\Bbb G}\bigl(\la,\iota(t)\bigr)&=&
\eps\cdot \mu_{{\Bbb G}_1}\bigl(\la, \iota_1(t)\bigr)
+
\mu_{{\Bbb G}_2}\bigl(\la,\iota_2(t)\bigr)\\
&\le& {p-a\cdot\delta\over r\cdot\delta}\bigl(p\cdot\rk E^\p-h^0(E^\p(n))\cdot r
\bigr) + a\cdot (p-i)\\
\\
&=& {p^2\rk E^\p\over r\delta}-{pa\rk E^\p\over r}-{ph^0(E^\p(n))\over \delta}+pa.
\end{eqnarray*}
Next, we multiply the last expression by the positive number $r\delta/p$
in order to obtain
\begin{eqnarray*}
&& p\rk E^\p-rh^0(E^\p(n))+\delta a\bigl(r-\rk E^\p\bigr)\\
&\le & \bigl(d+r(n+1-g)\bigr)\rk E^\p-r\bigl(d+K+\rk E^\p(n+1-g)\bigr)
+\delta a\bigl(r-1\bigr)\\
&=& d\bigl(\rk E^\p-r\bigr)+\delta a \bigl(r-1\bigr)-rK\q <\q 0,
\end{eqnarray*}
by our choice of $K$. This obviously contradicts the assumption
$\iota(t)\in {\Bbb G}^{\eps-ss}$. We can also assume that $d+K>0$.
Set $C_3:=(r-1)d/r +K$. Then our arguments show that
$\iota(t)\in {\Bbb G}^{\eps-ss}$ implies
$$
\mu_{\max}(E_t)\q\le\q {d\over r}+C_3,
$$
\sl independently \rm
of the number $n$ with which we performed the construction of ${\Bbb G}$.
An argument similar to the one used in the proof of Theorem~\ref{Sect1}
shows that a $\rho_{a,b,c}$-pair $(E,M,\tau)$ is $\delta$-(semi)stable,
if and only if
for every weighted filtration $(E^\bullet,\ul{\alpha})$, such that
$$
\mu(E_j)\q\ge\q {d\over r}-{\delta\cdot a\cdot (r-1)\over r}
\q=\q
{d\over r}-C_1,\qquad j=1,...,s,
$$
one has
$$
M\bigl(E^\bullet,\ul{\alpha}\bigr)
\q+\q \delta\cdot\mu_{{\rho_{a,b,c}}}\bigl(E^\bullet,\ul{\alpha};\tau\bigr)
\q (\ge) \q 0.
$$
Therefore, we choose $n$ so large that for every vector bundle
$E^\p$ with $d/r+C_3\ge \mu_{\max}(E^\p)$, $\mu(E^\p)\ge d/r-C_1$,
and $1\le\rk E^\p\le r-1$, one has that $E^\p(n)$ is globally generated
and $H^1(E^\p(n))$ vanishes.
\par
Now, let $(E^\bullet,\ul{\alpha})$ be a weighted filtration with
$$
\mu(E_j)\q\ge\q {d\over r}-C_1,\qquad j=1,...,s.
$$
Fix a basis $\ul{w}=(w_1,...,w_r)$ of $W:=\C^r$,
and let $W^\bullet: 0\subset W_{\ul{w}}^{(i_1)}\subset\cdots\subset
W_{\ul{w}}^{(i_s)}\subset W$ be the associated flag, $i_j:=\rk E_j$,
$j=1,...,s$.
Let $\ul{u}=(u_1,...,u_p)$ be a basis of $U$ such that there are
indices $l_1,...,l_s$ with $U_{\ul{u}}^{(l_j)}=H^0(E_j(n))$, $j=1,...,s$.
Define
$$
\widetilde{\ul{\gamma}}
\q:=\q
\sum_{j=1}^s \alpha_{j} \gamma_p^{(l_j)}.
$$
We also set, for $j=1,...,s+1$, $l_{s+1}:=p$, $l_0:=0$,
$i_{s+1}:=r$, $i_0:=0$,
$$
{\gr}_j(U,\ul{u}):=
U_{\ul{u}}^{(l_{j})}/U_{\ul{u}}^{(l_{j-1})}
=H^0(E_j/E_{j-1}(n)),\quad \hbox{and}
\quad {\gr}_j(W,\ul{w}):= W_{\ul{w}}^{(i_j)}/W_{\ul{w}}^{(i_{j-1})}.
$$
The fixed bases $\ul{w}$ for $W$ and $\ul{u}$ for $U$ provide
us with isomorphisms
$$
U\cong \bigoplus_{j=1}^{s+1} {\gr}_j(U,\ul{u}),\q \hbox{and}\q
W\cong \bigoplus_{j=1}^{s+1} {\gr}_j(W,\ul{w}).
$$
Let $J^a:=\{\, 1,...,s\,\}^{\times a}$. For every index $\ul{\iota}\in J^a$,
we set
$$
W_{\ul{\iota},\ul{w}}:={\gr}_{\iota_1}(W,\ul{w})
\otimes\cdots\otimes {\gr}_{\iota_a}(W,\ul{w}).
$$
Analogously, we define $U_{\ul{\iota},\ul{w}}$.
Moreover, for $k\in \{\,1,...,c\,\}$
and $\ul{\iota}\in J^a$, we let $W_{\ul{\iota},\ul{w}}^k$ be the
subspace of $W_{a,c}:={W^{\otimes a}}^{\oplus c}$ which is
$W_{\ul{\iota},\ul{w}}$ living in the
$k$-th copy of $W^{\otimes a}$ in $W_{a,c}$,
and similarly we define $U_{\ul{\iota},\ul{u}}^k$.
The spaces $W^k_{\ul{\iota},\ul{w}}$ and  $U_{\ul{\iota},\ul{u}}^k$,
$k\in\{\, 1,...,c\,\}$
and $\ul{\iota}\in J^a$, are eigenspaces for the actions
of the one parameter subgroups $\la(\ul{w},\gamma^{(i_j)})$
and $\la(\ul{u},\gamma_p^{(l_j)})$, respectively, $j=1,...,s$.
Define
\begin{equation}
\label{nu}
\nu_j(\ul{\iota})\q:=\q
\#\bigl\{\, \iota_i\le j\,|\,\ul{\iota}=(\iota_1,...,
\iota_a),\ i=1,...,a\,\bigr\}.
\end{equation}
Then $\la(\ul{w},\gamma^{(i_j)})$ acts on $W^k_{\ul{\iota},\ul{w}}$ with
weight $\nu_j(\ul{\iota})\cdot r-a\cdot i_j$, and
$\la(\ul{u},\gamma_p^{(l_j)})$ acts on $U_{\ul{\iota},\ul{u}}^k$
with weight $\nu_j(\ul{\iota})\cdot p-a\cdot l_j$.
\par
Let $Z_t:=H^0(\det(E_t)^{\otimes b}\otimes M_t\otimes\Oh _X(na))$.
Then $\iota_2(t)\in\Pe(\Hom(U_{a,c},Z_t)^\vee)$ can be represented
by a homomorphism
$$
L_t\colon U_{a,c}\lra Z_t.
$$
One readily verifies
\begin{equation}
\label{Filtiii}
\begin{array}{ll}
\mu_{{\Bbb G}_1}\bigl(\la(\ul{u},\widetilde{\ul{\gamma}}), [L_t]\bigr)\q=
&\hskip 6cm\q
\end{array}
\end{equation}
$$
-\min
\Bigl\{\, \sum_{j=1}^s\alpha_{j}
\bigl(\nu_j(\ul{\iota})\cdot p- a\cdot l_j\bigr)\,|\, k\in\{\,1,..,c\,\},
\ul{\iota}\in J^a: U^k_{\ul{\iota},\ul{u}}\not\subset\ker L_t\,
\Bigr\}.
$$
Next, we observe that we can choose a small open subset $X_0\subset X$
over which $E_t$ and $M_t$ are trivial
and there is an isomorphism $\psi\colon E_{t|X_0}\cong W\otimes \Oh _{X_0}$
with $\psi(E^\bullet_{|X_0})=W^\bullet\otimes\Oh _{X_0}$.
This trivialization and the $\rho_{a,b,c}$-pair $(E_t,M_t,\tau_t)$
provide us with
$$
l_t\colon W_{a,c}\otimes\Oh _{X_0}\lra
\Bigl(\bigwedge^r W\Bigr)^{\otimes b}\otimes\Oh _{X_0}.
$$
We observe that for every $k\in\{\,1,...,c\,\}$ and every $\ul{\iota}\in J^a$
\begin{equation}
\label{Filti}
W_{\ul{\iota},\ul{w}}^k\otimes\Oh _{X_0}\not\subset \ker l_t\q\Leftrightarrow
\q
U_{\ul{\iota},\ul{u}}^k\not\subset \ker L_t,
\end{equation}
and that
\begin{equation}
\label{Filtii}
\begin{array}{ll}
\mu_{\rho_{a,b,c}}\bigl(E^\bullet,\ul{\alpha};\tau_t\bigr)\q=
&\hskip 7cm\q\
\end{array}
\end{equation}
$$
 -\min
\Bigl\{\, \sum_{j=1}^s\alpha_{j}
\bigl(\nu_j(\ul{\iota})\cdot r- a\cdot i_j\bigr)\,|\, k\in\{\,1,..,c\,\},
\ul{\iota}\in J^a: W^k_{\ul{\iota},\ul{w}}\otimes\Oh _{X_0}
\not\subset\ker l_t\,
\Bigr\}.
$$
\par
Now, let $k_0\in \{\,1,...,c\,\}$ and $\ul{\iota}_0\in J^a$ be such that
the minimum in~(\ref{Filtiii}) is achieved by
$\sum_{j=1}^s\alpha_{j}
(\nu_j(\ul{\iota}_0)\cdot p- a\cdot l_j)$ and
$U^{k_0}_{\ul{\iota}_0,\ul{u}}
\not\subset\ker L_t$.
We obtain
\begin{eqnarray*}
0 &(\le)& \mu_{\Bbb G}^\eps\bigl(\la(\ul{u},\widetilde{\ul{\gamma}}),
\iota(t)\bigr)\\
&=& \eps\cdot \mu_{{\Bbb G}_1}\bigl(\la(\ul{u},
\widetilde{\ul{\gamma}}), \iota_1(t)\bigr)
+\mu_{{\Bbb G}_2}\bigl(\la(\ul{u},\widetilde{\ul{\gamma}}),
\iota_2(t)\bigr)\\
&=& \eps\cdot \sum_{j=1}^s\alpha_{j}\bigl(p\rk E_j-h^0(E_j(n))r\bigr)+
\sum_{j=1}^s\alpha_{j}
\bigl(\nu_j(\ul{\iota}_0)\cdot p- a\cdot l_j\bigr)
\\
&=& {p-a\delta\over r\delta}\sum_{j=1}^s\alpha_{j}
\bigl(p\rk E_j-h^0(E_j(n))r\bigr)
+\sum_{j=1}^s\alpha_{j}
\bigl(\nu_j(\ul{\iota}_0)\cdot p- a\cdot h^0(E_j(n))\bigr)\\
&=&\sum_{j=1}^s\alpha_{j}\biggl({p^2\rk E_j\over r\delta}-{pa\rk E_j\over r}
-{ph^0(E_j(n))\over\delta}\biggr)+\sum_{j=1}^s\alpha_{j}
\nu_j(\ul{\iota}_0)\cdot p.
\end{eqnarray*}
We multiply this inequality by $r\delta/p$ and find
$$
0 \q(\le)\q
\sum_{j=1}^s\alpha_{j}\bigl(p\rk E_j-r h^0(E_j(n))\bigr)+
\delta\sum_{j=1}^s
\alpha_{j}\bigl(\nu_j(\ul{{\iota}_0})r-a\rk E_j\bigr).
$$
Since $h^1(E_j(n))=0$, $j=1,...,s$, we have
$p\rk E_j-r h^0(E_j(n))=d\rk E_j-r\deg(E_j)$, $j=1,...,s$.
Moreover, $\rk E_j=i_j$, by definition, and 
$\mu_{{\rho_{a,b,c}}}(E^\bullet,\ul{\alpha};\tau_t)\ge
\sum_{j=1}^s
\alpha_{j}(\nu_j(\ul{{\iota}_0})r-ai_j)$,
by~(\ref{Filti}) and~(\ref{Filtii}), whence we finally see
$$
M\bigl(E^\bullet,\ul{\alpha}\bigr)+\delta\cdot\mu_{{\rho_{a,b,c}}}
\bigl(E^\bullet,\ul{\alpha};\tau_t\bigr)\q(\ge)\q 0,
$$
as required.
\subsubsection*{The implication
$(E_t,M_t,\tau_t)$ is $\delta$-(semi)stable
$\Rightarrow$ $t\in\iota^{-1}({\Bbb G}^{\eps-(s)s})$}
By the Hilbert-Mumford criterion, we have to show that for every basis
$\ul{u}=(u_1,...,u_p)$ of $U$ and
every weight vector $\widetilde{\ul{\gamma}}
=(\gamma_1,...,\gamma_p)$ with $\gamma_1\le\cdots\le\gamma_p$
and $\sum_{i=1}^{p}\gamma_i=0$
$$
\mu^\eps_{\Bbb G}\bigl(\la(\ul{u},\widetilde{\ul{\gamma}}),\iota(t)\bigr)
\q=\q
\eps\mu_{{\Bbb G}_1}\bigl(\la(\ul{u},\widetilde{\ul{\gamma}}),\iota_1(t)
\bigr)+
\mu_{{\Bbb G}_2}\bigl(\la(\ul{u},\widetilde{\ul{\gamma}}),\iota_2(t)
\bigr)\q
(\ge)\q 0.
$$
So, let $\ul{u}=(u_1,...,u_p)$ be an arbitrary basis for $U$
and $\widetilde{\ul{\gamma}}=\sum_{i=1}^{p-1}\beta_i\gamma_p^{(i)}$
a weight vector. Let $l_1,...,l_v$ be the indices with
$\beta_{l_h}\neq 0$, $h=1,...,v$.
For each $h\in\{\, 1,...,v\,\}$, let $E_{l_h}$ be the subbundle
of $E_t$ generated by $\Imm(U_{\ul{u}}^{(l_h)}\otimes\Oh _X(-n)\lra E_t)$.
Note that for $h^\p\ge h$ we will have $E_{l_{h^\p}}=E_{l_h}$
if and only if $U_{\ul{u}}^{(l_{h^\p})}\subset H^0(E_{l_h}(n))$.
We let $E^\bullet: 0=:E_0\subset E_1\subset\cdots\subset
E_s\subset E_{s+1}:=E$
be the filtration by the distinct vector bundles occurring among
the $E_{l_h}$'s.
\par
Recall that we know (\ref{Weights})
$$
\mu_{{\Bbb G}_1}\bigl(\la(\ul{u},\widetilde{\ul{\gamma}}),\iota_1(t)
\bigr)=\sum_{h=1}^v \beta_{l_h}\bigl(p\rk E_{l_h}- l_h r\bigr)
\ge \sum_{h=1}^v \beta_{l_h}\bigl(p\rk E_{l_h}-h^0(E_{l_h}(n)) r\bigr).
$$
Set, for $j=1,...,s$,
$$
\alpha_{j}\q:=\q \sum_{h: E_{l_h}=E_j}\beta_{l_h},
$$
so that we see
\begin{equation}
\label{weighti}
\mu_{{\Bbb G}_1}\bigl(\la(\ul{u},\widetilde{\ul{\gamma}}),\iota_1(t)
\bigr)
\q\ge\q
\sum_{j=1}^s\alpha_{j}\bigl(p\rk E_j-h^0(E_j(n))r\bigr).
\end{equation}
\par
Next, we define for $j=0,...,s$
$$
h(j)\q:=\q \max\bigl\{\, h=1,...,v\,|\, U_{\ul{u}}^{(l_h)}
\subset h^0(E_j(n))\bigr\}.
$$
With these conventions, $h=h(j)+1$ is the minimal index, such that
$U^{(l_h)}_{\ul{u}}\otimes\Oh _X(-n)$ generically generates $E_{j+1}$,
$j=0,...,s$.
We now set
$$
\widetilde{\gr}_j\bigl(U,\ul{u}\bigr)
\q:=\q U_{\ul{u}}^{(l_{h(j-1)+1})}/
U^{(l_{h(j-1)})}_{\ul{u}},\qquad j=1,...,s+1.
$$
The space $\bigoplus_{j=1}^{s+1} \widetilde{\gr}_j(U,\ul{u})$
can be identified with a subspace of $U$, via $\widetilde{\gr}_j
(U,\ul{u})\cong\langle\, l_{h(j-1)}+1,...,\allowbreak l_{h(j-1)+1}\,\rangle$,
$j=1,...,s$.
\par
For any index tuple $\ul{\iota}=(\iota_1,...,\iota_a)\in J^a:=
\{\,1,...,s\,\}^{\times a}$, we define
$$
\widetilde{U}_{\ul{\iota},\ul{u}}\q:=\q \widetilde{\gr}_{\iota_1}(U,\ul{u})
\otimes\cdots\otimes \widetilde{\gr}_{\iota_a}(U,\ul{u}).
$$
Again, for $\ul{\iota}\in J^a$ and $k\in\{\,1,...,c\,\}$,
$\widetilde{U}^k_{\ul{\iota},\ul{u}}$ will be
$\widetilde{U}_{\ul{\iota},\ul{u}}$ viewed as a subspace of the $k$-th
summand of $U_{a,c}$.
\par
The effect of our definition of the $h(j)$'s is that
the spaces $\widetilde{U}^k_{\ul{\iota},\ul{u}}$,
$\ul{\iota}\in J^a$ and $k\in\{\,1,...,c\,\}$,
are eigenspaces \sl for all \rm the one parameter subgroups
$\la(\ul{u}, \gamma_p^{(l_h)})$, $h=1,...,v$,
with respect to the weight
$\nu_j(\ul{\iota})p-a l_h$, $\nu_j(\ul{\iota})$ as in (\ref{nu}).
\par
Now, let $\ul{w}=(w_1,...,w_r)$ be a basis for $W$ and
$W^\bullet:0\subset W_{\ul{w}}^{(i_1)}\subset\cdots\subset W_{\ul{w}}^{(i_s)}
\subset W$, $i_j:=\rk E_j$, $j=1,...,s$, the corresponding flag.
Then, the spaces $W^k_{\ul{\iota},
\ul{w}}$, $\ul{\iota}\in J^a$ and $k\in\{\,1,...,c\,\}$, are defined
as before.
We can find a small open set $X_0\subset X$, such that
\begin{itemize}
\item $M_t$ and $E_t$ are trivial over $X_0$,
\item there is an isomorphism $\psi\colon E_{t|X_0}\lra W\otimes \Oh _{X_0}$
      with $\psi(E^\bullet_{|X_0})=W^\bullet\otimes \Oh _{X_0}$,
\item $E_{t|X_0}\cong \bigoplus_{j=1}^{s+1} (E_j/E_{j-1})_{|X_0}$,
\item the homomorphism $\bigl(
\bigoplus_{j=1}^{s+1} \widetilde{\gr}_j(U,\ul{u})\bigr)\otimes\Oh _{X_0}(-n)
\lra \bigoplus_{j=1}^{s+1} (E_j/E_{j-1})_{|X_0}$ is surjective.
\end{itemize}
\par
As before, let $Z_t:=H^0(\det(E_t)^{\otimes b}\otimes M_t\otimes\Oh _X(na))$,
so that
$\iota_2(t)\in\Pe(\Hom(U_{a,c},\allowbreak Z_t)^\vee)$ induces a homomorphism
$$
\widetilde{L}_t\colon \bigoplus \widetilde{U}^k_{\ul{\iota},\ul{u}}\lra Z_t.
$$
Letting
$$
l_t\colon W_{a,c}\otimes\Oh _{X_0}\lra
\Bigl(\bigwedge^r W\Bigr)^{\otimes b}\otimes\Oh _{X_0}
$$
be the resulting homomorphism, we find
that for every $k\in\{\,1,...,c\,\}$ and every $\ul{\iota}\in J^a$
\begin{equation}
\label{Filt2i}
W_{\ul{\iota},\ul{w}}^k\otimes\Oh _{X_0}\not\subset \ker l_t\q\Leftrightarrow
\q
\widetilde{U}_{\ul{\iota},\ul{u}}^k\not\subset \ker \widetilde{L}_t.
\end{equation}
By Theorem~\ref{Sect1}, we have
\begin{equation}
\label{Filt2ii}
\sum_{j=1}^s\alpha_{j}\bigl(p\rk E_j-h^0(E_j(n))r\bigr)
+\delta\cdot
\mu_{\rho_{a,b,c}}\bigl(E^\bullet, (\alpha_1,...,\alpha_s);\tau_t\bigr)
\q(\ge)\q 0.
\end{equation}
Now, we choose $k_0\in\{\, 1,...,c\,\}$ and $\ul{\iota}_0\in J^a$
with $W^{k_0}_{\ul{\iota}_0,\ul{w}}\otimes\Oh _{X_0}\not\subset\ker l_t$
and $\mu_{\rho_{a,b,c}}(E^\bullet,\allowbreak(\alpha_1,...,\alpha_s);\tau_t)
\allowbreak
=\sum_{j=1}^s \alpha_{j}(\nu_j(\ul{\iota}_0)r-a\rk E_j)$.
Plugging this into (\ref{Filt2ii}) and multiplying by $p/(r\delta)$
yields
\begin{eqnarray*}
0 &(\le)&
\sum_{j=1}^s\alpha_{j}\biggl({p^2\rk E_j\over r\delta}-{pa\rk E_j\over r}
-{ph^0(E_j(n))\over\delta}\biggr)+\sum_{j=1}^s\alpha_{j}
\nu_j(\ul{\iota})\cdot p\\
&=&\eps\sum_{j=1}^s\alpha_j\bigl(p\rk E_j-h^0(E_j(n))r\bigr)
+\sum_{j=1}^s\alpha_{j}
\bigl(\nu_j(\ul{\iota}_0)\cdot p- a\cdot h^0(E_j(n))\bigr).
\end{eqnarray*}
By our definition of the $\alpha_{j}$, and (\ref{Filt2i}),
we know
\begin{eqnarray*}
\mu_{{\Bbb G}_2}\bigl(\la(\ul{u},\widetilde{\ul{\gamma}}),\iota_2(t)\bigr)
&\ge& \sum_{h=1}^{v} \beta_{l_h}\bigl(\nu_{j(h)}(\ul{\iota}_0)p-a l_h\bigr)
\\
&\ge& \sum_{j=1}^s \alpha_{j}\bigl(\nu_j(\ul{\iota}_0)p-ah^0(E_j(n))\bigr).
\end{eqnarray*}
Here, we have set $j(h)$ to be the element $j\in\{\, 1,...,s\,\}$
with $E_{l_h}=E_j$.
This together with (\ref{weighti}) finally shows
$\mu^\eps_{\Bbb G}(\la(\ul{u},\widetilde{\ul{\gamma}}),\iota(t))(\ge) 0$.

\subsubsection*{The identification of the polystable points}
By the Hilbert-Mumford criterion, a point $\iota(t)$ is polystable if and only if it is semistable and,
for every one parameter subgroup $\la$ of $\SL(U)$ with $\mu^\epsilon_{\Bbb G}(\la,\iota(t))=0$,
$\lim_{z\rightarrow\infty}\la(z)\cdot \iota(t)$ lies in the orbit of $\iota(t)$.
\par
Now, let $\ul{u}=(u_1,...,u_p)$ be a basis for $U$ and
$\widetilde{\ul{\gamma}}=\sum_{j=1}^s \beta_{l_j}\gamma_p^{(l_j)}$
be a weight vector with $\beta_{l_j}\neq 0$ and $l_j\in\{\,1,...,p-1\,\}$
such that $\mu^\eps_{\Bbb G}(\la(\ul{u},\widetilde{\ul{\gamma}}),\iota(t))
=0$.
Then, our previous considerations show that the following must
be satisfied
\begin{itemize}
\item $U_{\ul{u}}^{(l_j)}=H^0(E_{l_j}(n))$, $j=1,...,s$,
\item $E_{l_j}(n)$ is generated by global sections and $H^1(E_{l_j}(n))=0$.
\end{itemize}
Set $E_j:=E_{l_j}$, $i_j:=\rk E_j$, $\alpha_{j}:=
\beta_{l_j}$, $j=1,...,s$,
and choose a basis $w_1,...,w_r$ for $W$.
As before, we associate to these data a flag $W^\bullet$.
Consider the weighted filtration
$(E^\bullet,\ul{\alpha})$ with
$E^\bullet:0\subset E_1\subset\cdots\subset E_s\subset E_t$
and $\ul{\alpha}=(\alpha_1,...,\alpha_s)$,
so that the condition
$\mu^\eps_{\Bbb G}(\la(\ul{u},\widetilde{\ul{\gamma}}),\iota(t))=0$
becomes equivalent to
$$
M\bigl(E^\bullet,\ul{\alpha}\bigr)
+\delta\cdot\mu_{\rho_{a,b,c}}\bigl(E^\bullet,\ul{\alpha}; \tau_t\bigr)
\q=\q
0.
$$
Let $t_\infty:=\lim_{z\rightarrow\infty}\la(z)\cdot t$ and $(E_{t_\infty},
M_{t_{\infty}}, \tau_{t_\infty})$ be the corresponding $\rho_{a,b,c}$-pair.
Then, clearly $M_{t_\infty}\cong M_t$, and it is well known that
$E_{t_\infty}\cong \bigoplus_{j=1}^{s+1} E_j/E_{j-1}$.
Let $U_{a,c}:=\bigoplus U^{\widetilde{g}_i}$ be the decomposition
of $U_{a,c}$ into eigenspaces w.r.t.\ the $\C^*$-action coming from
$\la(\ul{u},\widetilde{\ul{\gamma}})$, and $\widetilde{g}_{i_0}
=-\mu_{{\Bbb G}_2}(\la(\ul{u},\widetilde{\ul{\gamma}}),\iota_2(t))$.
If $L_t\colon U_{a,c}\allowbreak\lra Z_t$ and $L_{t_\infty}\colon U_{a,c}
\allowbreak \lra
Z_{t_\infty}= Z_t$ are the homomorphisms representing
$t$ and $t_\infty$, respectively, then $L_{t_\infty}$ is just
the restriction of $L_t$ to $U^{\widetilde{g}_{i_0}}$ extended
by zero to the other weight spaces.
As we have seen before, the condition that
$L_{t_\infty}$ be supported only on $U^{\widetilde{g}_{i_0}}$
is equivalent to the fact that over each open subset
$X_0$ over which $\tau_{t_\infty}$ is surjective
and we have a trivialization $\psi\colon E_{t_\infty|X_0}
\lra W\otimes\Oh _{X_0}$ with $\psi(E^\bullet_{|X_0})= W^\bullet
\otimes\Oh _{X_0}$, the induced morphism $X_0\lra \Pe(W_{a,c})$
lands in $\Pe(W^{g_0})$,
where $W^{g_0}$ is
the eigenspace for the weight $g_0:=-\mu_{\rho_{a,b,c}}(
E^\bullet,\ul{\alpha};\tau_{t_\infty})$.
Thus, we have shown that $(E_t,M_t,\tau_t)$ being $\delta$-polystable
implies that $\iota(t)$ is a polystable point.
The converse is similar.
\subsubsection*{The properness of the Gieseker map}
In this section, we will prove that the Gieseker morphism
$\iota$ is proper, using the (discrete) valuative criterion.
\par
Thus, let $(C,0)$ be the spectrum of a DVR $R$ with quotient
field $K$.
Suppose we are given a morphism $h\colon
C\lra {\Bbb G}^{\eps-ss}$ which lifts
over $\Spec K$ to ${\frak T}$.
This means that we are given a quotient family $(q_K\colon
U\otimes\pi_X^*\Oh _X(-n)\lra E_K, \kappa_K,\tau_K)$
of $\rho_{a,b,c}$-pairs parameterized by $\Spec K$ (we left out
${\frak N}_K$, because it is trivial).
This can be extended to a certain family
$(\widetilde{q}_C\colon U\otimes\pi_X^*\Oh _X(-n)\lra\widetilde{E}_C,
\kappa_C, \tau_C)$, consisting of
\begin{itemize}
\item a surjection $\widetilde{q}_C$
      onto the flat family $\widetilde{E}_C$, where
      $\widetilde{E}_{C|\{0\}\times X}$ may have torsion
\item the continuation $\kappa_C$ of $\kappa_K$ into $0$
\item a homomorphism $\tau_C\colon
      \widetilde{E}_C^{\otimes a^{\oplus c}}
      \lra \det(\widetilde{E})^{\otimes b}\otimes\EL[\kappa_C]$
      whose restriction to $\{0\}\times X$ is non trivial and
      whose restriction to $\Spec K\times X$ differs from $\tau_K$
      by an element in $K^*$.
\end{itemize}
The resulting datum $L\colon U_{a,c}\lra
\pi_{C*}(\det(\widetilde{E}_{C})^{\otimes b}\otimes \EL[\kappa_C]
\otimes \pi_X^*\Oh _X(na))$
defines a
morphism $C\lra {\Bbb G}_2$ which coincides with the second component
$h_2$ of $h$.
\par
Set $E_C:=\widetilde{E}_C^{\vee\vee}$. This is a reflexive sheaf
on the smooth surface $C\times X$, whence it is locally free and
thus flat over $C$. Therefore, we have a
family
$$
q_C\colon U\otimes \pi_X^*\Oh _X(-n) \lra E_C
$$
where the kernel of the homomorphism
$
U\otimes\Oh _X(-n) \lra E_{C|\{0\}\times X}
$
is isomorphic to the torsion ${\cal T}$ of $\widetilde{E}_{C|\{0\}\times X}$.
One gets a homomorphism $\bigwedge^r U\otimes\Oh _C
\lra \pi_{C*}(\det(\widetilde{E}_C)\otimes\pi_X^*\Oh _X(rn))$ which defines
a morphism $C\lra {\Bbb G}_1$ which coincides with the first component
$h_1$ of $h$.
\par
Set $E_0:=E_{C|\{0\}\times X}$.
Our claim is that $H^0(q_{C|\{0\}\times X}\otimes\id_{\pi_X^*\Oh _X(n)})
\colon U\lra H^0(E_0(n))$ must be injective.
This implies, in particular, that $\widetilde{E}_{C|\{0\}\times X}$
is torsion free and, hence, $E_C=\widetilde{E}_C$ and $q_C=\widetilde{q}_C$.
If $H:=\ker(H^0(q_{C|\{0\}\times X}\otimes\id_{\pi_X^*\Oh _X(n)}))$
is non trivial, we choose a basis $u_1,...,u_j$
for $H$ and complete it to a basis $\ul{u}=(u_1,...,u_p)$ of $U$.
Set $\ol{H}=\langle\, u_{j+1},...,u_p\,\rangle$.
We first note (\ref{Weights})
$$
\mu_{{\Bbb G}_1}\bigl(\la(\ul{u},\gamma_p^{(j)}),h_1(0)\bigr)\q=\q
-jr.
$$
The spaces $H_l:=H^{\otimes l}\otimes \ol{H}^{\otimes (a-l)}$, $l=1,...,a$,
are the eigenspaces of $U^{\otimes a}$ for the $\C^*$-action
coming from $\la(\ul{u},\gamma_p^{(j)})$. Let $H^k_l$
be $H_l$ embedded into the $k$-th component of $U_{a,c}$,
$k=1,...,c$, $l=1,...,a$. For every $k\in\{\,1,...,c\,\}$ and
every $l\in\{\,1,...,a\,\}$, $H^k_l\otimes \Oh _X(-n)$ generates
a torsion subsheaf of $\widetilde{E}^{\otimes a^{\oplus c}}_{C|\{0\}\times X}$, so that
$H^k_l\subset\ker L$. This implies
$$
\mu_{{\Bbb G}_2}\bigl(\la(\ul{u},\gamma_p^{(j)}),h_2(0)\bigr)\q
\q=\q
\mu_{{\Bbb G}_2}\bigl(\la(\ul{u},\gamma_p^{(j)}), [L]\bigr)
\q=\q -aj,
$$
and thus
$$
\mu^\eps_{\Bbb G}\bigl(\la(\ul{u},\gamma_p^{(j)}), h(0)\bigr)
\q=\q -\eps j r-a j\q<\q 0,
$$
in contradiction to the assumption $h(0)\in {\Bbb G}^{\eps-ss}$.
\par
We identify $U$ with its image in $H^0(E_0(n))$. Let $K$ be a positive
constant
such that $rK>\max\{\, d(s-r)+\delta a (r-1)\,|\,s=1,...,r-1\,\}$.
We assert that for every non-trivial and proper quotient bundle $Q$
of $E_0$ we must have $\deg Q\ge- K-(r-1)g$. For this, let
$Q$ be the minimal destabilizing quotient bundle.
Set $E^\p:=\ker(E\lra Q)$. It suffices to show
that  $\deg Q<- K-(r-1)g$ implies $\dim(H^0(E^\p(n))\cap U)
\ge d+K+\rk E^\p(n+1-g)$, because then a previously given argument applies.
Note that we have an exact sequence 
$$
\begin{CD}
0 @>>> H^0(E^\p(n))\cap U @>>> U @>>>
H^0(Q(n)).
\end{CD}
$$
Assume first that $h^0(Q(n))=0$. Thus, $\dim(H^0(E^\p(n))\cap U)
=p = d+\rk E^\p(n+1-g) +(r-\rk E^\p)(n+1-g)\ge d+\rk E^\p(n+1-g)+n+1-g$.
Since we can assume $n+1-g>K$, this is impossible.
\par
Therefore, the Le Potier-Simpson estimate gives
$h^0(Q(n))\le \deg Q+\rk Q(n+1)$ and thus
\begin{eqnarray*}
\dim\bigl(H^0(E^\p(n))\cap U\bigr)&\ge & p-h^0(Q(n))\\
&\ge & d+r(n+1-g) -\deg Q - \rk Q(n+1)\\
&=& d-\deg Q -g\rk Q +(r-\rk Q)(n+1-g)\\
&\ge& d-\deg Q-g(r-1)+\rk E^\p(n+1-g).
\end{eqnarray*}
This gives the claim.
We see
$$
\mu_{\min}(E_0)\q\ge\q {-K-(r-1)g\over \rk Q}\q\ge\q -K-(r-1)g.
$$
This bound does not depend on $n$. Since the family of isomorphy
classes of vector bundles $G$  of degree $d$ and rank $r$
with $\mu_{\min}(G)\ge -K-(r-1)g$ is bounded, we can choose
$n$ so large that $H^1(G(n))=0$ for every such vector bundle.
In particular, $H^1(E_0(n))=0$, i.e., $U=H^0(E_0(n))$.
This means that the family $(\widetilde{q}_C, \kappa_C, \tau_C)$
we started with is a quotient family of $\rho_{a,b,c}$-pairs
parameterized by $C$ and thus defines a morphism from $C$
to ${\frak T}$ which lifts $h$.
By Theorem~\ref{Quot} i), this morphism factorizes through
${\frak T}^{\delta-ss}$, and we are done.
\section{Examples}
This section is devoted to the treatise of the known examples within our
general context. First, we discuss two important methods of
simplifying the stability concept. Second, we will consider some
easy specializations of the moduli functors. Then, we briefly
discuss the variation of the stability parameter and prove an
``asymptotic irreducibility" result. Afterwards, we turn to the
examples. In the examples, we will show how many of the known
stability concepts and constructions of the moduli spaces over
curves can be obtained via our construction. In two cases we will
see that our results give a little more than previous
constructions. We have also added the stability concept for conic
bundles of rank $4$. The main aim of the examples is to
illustrate that the complexity of the stability concept only
results from the complexity of the input representation
$\rho\colon\GL(r)\lra\GL(V)$ and to illustrate how the
understanding of $\rho$ can be used to simplify the stability
concept.
\subsection{Simplifications of the stability concept}
In this part, we will formulate several ways of restating the concept of
$\delta$-semistability in different, easier ways which will be used
in the study of examples to recover the known notions
of semistability.
The first one uses a well-known additivity property to reduce the
stability conditions to conditions on subbundles.
The second one generalizes this to a method working for all representations.
This provides the mechanism alluded to in the introduction.
The third one is a method to express the concept of $\delta$-semistability
for $\rho$-pairs
associated with a direct sum $\rho=\rho_1\oplus\cdots\oplus \rho_n$
of representations in a certain sense in terms of the semistability
concepts corresponding to the summands $\rho_i$.
Further methods of simplifying the semistability concept
will be discussed in the examples.
\subsubsection*{A certain additivity property}
Let $\rho\colon \GL(r)\lra \GL(V)$ be a representation such that the following property holds
true:
For any basis $\ul{w}=(w_1, ..., w_r)$ of $\C^r$, any two weight
vectors $\ul{\gamma}_1$ and $\ul{\gamma}_2$, and any point
$[l]\in\Pe(V)$
\begin{equation}
\label{addi}
\mu_{\rho}\bigl(\la(\ul{w},\ul{\gamma}_1+\ul{\gamma}_2), [l]\bigr)
\q=\q
\mu_{\rho}\bigl(\la(\ul{w},\ul{\gamma}_1),[l]\bigr)+
\mu_{\rho}\bigl(\la(\ul{w},\ul{\gamma}_2),[l]\bigr).
\end{equation}
Now, let $(E,M,\tau)$ be a $\rho$-pair and $\delta$ a positive
rational number.
For every weighted filtration $(E^\bullet,\ul{\alpha})$,
$E^\bullet:0\subset E_1\subset\cdots\subset E_s\subset E$,
the definition of $\mu_\rho(E^\bullet, \ul{\alpha};\tau)$ and (\ref{addi})
imply
$$
M\bigl(E^\bullet,\ul{\alpha}\bigr)+\delta\mu_\rho\bigl(
E^\bullet,\ul{\alpha}; \tau\bigr)\hskip 5cm
$$
$$
=\q
\sum_{j=1}^s \alpha_{j}\Bigl(
\bigl(d\rk E_j-r\deg E_j\bigr)+\delta\mu_\rho\bigl(E_j,
\tau\bigr)\Bigr).
$$
We see that the semistability condition becomes a condition
on subbundles of $E$:
The $\rho$-pair $(E,M,\tau)$ is $\delta$-(semi)stable, if and only if for
every non trivial proper subbundle $E^\p$ of $E$
one has
\begin{equation}
\label{addii}
\mu(E^\p)\q(\le)\q \mu(E)+{\mu_\rho\bigl(E^\p,\tau\bigr)
\over\rk E^\p \rk E}.
\end{equation}
\subsubsection*{The general procedure}
Let $\rho\colon \GL(r)\lra \GL(V)$ be a representation on $V$
and $\rho^\p\colon\allowbreak
\SL(r)\lra \GL(V)$ its restriction to $\SL(r)$. We fix a basis $\ul{w}
=(w_1,...,w_r)$ of $\C^r$. This basis determines a maximal torus $T\subset
\SL(r)$. First, we observe that the Hilbert-Mumford criterion can be restated
in the following form:
A point $[l]\in \Pe(V)$ is $\rho^\p$-(semi)stable, if and only if for every
element $g\in\SL(r)$ and every weight vector $\ul{\gamma}=(\gamma_1,...,
\gamma_r)$ with $\gamma_1\le\cdots\le\gamma_r$ and $\sum \gamma_i=0$
\begin{equation}
\label{addni}
\mu_{\rho}\bigl(\la(\ul{w},\ul{\gamma}),g\cdot[l]\bigr)\q (\ge)\q 0.
\end{equation}
The representation $\rho_{|T}\colon T\lra \GL(V)$ yields
a decomposition
$$
V\q =\q \bigoplus_{\chi\in X(T)} V_\chi
$$
with
$$
V_\chi\q:=\q \bigl\{\, v\in V\, |\, \rho(t)(v)=\chi(t)\cdot v\q\forall t\in T
\,\bigr\}.
$$
The set $\mathop{\rm ST}(\rho):=\{\, \chi\in X(T)\,|\, V_\chi\neq\langle
0\rangle\,\}$ is the \it set of states of $\rho$\rm.
We look at the rational polyhedral cone
$$
C:=\bigl\{\, (\gamma_1,...,\gamma_r)\,|\, \gamma_1\le\cdots\le\gamma_r,
\sum\gamma_i=0\,\bigr\}=\R_{\ge 0}\cdot \gamma^{(1)}+\cdots+\R_{\ge 0}\cdot
\gamma^{(r-1)}.
$$
For every subset $A\subset \mathop{\rm ST}(\rho)$, we obtain a decomposition
$$
C=\bigcup_{\chi\in A} C^\chi_A\q \hbox{with}\q
C^\chi_A:=\bigl\{\,\ul{\gamma}\in C\,|\, \langle\la(\ul{w},\ul{\gamma}),
\chi\rangle\le \langle\la(\ul{w},\ul{\gamma}),
\chi^\p\rangle\ \forall \chi^\p\in A\,\bigr\}.
$$
Here, $\langle.,.\rangle$ is the natural pairing between one parameter
subgroups and characters.
The cones $C^\chi_A$ are also rational polyhedral cones and one has
$$
C^\chi_A\cap C^{\chi^\p}_A\q=\q C^\chi_A\ \cap\ \bigl\{\,\ul{\gamma}\,|\,
\langle \la(\ul{w},\ul{\gamma}),\chi-\chi^\p\rangle=0\,\bigr\},
$$
so that two cones intersect in a common face.
Therefore, for each $A$, we get a fan decomposition of $C$.
For each edge of a cone $C^\chi_A$, there is a minimal integral generator.
For $A\subset \mathop{\rm ST}(\rho)$ and $\chi\in A$, we let
$K^\chi_A$ be the set of those generators and $K_A=\bigcup_{\chi\in A}
K^\chi_A$. The set $K_A$ obviously contains $\{\, \gamma^{(1)},...,
\gamma^{(r-1)}\,\}$, and we call $A$ \it critical\rm, if $K_A$ is strictly
bigger than $\{\, \gamma^{(1)},...,\gamma^{(r-1)}\,\}$.
Now, for each point $[l]\in \Pe(V)$, we set $\mathop{\rm ST}(l):=
\{\, \chi\,|\, l_{|V_\chi}\not\equiv 0\,\}$. Moreover, an element
$g\in \SL(r)$ is called \it critical for $[l]$\rm, if the set
$\mathop{\rm ST}(g\cdot l)$ is critical.
\par
We observe that for a point $[l]\in\Pe(V)$ and a weight vector
$\ul{\gamma}\in C$ one has
$$
\mu_{\rho}\bigl(\la(\ul{w},\ul{\gamma}),[l]\bigr)
\q=\q -\min\bigl\{\, \langle\la(\ul{w},\ul{\gamma}),\chi\rangle\,|\,
\chi\in\mathop{\rm ST}(l)\,\bigr\}.
$$
This means that Equation~(\ref{addi}) remains valid, if there exists a
character $\chi\in\mathop{\rm ST}(l)$, such that $C^\chi_{\mathop{\rm ST}(l)}$
contains both $\ul{\gamma}_1$ and $\ul{\gamma}_2$.
We infer
\begin{Cor}
\label{addnii}
A point $[l]\in \Pe(V)$ is $\rho^\p$-(semi)stable if and only if it satisfies
the following two conditions:
\begin{enumerate}
\item[\rm 1.] For every element $g\in \SL(r)$ and every $i\in\{\,1,...,r-1\,\}$
      $$
      \mu_\rho\bigl(\la(\ul{w},\gamma^{(i)}), g\cdot[l]\bigr)\q (\ge)\q 0.
      $$
\item[\rm 2.] For every $g\in\SL(r)$ which is critical for $[l]$ and every weight vector
      $\ul{\gamma}
      \in K_{\mathop{\rm ST}(g\cdot l)}\setminus\{\, \gamma^{(1)},...,
      \gamma^{(r-1)}\,\}$
      $$
      \mu_\rho\bigl(\la(\ul{w},\ul{\gamma}),g\cdot [l]\bigr)\q(\ge)\q 0.
      $$
\end{enumerate}
In particular, it suffices to test {\rm (\ref{addni})} for the
weight vectors belonging to the finite set
$$
K_\rho\q:=\q\bigcup_{A\subset\mathop{\rm ST}(\rho)} K_A.
$$
\end{Cor}
\begin{Rem}
A similar procedure works for all semisimple groups $G$. Indeed, one fixes
a pair $(B,T)$ consisting of a Borel subgroup of $G$ and a maximal torus
$T\subset B$. With analogous arguments, one obtains decompositions of the
Weyl chamber $W(B,T)$. See \cite{ABCD} for a precise discussion.
\end{Rem}
Let's now turn to the $\rho$-pairs. Let $W^\bullet$ be the complete
flag $0\subset\langle w_1\rangle\subset\cdots\subset \langle w_1,...,\allowbreak w_{r-1}
\rangle\subset \C^r$. For a $\rho$-pair $(E,M,\phi)$ and a filtration
$0\subset E_1\subset \cdots\subset E_{r-1}\subset E$ with $\rk E_i=i$, $i=1,...,r-1$,
we define $\mathop{\rm ST}(E^\bullet)$ as follows: Choose an
open subset $U$ and
a trivialization $\psi\colon E_{|U}\lra \Oh _U^{\oplus r}$ with
$\psi(E^\bullet_{|U})=W^\bullet\otimes\Oh _U$. Then, for each $\chi\in\mathop{\rm ST}
(\rho)$, there is a rational map
$$
U\lra \Pe(E_{\rho|U})\cong \Pe(V)\times U\lra \Pe(V)\dasharrow\Pe(V_\chi).
$$
An element $\chi\in \mathop{\rm ST}(\rho)$ now belongs to  $\mathop{\rm ST}(E^\bullet)$,
if and only if this rational map is defined on a non-empty subset of $U$.
As before, one verifies that $\mathop{\rm ST}(E^\bullet)$ is well-defined.
The filtration $E^\bullet$ is called \it critical for $\phi$\rm, if
$\mathop{\rm ST}(E^\bullet)$ is critical.
Corollary~\ref{addnii} now shows
\begin{Thm}
\label{addniii} {\rm i)} The $\rho$-pair $(E,M,\phi)$ is
$\delta$-(semi)stable, if and only if it meets the following two
requirements:
\begin{enumerate}
\item[\rm 1.] For every proper non-trivial subbundle $E^\p$ of $E$
      $$
      \mu(E^\p)\q(\le)\q \mu(E)+\frac{\mu_\rho(E^\p,\phi)}{\rk E^\p\rk E}.
      $$
\item[\rm 2.] For every filtration $E^\bullet$ which is critical for $\phi$ and
      every
      element $\sum_{j=1}^s \alpha_j \gamma^{(i_j)}\in
      K_{\mathop{\rm ST}(E^\bullet)}\setminus\{\, \gamma^{(1)},...,\allowbreak
      \gamma^{(r-1)}\,\}$, $\alpha_j>0$, $i_j:=\rk E_j$, $j=1,...,s$,
      $$
      \begin{array}{lcl}
      M\bigl(0\subset E_{1}\subset \cdots \subset E_{s}\subset E,
      (\alpha_1,...,\alpha_s)\bigr)\ + 
      \\
	  \\
      \delta\cdot
      \mu_\rho\bigl(0\subset E_{1}\subset \cdots \subset E_{s}\subset E,
      (\alpha_1,...,\alpha_s);\phi\bigr)&(\ge)& 0.
      \end{array}
      $$
\end{enumerate}
\end{Thm}
\subsubsection*{Direct sums of representations}
Let $\rho_i\colon \GL(r)\lra \GL(V_i)$ be representations of the
general linear group and assume there is an integer $\alpha$ with
$\rho_i(z\cdot\id_{\C^r})=z^\alpha\cdot\id_V$ for all $z\in\C^*$,
$i=1,...,t$. Define $\rho:= \rho_1\oplus\cdots\oplus\rho_t$. Note
that for every rank $r$ vector bundle $E$ one has
$E_\rho=E_{\rho_1}\oplus\cdots\oplus E_{\rho_t}$. The following
result is a counterpart to Theorem~\ref{ch2} in the first part.
\begin{Prop}
\label{ch3}
Let $(E,M,\tau)$ be a $\rho$-pair of type $(d,r,m)$ and $\delta\in\Q_{>0}$.
Then the following conditions are equivalent:
\begin{enumerate}
\item[\rm 1.]
$(E,M,\tau)$ is $\delta$-semistable ($\delta$-polystable).
\item[\rm 2.] There exist pairwise distinct
indices $\iota_1,...,\iota_s\in\{\,1,...,t\,\}$, $s\le t$, such that
$$
j\in\bigl\{\,\iota_1,...,\iota_s\,\bigr\}
\q \Rightarrow\ (\Leftrightarrow)\q
\tau_{|E_{\rho_{j}}}\colon E_{\rho_{j}}
\lra M  \hbox{ is non-zero,}
$$
and positive rational numbers $\sigma_1,...,\sigma_s$
with $\sum_{j=1}^s\sigma_j=1$ such that for every
weighted filtration $(E^\bullet,\ul{\alpha})$
$$
M\bigl(E^\bullet,\ul{\alpha})+\delta\biggl(
\sum_{j=1}^s \sigma_j\mu_{\rho_{\iota_j}}\bigl(
E^\bullet,\ul{\alpha};\tau_{|E_{\rho_{\iota_j}}}\bigr)\biggr)\q\ge\q0.
$$
(And if equality holds
\begin{itemize}
\item $E\cong \bigoplus_{j=1}^{s+1} E_j/E_{j-1}$
\item the $\rho_{\iota_j}$-pair $(E,M,\tau|_{E_{\rho_{\iota_j}}})$ is equivalent to the 
${\rho_{\iota_j}}$-pair 
$(E,M,\tau|_{E_{\rho_{\iota_j}}^{\gamma_j}})$,
$\gamma_j:=-\mu_{\rho_{\iota_j}}(E^\bullet,\ul{\alpha};\tau|_{E_{\rho_{\iota_j}}})$, $j=1,...,s$ 
(Compare with definition of polystability above).
\end{itemize}
Here, $\ul{w}$ is a basis for $\C^r$,
$W^\bullet: 0\subset W_{\ul{w}}^{(\rk E_1)}\subset\cdots
\subset W_{\ul{w}}^{(\rk E_s)}\subset W$, and $\ul{\gamma}:=\sum_{j=1}^s
\alpha_j\gamma^{(\rk E_j)}$.)
\item[\rm 3.] There exist pairwise distinct
indices $\iota_1,...,\iota_s\in\{\,1,...,t\,\}$, $s\le t$, such that
$$
j\in\bigl\{\,\iota_1,...,\iota_s\,\bigr\}
\q \Rightarrow\ (\Leftrightarrow)\q
\tau_{|E_{\rho_{j}}}\colon E_{\rho_{j}}
\lra M  \hbox{ is non-zero,}
$$
and positive rational numbers $\sigma_1,...,\sigma_s$
with $\sum_{j=1}^s\sigma_j=1$ such that for every
positive integer $\nu$ with $\nu\sigma_j\in\Z_{>0}$, $j=1,...,s$,
the associated $(\rho_{\iota_1}^{\otimes \nu\sigma_1}\otimes\cdots\otimes
\rho_{\iota_s}^{\otimes \nu\sigma_s})$-pair
$$
\bigl(E,M^{\otimes (\nu\sigma_1+\cdots+\nu\sigma_s)},
\tau_{|E_{\rho_{\iota_1}}}^{\otimes \nu\sigma_1}\otimes\cdots\otimes
\tau_{|E_{\rho_{\iota_s}}}^{\otimes \nu\sigma_s}\bigr)
$$
of type $(d,r,\nu m)$
is $(\delta/\nu)$-semistable ($(\delta/\nu)$-polystable).
\end{enumerate}
\end{Prop}
\begin{proof}
To see the equivalence between 2.\ and 3., observe
that  $\Oh (\nu\sigma_1,...,\nu\sigma_s)$
provides an equivariant embedding of
$\Pe(V_{\iota_1})\times\cdots\times\Pe(V_{\iota_s})$
into $\Pe(S^{\nu\sigma_1}V_{\iota_1}\otimes\cdots\otimes S^{\nu\sigma_s}
V_{\iota_s})$. Via the canonical surjection
$V_{\iota_1}^{\otimes \nu\sigma_1}\otimes\cdots\otimes
V_{\iota_s}^{\otimes \nu\sigma_s}\lra
S^{\nu\sigma_1}V_{\iota_1}\otimes\cdots\otimes S^{\nu\sigma_s}
V_{\iota_s}$, the latter space becomes embedded into
$\Pe(V_{\iota_1}^{\otimes \nu\sigma_1}\otimes\cdots\otimes
V_{\iota_s}^{\otimes \nu\sigma_s})$, so that we have an equivariant
embedding
$\iota\colon\Pe(V_{\iota_1})\times\cdots\times\Pe(V_{\iota_s})
\hookrightarrow
\Pe(V_{\iota_1}^{\otimes \nu\sigma_1}\otimes\cdots\otimes
V_{\iota_s}^{\otimes \nu\sigma_s})$.
Since for every point
$x=(x_1,...,x_s)\in \Pe(V_{\iota_1})\times\cdots\times\Pe(V_{\iota_s})$
and every one parameter subgroup $\la\colon \C^*\lra \GL(r)$
$$
\sum_{j=1}^s\sigma_j\mu_{\rho_{\iota_j}}\bigl(\la, x_j\bigr)
\q=\q
{1\over \nu}\cdot\mu_{\rho_{\iota_1}^{\otimes \nu\sigma_1}\otimes\cdots
\otimes\rho_{\iota_s}^{\otimes \nu\sigma_s}}\bigl(\la,\iota(x)\bigr),
$$
the claimed equivalence is easily seen.
\par
For the equivalence between 1.\ and 3., we have to go into the GIT
construction of the moduli space of $\delta$-semistable $\rho$-pairs.
We choose $a,b,c$, such that $\rho$ is a direct summand
of $\rho_{a,b,c}$. Therefore, $\rho_i$ is also a direct summand
of $\rho_{a,b,c}$, $i=1,...,t$, so that we can assume
$\rho_i=\rho_{a,b,c}$ for $i=1,...,t$.
For a tuple $(\iota_1,...,\iota_s)$, positive
rational numbers $\sigma_1,...,\sigma_s$, and
$\nu\in \N$ as in the statement, we thus find
$$
\rho_{\iota_1}^{\otimes \nu\sigma_1}\otimes\cdots
\otimes\rho_{\iota_s}^{\otimes \nu\sigma_s}
\q=\q
\rho_{\nu a,\nu b,c^\p}
$$
for some $c^\p>0$. Recall that in our GIT construction of the
moduli space of $\delta$-semistable $\rho_{a,b,c}$ pairs of type
$(d,r,m)$, we had to fix some natural number $n$ which was large
enough. Being large enough depended on constants $C_1$, $C_2$,
$C_3$, and $K$ which in turn depended only on $d$, $r$, $a$, and
$\delta$. One now checks that $d$, $r$, $\nu a$, and $\delta/\nu$
yield exactly the same constants, so that the construction will
work also --- for all $\nu$ and all $c^\p$ --- for
$(\delta/\nu)$-semistable $\rho_{\nu a,\nu b,c^\p}$-pairs of type
$(d,r,\nu m)$. Fix such an $n$. We can now argue as follows. Set
$p:=d+r(n+1-g)$, and let $U$ be a complex vector space of
dimension $p$. Given a $\delta$-semistable $\rho_{a,b,c}$-pair
$(E,M,\tau)$ of type $(d,r,m)$, we can write $E$ as a quotient
$q\colon U\otimes\Oh _X(-n)\lra E$ where $H^0(q(n))$ is an
isomorphism. Set $Z:=\Hom(\bigwedge^r U, H^0(\det E(rn)))$ and
$W:= \Hom(U_{a,c}, H^0(\det E^{\otimes b}\otimes M\otimes
\Oh _X(na)))$. Then $(q\colon U\otimes \Oh _X(-n)\lra E, M,\tau)$
defines a Gieseker point $([z],[w_1,...,w_t])\in\Pe(Z^\vee)\times
\Pe({W^{\vee}}^{\oplus t})$ which is semistable for the
linearization of the $\SL(U)$-action in $\Oh (\eps, 1)$ with
$\eps=(p-a\delta)/(r\delta)$. By Theorem~\ref{ch2}, we find
indices $\iota_1,...,\iota_s$ and positive rational numbers
$\sigma_1,...,\sigma_s$ with $\sum_{j=1}^s\sigma_j=1$, such that
$w_{\iota_j}\neq 0$, $j=1,...,s$, and the point $([z],
[w_{\iota_1}],...,[w_{\iota_s}])\in \Pe(Z^\vee)\times
\Pe(W^\vee)\times\cdots\times\Pe(W^\vee)$ is semistable w.r.t.\ the
linearization of the $\SL(U)$-action in
$\Oh (\eps,\sigma_1,...,\allowbreak \sigma_s)$. As before, there is an embedding
$\iota:\Pe(Z^\vee)\times \Pe(W^\vee)\times\cdots\times\Pe(W^\vee)
\hookrightarrow \Pe(Z^\vee)\times \Pe({W^\vee}^{\otimes \nu})$ such
that the pullback of $\Oh (\nu\eps, 1)$ is $\Oh (\nu\eps,
\nu\sigma_1,..., \nu\sigma_s)$. The point $y:=\iota([z],
[w_{\iota_1}],...,[w_{\iota_s}])$ is thus semistable w.r.t.\ the
linearization in $\Oh (\nu\eps, \nu\sigma_1,..., \nu\sigma_s)$.
Now, the second component of $y$ is defined by the homomorphism
$U_{\nu a, c^\p}=U_{a,c}^{\otimes \nu} \lra H^0(\det E^{\otimes
b}\otimes M\otimes\Oh _X(na))^{\otimes \nu}$ obtained from $q$ and
the components $\tau_{|E_{\rho_{\iota_j}}}$, $j=1,...,s$.
Composing this homomorphism with the natural map
$H^0(\det E^{\otimes b}\otimes M\otimes\Oh _X(na))^{\otimes \nu}
\lra H^0(\det E^{\otimes \nu b}\otimes M^{\otimes \nu}
\otimes\Oh _X(n\nu a))$, we find a point $y^\p\in
\Pe(Z^\vee)\times \Pe({W^\p}^\vee)$, $W^\p:=\Hom(U_{\nu a,c^\p},
H^0(\det E^{\otimes \nu b}\otimes M^{\otimes \nu}
\otimes\Oh _X(n\nu a)))$.
The point $y^\p$ is semistable w.r.t.\ the linearization of the
$\SL(U)$-action in $\Oh (\nu \eps, 1)$.
By construction, $y^\p$ is the Gieseker point of
the quotient $\rho_{\nu a,\nu b,c^\p}$-pair
$(q\colon U\otimes \Oh _X(-n)\lra E, M^{\otimes \nu},\tau_{|E_{\rho_{\iota_1}}}
^{\otimes \nu\sigma_1}\otimes\cdots\otimes\tau_{|E_{\rho_{\iota_s}}}^{\otimes
\nu\sigma_s})$.
Since $\nu\eps=(p-(\nu a)(\delta/\nu))/(r\delta/\nu)$, we infer
that  $(E, M^{\otimes \nu},\tau_{|E_{\rho_{\iota_1}}}
^{\otimes \nu\sigma_1}\otimes\cdots\otimes\tau_{|E_{\rho_{\iota_s}}}^{\otimes
\nu\sigma_s})$ is $(\delta/\nu)$-semistable.
The converse and the polystable part are similar.
\end{proof}
\subsection{Some features of the moduli spaces}
Here, we will discuss several properties of the moduli spaces
which we have constructed.
\subsubsection*{Trivial specializations}
Let $\rho\colon \GL(r)\lra \GL(V)$ be a representation. Very
often, one fixes the determinant of the vector bundles under
consideration. So, let $L_0$ be a line bundle of degree $d$. If
we want to consider only $\rho$-pairs $(E,M,\tau)$ of type
$(d,r,m)$ with $\det E\cong L_0$, we say that the \it type of
$(E,M,\tau)$ \rm is $(L_0,r,m)$. We then obtain a closed
subfunctor $\ul{\hbox{M}}(\rho)^{\delta-(s)s}_{L_0/r/m}$ of
$\ul{\hbox{M}}(\rho)^{\delta-(s)s}_{d/r/m}$. Note that our
construction shows that we have a morphism ${\cal
M}(\rho)^{\delta-(s)s}_{d/r/m}\lra \Jac^d$, $[E,M,\tau]\lma [\det
E]$. Let ${\cal M}(\rho)^{\delta-(s)s}_{L_0/r/m}$ be the fibre
over $[L_0]$. This is then the moduli space for
$\ul{\hbox{M}}(\rho)^{\delta-(s)s}_{L_0/r/m}$.
\par
In the applications, the line bundle $M$ is traditionally fixed.
Having fixed a line bundle $M_0$ of degree $m$, we will speak of
\it $\rho$-pairs $(E,\tau)$ of type $(d,r,M_0)$\rm. This yields a
moduli functor $\ul{\hbox{M}}(\rho)^{\delta-(s)s}_{d/r/M_0}$
which is also a closed subfunctor of
$\ul{\hbox{M}}(\rho)^{\delta-(s)s}_{d/r/m}$. Its moduli space,
denoted by ${\cal M}(\rho)^{\delta-(s)s}_{d/r/M_0}$, is the fibre
over $[M_0]$ of the morphism ${\cal
M}(\rho)^{\delta-(s)s}_{d/r/m}\lra \Jac^m$, $[E,M,\tau]\lma [M]$.
\par
If we want to fix both $L_0$ and $M_0$, we speak of \it
$\rho$-pairs $(E,\tau)$ of type $(L_0,r,M_0)$\rm. The
corresponding moduli spaces are denoted by ${\cal
M}(\rho)^{\delta-(s)s}_{L_0/r/M_0}$.
\subsubsection*{Variation of $\delta$}
Given $\rho\colon \GL(r)\lra\GL(V)$, $d,r,m\in\Z$, $r>0$, we get
a whole family of moduli spaces ${\cal
M}(\rho)^{\delta-(s)s}_{d/r/m}$ parameterized by
$\delta\in\Q_{>0}$. This phenomenon was first studied by Thaddeus
in the proof of the Verlinde formula \cite{Th}. The papers
\cite{DH} and \cite{Th2} study the corresponding abstract GIT
version. Using these, one makes the following observations
\begin{enumerate}
\item There is an increasing sequence $(\delta_\nu)_{\nu\ge 0}$,
      $\delta_\nu\in \Q_{>0}$,
      $\nu=0,1,2,...$, which is discrete in $\R$,
      such that the concept of $\delta$-(semi)stability is
      constant within each interval $(\delta_\nu,\delta_{\nu+1})$, $\nu=0,1,2,...$,
      and, for given $\nu$,
      $\delta$-semistability for $\delta\in (\delta_\nu,\delta_{\nu+1})$
      implies $\delta_\nu$- and $\delta_{\nu+1}$-semistability and both
      $\delta_\nu$- and $\delta_{\nu+1}$-stability imply $\delta$-stability.
      In particular, there are maps ${\cal M}(\rho)^{\delta-ss}_{d/r/m}
      \lra
      {\cal M}(\rho)^{\delta_{\nu(+1)}-ss}_{d/r/m}$
      (``chain of flips", \cite{Th}).
\item For $\delta\in (0,\delta_0)$ and $(E,M,\tau)$ a
      $\delta$-semistable $\rho$-pair, the vector bundle $E$ must be semistable,
      and there is a morphism ${\cal M}(\rho)^{\delta-ss}_{d/r/m}\lra {\cal M}^{ss}_{d/r}$
      to the moduli space of semistable bundles of degree $d$ and rank $r$.
      Conversely, if $E$ is a stable bundle, then $(E,M,\tau)$ will be
      $\delta$-stable.
\item In the studied examples, there are only finitely many critical values,
      i.e., there is a $\delta_\infty$, such that the concept of $\delta$-semistability
      is constant in $(\delta_\infty,\infty)$. We refer to \cite{Bretal},
      \cite{Th}, \cite{OST}, \cite{SchHit2} and the examples
      for explicit discussions of this
      phenomenon. It would be interesting to know whether this is
      true in general or not, i.e., to check it for $\rho_{a,b,c}$.
\end{enumerate}
We note that in two of the examples, namely the example of
oriented framed modules and the example of Hitchin pairs, only a
parameter independent stability concept has been treated so far.
Our discussions will
therefore complete the picture in view of the above observations.
\subsubsection*{Asymptotic irreducibility}
Fix the representation $\rho$, the integers $d$ and $r$ as well
as the stability parameter $\delta\in\Q_{>0}$. Suppose that
$\rho$ is a direct summand of the representation $\rho_{a,b,c}$.
Since the estimate in Theorem~\ref{bound2} does not depend on the
integer $m$, we conclude that the set ${\cal S}$ of isomorphy
classes of vector bundles $E$, such that there exist an $m\in\Z$
and a $\delta$-semistable $\rho$-pair $(E,M,\tau)$ of type
$(d,r,m)$ is still bounded. The same goes for the set ${\cal
S}_\rho$ of vector bundles of the form $E_\rho$ with $[E]\in{\cal
S}$. Thus, there is a constant $m_0$, such that for every $m\ge
m_0$ and every $\delta$-semistable $\rho$-pair $(E,M,\tau)$ of
type $(d,r,m)$, one has
$$
{\rm Ext}^1(E_\rho,M)\q=\q H^1(E^\vee_\rho\otimes M)\q=\q 0.
$$
Our construction and standard arguments \cite{LP}, \S 8.5, now
show that the natural parameter space for $\delta$-semistable
$\rho$-pairs of type $(d,r,m)$ is a projective bundle over the
product of a smooth, irreducible, and quasi-projective quot scheme
and the Jacobian of degree $m$ line bundles. In particular, it is
smooth and irreducible. We infer
\begin{Thm}
Given the data $\rho$, $d$, $r$, and $\delta$ as above, there
exists a constant $m_0$, such that the moduli space ${\cal
M}(\rho)^{\delta-ss}_{d/r/m}$ is a normal and irreducible
quasi-projective variety for every $m\ge m_0$.
\end{Thm}
\begin{Rem}
Given $m^\p$, $m$ with $m-m^\p=l>0$ and a point $p_0\in X$, the
assignment $(E,M,\tau)\lma (E,M(lp_0),\tau^\p)$ with
$\tau^\p\colon E_\rho\lra M\subset M(lp_0)$ induces a closed
embedding
$$
{\cal M}(\rho)^{\delta-ss}_{d/r/m^\p} \hookrightarrow {\cal
M}(\rho)^{\delta-ss}_{d/r/m}.
$$
\end{Rem}
\subsection{Extension pairs}
Fix positive integers $0<s<r$, and let $F$ be the Grassmannian of
$s$-dimensional quotients of $\C^r$. An $F$-pair is thus a pair
$(E,q\colon E\lra Q)$ where $E$ is a vector bundle of rank $r$
and $q$ is a homomorphism onto a vector bundle $Q$ of rank $s$.
Setting $K:=\ker q$, we obtain a pair $(E,K)$ with $E$ as before
and $K\subset E$ a subbundle of rank $r-s$.
These objects were introduced by Bradlow and Garc\'\i a-Prada \cite{BG5}
as holomorphic extensions and called
(smooth) extension pairs in \cite{UW}.  In that work, $q$ is not required to be
surjective.
\par
We embed $F$ via the Pluecker embedding into $\Pe(\bigwedge^s\C^r)$, i.e.,
we consider the representation $\rho\colon \GL(r)\lra \GL(\bigwedge^s\C^r)$.
To describe the notion of $(\delta,\rho)$-semistability, we observe
that for points $[v]\in F\subset \Pe(\bigwedge^s\C^r)$, bases
$\ul{w}$ of $\C^r$, and weight vectors $\ul{\gamma}_1$ and $\ul{\gamma}_2$,
Equation~(\ref{addi}) holds true.
Furthermore, for a point $[v\colon \C^r\lra\C^s]\in F$,
a basis $\ul{w}=(w_1,...,w_r)$ of $\C^r$, and $i\in\{\, 1,...,r-1\,\}$
$$
\mu_{\rho}\bigl(\la(\ul{w}, \gamma^{(i)}), [v]\bigr)
\q=\q
i\dim \ker v-r\dim \bigl(\langle\, w_1,...,w_i\,\rangle\cap\ker v\bigr)
.
$$
Therefore, according to~(\ref{addii}),
an $F$-pair $(E,q\colon E\lra Q)$ is $(\delta,\rho)$-(semi)stable,
if and only if for every non trivial proper subbundle $E^\p$ of $E$
one has
$$
\mu(E^\p)+\delta{\rk(E^\p\cap\ker q)\over\rk E^\p}
\q (\le)\q
\mu(E)+\delta {\rk\ker q\over\rk E}.
$$
This is the same notion \cite{UW} provides for the extension pair
$(E,\ker q)$.
\subsection{Framed modules}
\label{FMOD}
The case of framed modules is one of the most thoroughly studied
examples of a decorated vector bundle problem (see, e.g., \cite{Br},
\cite{HT1}, \cite{Th}, \cite{Lu}, \cite{HL2}, \cite{HLF2}).
\par
First, we fix a positive integer $r$, an integer $d$, and a line
bundle $M_0$ on $X$ and look at the $\rho$-pairs of type
$(d,r,M_0)$ associated with the representation $\rho\colon
\GL(r)\lra \GL(\Hom(\C^s,\allowbreak \C^r))$, i.e., at pairs
$(E,\phi)$ consisting of a vector bundle $E$ of degree $d$ and
rank $r$ and a homomorphism $\phi\colon E\lra M_0^{\oplus s}$.
For the representation $\rho$, the Additivity
Property~(\ref{addi}) is clearly satisfied, and given a
non-trivial proper subbundle $E^\p$ of $E$ one has
$\mu_\rho(E^\p,\phi)=-\rk E^\p$ or $r-\rk E^\p$ if $E^\p\subset
\ker \phi$ or $\not \subset\ker \phi$, respectively.
\par
Given $\delta\in\Q_{>0}$, Equation~(\ref{addii}) thus shows that
$(E,\phi)$ is $\delta$-(semi)stable, if for every non-trivial
proper subbundle $E^\p$ of $E$
\begin{eqnarray*}
\mu(E^\p)\q(\le)\q \mu(E)-{\delta\over \rk E}, &\hbox{if}& E^\p\subset\ker\phi
\\
\mu(E^\p)-{\delta\over\rk E^\p}
\q(\le)\q \mu(E)-{\delta\over \rk E}, &\hbox{if}& E^\p\not\subset\ker\phi.
\end{eqnarray*}
Finally, one has the following result on the stability parameter $\delta$:
\begin{Lem}
Fix integers $d,r$, $r>0$, and a line bundle $M_0$.
The set of isomorphy classes of vector bundles $E$ for which there exist
a parameter $\delta\in\Q_{>0}$ and a $\delta$-semistable
$\rho$-pair of type $(d,r,M_0)$ of the form $(E,\phi)$
is bounded.
\end{Lem}
This is proved as Prop.~2.2.2.\ in \cite{OST}. From this
boundedness result, it follows easily that the set of isomorphy
classes of vector bundles of the form $\ker\phi$, $(E,\phi)$ a
$\rho$-pair of type $(d,r,M_0)$ for which there exists a
$\delta\in\Q_{>0}$ w.r.t.\ which it becomes semistable is
bounded, too. We infer
\begin{Cor}
\label{FMODi}
There exists a positive rational number $\delta_\infty$ such that
for every $\delta\ge\delta_\infty$ and every $\rho$-pair
$(E,\phi)$ of type $(d,r,M_0)$, the following conditions
are equivalent
\begin{enumerate}
\item[\rm 1.]
$(E,\phi)$ is $\delta$-(semi)stable.
\item[\rm 2.] $\phi$ is injective.
\end{enumerate}
\end{Cor}
Now, fix a vector bundle $E_0$ on $X$. Recall that a \it framed module
of type $(d,r,E_0)$ \rm is a pair $(E,\psi)$
consisting of a vector bundle $E$ of degree
$d$ and rank $r$ and
a non-zero homomorphism $\psi\colon E\lra E_0$.
Fix a sufficiently ample line bundle $M_0$ on $X$ and an embedding
$\iota\colon E_0\subset M_0^{\oplus s}$ for some $s$.
Therefore, any framed module $(E,\psi)$ of type $(d,r,E_0)$ gives
rise to the $\rho$-pair $(E,\phi:=\iota\circ\psi)$ of type
$(d,r,M_0)$, and the
$\rho$-pair $(E,\phi)$ is $\delta$-(semi)stable, if and only if
$(E,\psi)$ is a $\delta$-(semi)stable framed module
in the sense of \cite{HL2}.
Finally, a family of framed modules of type $(d,r,E_0)$
parameterized by $S$ is a triple $(E_S,\psi_S,{\frak N}_S)$
consisting of a rank $r$ vector
bundle $E_S$ on $S\times X$, a line bundle ${\frak N}_S$ on $S$,
and a homomorphism
$\psi_S\colon E_S\lra\pi_X^*E_0\otimes {\frak N}_S$ which is non trivial
on every fibre $\{s\}\times X$, $s\in S$.
Associate to such a family $(E_S,\psi_S, {\frak N}_S)$ the family
$(E_S,\kappa_S,{\frak N}_S,\phi_S)$ of $\rho$-pairs of type $(d,r,M_0)$
where $\kappa_S(s)=[M_0]$ for all $s\in S$
and $\phi_S=(\pi_X^*(\iota)\otimes\id_{\pi_S^*{\frak N}_S})\circ \psi_S$.
This exhibits the functor associating to a scheme $S$ the set of
equivalence classes of families of $\delta$-(semi)stable framed
modules of type $(d,r,M_0)$ as the subfunctor
of $\ul{\hbox{M}}(\rho)_{d/r/M_0}^{\delta-(s)s}$ of those
families $(E_S,\kappa_S,{\frak N}_S,\phi_S)$ where
$\kappa_S$ is the constant morphism $s\lma [M_0]$,
and the composite $E_S\lra \pi_X^*(M_0^{\oplus s})\otimes \pi_S^*{\frak N}_S
\lra
\pi_X^*(M_0^{\oplus s}/\iota(E_0))\otimes\pi_S^*{\frak N}_S$ vanishes.
Since all these conditions are closed conditions, the moduli spaces
of $\delta$-(semi)stable
framed modules on curves (\cite{Th}, \cite{HL2})
become closed subschemes of our moduli spaces ${\cal M}(\rho)^{\delta-(s)s}
_{d/r/M_0}$.
\begin{Rem}
We have used a slightly different, more general notion of family
than \cite{HL2}. This choice only destroys the property of being
a fine moduli space and does not affect the construction of the
moduli space of framed modules.
\end{Rem}
\subsection{Oriented framed modules}
We begin with the representations $\rho_1\colon \GL(r)\allowbreak \lra
\GL(\Hom(\C^s,S^r\C^r))$ and
$\rho_2\colon\GL(r)\allowbreak\lra\GL(\bigwedge^r\C^r)$,
and set $\rho:=\rho_1\oplus\rho_2$.
Fix line bundles $L_0$ and $M_0$.
Then, a $\rho$-pair of type $(L_0,r,M_0)$ is a triple
$(E,\phi,\sigma)$, consisting of a vector bundle $E$ of rank $r$
with $\det E\cong L_0$, a homomorphism
$\phi\colon S^r E\lra M_0^{\oplus s}$, and a homomorphism $\sigma
\colon\det E\lra M_0$.
Next, assume we are given a line bundle $N_0$ with $N_0^{\otimes r}=M_0$
and $t$ such that $s=\#\{\, (i_1,...,i_t)\,|\,i_j\in\{\, 0,...,r\,\},\
j=1,...,t,
\hbox{ and } \sum_{j=1}^ti_j=r\,\}$, i.e., $S^rN_0^{\oplus t}\cong
M_0^{\oplus s}$. Then, to any triple $(E,\psi,\sigma)$
where $E$ is a vector bundle of rank $r$ with $\det E\cong L_0$ and
$\psi\colon E\lra N_0^{\oplus t}$ and $\sigma\colon \det E\lra N_0^{\otimes r}$
are homomorphisms, we can associate the $\rho$-pair
$(E,S^r\psi,\sigma)$ of type $(L_0,r,M_0)$.
Observe that for any weighted filtration $(E^\bullet,\ul{\alpha})$
one has
$$
\mu_{\rho_2}\bigl(E^\bullet,\ul{\alpha};\sigma\bigr)=0
\qquad \hbox{and}\qquad
\mu_{\rho_1}\bigl(E^\bullet,\ul{\alpha};S^r\psi\bigr)
=r\cdot \mu_{\rho_1^\p}(E^\bullet,\ul{\alpha};\psi\bigr)
$$
where $\rho_1^\p\colon\GL(r)\lra\GL(\Hom(\C^t,\C^r))$.
Therefore, Proposition~\ref{ch3} and the discussion of framed modules
show
\begin{Lem}
\label{ostii}
Let $(E,\psi,\sigma)$ be a triple
where $E$ is a vector bundle of rank $r$ with $\det E\cong L_0$ and
$\psi\colon E\lra N_0^{\oplus t}$ and $\sigma\colon \det E\lra N_0^{\otimes r}$
are homomorphisms, and $\delta\in\Q_{>0}$.
Then, the following conditions are equivalent:
\begin{enumerate}
\item[\rm 1.]
The associated $\rho$-pair $(E,S^r\psi,\sigma)$ of type $(L_0,r,M_0)$
is $\delta$-semistable.
\item[\rm 2.]
One of the following three conditions is verified:
\begin{enumerate}
\item[\rm i.]
$E$ is a semistable vector bundle.
\item[\rm ii.] The homomorphisms $\psi$ and $\sigma$ are
               non-zero and there exists a positive rational number
               $\delta^\p\le r\delta$, such that $(E,\psi)$
               is a $\delta^\p$-semistable $\rho_1^\p$-pair of type
               $(L_0,r,N_0)$.
\item[\rm iii.] The homomorphism $\sigma$ vanishes and $(E,\psi)$
               is an $(r\cdot \delta)$-semistable $\rho_1^\p$-pair of type
               $(L_0,r,N_0)$.
\end{enumerate}
\end{enumerate}
\end{Lem}
We omit the ``polystable version" of this Lemma. In particular,
for $r\delta>\delta_\infty$ (cf.~Corollary~\ref{FMODi}), one finds
\begin{Cor}
\label{osti}
Let $(E,\psi,\sigma)$ be a triple
where $E$ is a vector bundle of rank $r$ with $\det E\cong L_0$ and
$\psi\colon E\lra N_0^{\oplus t}$ and $\sigma\colon \det E\lra N_0^{\otimes r}$
are homomorphisms, and $\delta>\delta_\infty/r$.
Then, the following conditions are equivalent:
\begin{enumerate}
\item[\rm 1.]
The associated $\rho$-pair $(E,S^r\psi,\sigma)$ of type $(L_0,r,M_0)$
is $\delta$-semistable.
\item[\rm 2.]
One of the following three conditions is verified:
\begin{enumerate}
\item[\rm i.]
$E$ is a semistable vector bundle.
\item[\rm ii.] The homomorphisms $\psi$ and $\sigma$ are
               non-zero and there exists a positive rational number
               $\delta^\p$, such that $(E,\psi)$
               is a $\delta^\p$-semistable $\rho_1^\p$-pair of type
               $(L_0,r,N_0)$.
\item[\rm iii.] The homomorphism $\sigma$ vanishes and $\psi$ is injective.
\end{enumerate}
\end{enumerate}
\end{Cor}
Now, we turn to the moduli problem we would like to treat.
For this, we fix a line bundle $L_0$ and a vector bundle
$E_0$. Then, an \it oriented framed module of type $(L_0,r,E_0)$ \rm
is a triple $(E,\eps,\psi)$ where $E$ is a vector bundle
of rank $r$ with $\det E\cong L_0$ and $\eps\colon \det E\lra L_0$
and $\psi\colon E\lra E_0$ are homomorphisms, not both zero.
The corresponding moduli problem was treated in \cite{OST}.
Over curves, we can recover it from our theory in the following way:
If $N_0$ is sufficiently ample, there are embeddings $\iota_1\colon
L_0\subset N_0^{\otimes r}$ and $\iota_2\colon E_0\subset N_0^{\oplus t}$.
Thus, setting $M_0:=N_0^{\otimes r}$,
we can define $\sigma:=\iota_1\circ\eps\colon \det E\lra M_0$
and $\phi:=S^r(\iota_2\circ\psi)\colon S^r E\lra M_0^{\oplus s}=S^r
N_0^{\oplus t}$
in order to get the
$\rho$-pair $(E,\phi,\sigma)$ of type $(L_0,r,M_0)$.
By Corollary~\ref{osti}, for $\delta\ge\delta_\infty/r$, the
$\rho$-pair $(E,\phi,\sigma)$ is $\delta$-semistable if and
only if $(E,\eps,\psi)$ is a semistable oriented framed module
in the sense of \cite{OST}.
\begin{Rem}
The corresponding stability concept can be recovered
via Proposition~\ref{ch3} and the characterisation
``stable=polystable+simple" (Remark~\ref{st=pst+s} ii).
\end{Rem}
We conclude by observing that applying Lemma~\ref{ostii} yields new semistability
concepts for oriented framed modules.
\subsection{Hitchin pairs}
The theory of Hitchin pairs or Higgs bundles is also a famous
example of a decorated vector bundle problem (\cite{Hi}, \cite{Si},
\cite{Fa}, \cite{Ni}, \cite{Yo}, \cite{Ha}, \cite{SchHit}).
\par
To begin with, we fix integers $d$ and $r>0$, a line bundle
$M_0$, and the representation
$\rho\colon \GL(r)\lra \GL(\End(\C^r)\oplus\C)$.
In this case, a $\rho$-pair of type $(d,r,M_0)$ is a triple
$(E,\phi,\sigma)$ consisting of a vector bundle $E$ of degree $d$ and rank
$r$, a twisted endomorphism $\phi\colon E\lra E\otimes M_0$, and
a section $\sigma\colon\Oh _X\lra M_0$.
\begin{Lem}
\label{Hit}
There is a positive rational number $\delta_\infty$, such that for
all $\delta\ge\delta_\infty$ and all $\rho$-pairs $(E,\phi,\sigma)$
of type $(d,r,M_0)$ the following conditions are equivalent:
\begin{enumerate}
\item[\rm 1.]
$(E,\phi,\sigma)$ is a $\delta$-(semi)stable $\rho$-pair
\item[\rm 2.]
for every non trivial subbundle $E^\p$ of $E$ with $\phi(E^\p)\subset
E^\p\otimes M_0$
$$
\mu(E^\p)\q(\le)\q \mu(E),
$$
and either $\sigma\neq 0$ or $\phi$ is not nilpotent, i.e.,
$(\phi\otimes \id_{M_0^{\otimes r-1}})\circ\cdots\circ\phi\neq 0$.
\end{enumerate}
\end{Lem}
\begin{proof} First, assume 1. Let $f\colon\C^r\lra\C^r$ be a
homomorphism. Call a sub vector space $V\subset \C^r$ \it
$f$-superinvariant\rm, if $V\subset\ker f$ and $f(\C^r)\subset V$.
\begin{Lem}
Let $[f,\eps]\in\Pe(\Hom(\C^r,\C^r)\oplus \C)$. Given a basis
$\ul{w}=(w_1,...,w_r)$ of $W$ and $i\in\{\, 1,...,r-1\,\}$, set
$W_{\ul{w}}^{(i)}:=\langle\, w_1,...,w_i\,\rangle$.
Then
\par
{\rm i)}
$\mu_\rho\bigl(\la(\ul{w},\gamma^{(i)}),[f,\eps]\bigr)=r$, if  $W_{\ul{w}}^{(i)}$
is not $f$-invariant.
\par
{\rm ii)} $\mu_\rho\bigl(\la(\ul{w},\gamma^{(i)}),[f,\eps]\bigr)=-r$,
if $W_{\ul{w}}^{(i)}$ is  $f$-superinvariant and $\eps=0$.
\par
{\rm iii)}
$\mu_\rho\bigl(\la(\ul{w},\gamma^{(i)}),[f,\eps]\bigr)=0$ in all the
other cases.
\end{Lem}
Now, let $(E,\phi,\sigma)$ be a $\rho$-pair of type $(d,r,M_0)$.
For any subbundle $E^\p$ of $E$ with $\phi(E^\p)\subset E^\p\otimes M_0$,
we find
$\mu_\rho(E^\p, (\phi,\sigma))\le 0$.
\begin{Cor}
Let $\delta\in \Q_{>0}$  and $(E,\phi,\sigma)$
a $\delta$-(semi)stable $\rho$-pair of type $(d,r,\allowbreak M_0)$.
Then $\mu(E^\p)(\le)\mu(E)$ for every non-trivial proper subbundle
$E^\p$ of $E$ with $\phi(E^\p)\subset E^\p\otimes M_0$.
\end{Cor}
This condition implies that
for every $\delta>0$, every $\delta$-semistable $\rho$-pair
$(E,\phi,\sigma)$ of type $(d,r,M_0)$, and every subbundle $E^\p$
of $E$
\begin{equation}
\label{nit}
\mu(E^\p)\q\le\q\max\biggl\{\, \mu(E),\mu(E)+{(r-1)^2\over r}\deg M_0\,\biggr\}.
\end{equation}
See, e.g., \cite{Ni}.
Therefore, the set of isomorphy classes of bundles $E$,
such that there exist a positive rational number $\delta$
and a $\delta$-semistable
$\rho$-pair of type $(d,r,M_0)$ of the form $(E,\phi,\sigma)$, is bounded.
\par
Now, the only thing we still have to show is that for every sufficiently large
positive rational number $\delta$ and every $\delta$-semistable
$\rho$-pair $(E,\phi,\sigma)$ of type $(d,r,M_0)$, such that $\sigma=0$,
the homomorphism $\phi$ can't be nilpotent.
First, let $(E,\phi,\sigma)$ be a $\rho$-pair of type $(d,r,M_0)$,
such that there exists a positive rational number $\delta$ w.r.t.\
which $(E,\phi,\sigma)$ is semistable and
such that $\phi$ is nilpotent.
Then, there is a filtration
$$
0=:E_0\subset E_1\subset \cdots\subset E_{s-1}\subset E_s:=E
$$
with $E_j\otimes M_0=\phi(E_{j+1})$, $j=0,...,s-1$.
It is clear by the boundedness result that the $E_j$'s occurring
in this way live
in bounded families, so that we can find a positive constant $C$
with
$$
d\rk E_j-\deg E_j r\q<\q C,\q j=1,...,s-1
$$
for all such filtrations.
One checks $\mu_\rho(E^\bullet, (1,...,1);(\phi,\sigma))=-r$, so that the
semistability assumption yields
\begin{eqnarray*}
0&\le& M\bigl(E^\bullet, (1,...,1)\bigr)+\delta
\mu_{\rho}\bigl(E^\bullet, (1,...,1); (\phi,\sigma)\bigr)
\\
&\le& (r-1) C-\delta r.
\end{eqnarray*}
This is impossible if $\delta\ge C$.
\par
To see the converse, let $(E,\phi,\sigma)$ be a $\rho$-pair satisfying 2.
Let $m_0:=\max\{\, 0,\allowbreak\deg M_0 (r-1)^2/r\,\}$.
Then, as before, $\mu(E^\p)\le\mu(E)+m_0$ for every non-trivial
proper subbundle $E^\p$
of $E$, i.e., $d\rk E^\p-r\deg E^\p\ge -m_0 r\rk E^\p\ge -m_0 (r-1)r$.
First, consider a weighted filtration $(E^\bullet,\ul{\alpha})$ such that
$\phi(E_j)\subset E_j\otimes M_0$, $j=1,...,s$.
Then, the condition that $\phi$ be not nilpotent if $\sigma=0$
implies $\mu_\rho(E^\bullet,\ul{\alpha};(\phi,\sigma))=0$,
so that $M(E^\bullet,\ul{\alpha})(\ge)0$ follows from 2.
Second, suppose that we are given a weighted filtration
$(E^\bullet, \ul{\alpha})$ such that, say, $E_{j_1},...,E_{j_t}$ are not
invariant under $\phi$, i.e., $\phi(E_{j_i})\not\subset E_{j_i}\otimes M_0$,
$i=1,...,t$, and $t>0$.
Let $\alpha:=\max\{\, \alpha_{j_1},...,\alpha_{j_t}\,\}$.
One readily verifies $\mu_{\rho}(E^\bullet,\ul{\alpha};(\phi,\sigma))\ge \alpha\cdot r$.
We thus find
\begin{eqnarray*}
M\bigl(E^\bullet,\ul{\alpha})
+\delta\mu_{\rho}\bigl(E^\bullet,\ul{\alpha}\bigr)
&\ge &\sum_{i=1}^t \alpha_{j_i} \bigl(d\rk E_{j_i}-r\deg E_{j_i}\bigr)
+ r\alpha\delta
\\
&\ge & - (r-1) r m_0\sum_{i=1}^t \alpha_{j_i} + r\alpha\delta
\\
&\ge & \bigl(-(r-1)^2rm_0+r\delta\bigr)\alpha,
\end{eqnarray*}
so that $M\bigl(E^\bullet,\ul{\alpha})+\delta\mu_{\rho}\bigl(E^\bullet,\ul{\alpha}\bigr)$
will be positive if we choose $\delta>(r-1)^2 m_0$.

\end{proof}
\begin{Ex}
\label{pari}
For small values of $\delta$, the concept of $\delta$-(semi)stability
seems to become rather difficult.
However, in the rank two case we have:
A $\rho$-pair $(E,\phi,\sigma)$ of type $(d,2,M_0)$ is $\delta$-(semi)stable
if for every line subbundle $E^\p$ of $E$ one has
\begin{enumerate}
\item $\deg E^\p(\le) d/2+\delta$,
\item $\deg E^\p(\le) d/2$ if $E^\p$ is invariant under $\phi$,
\item $\deg E^\p(\le) d/2-\delta$ if $E^\p=\ker\phi$,
      $\phi(E)\subset E^\p\otimes M_0$, and $\sigma=0$.
\end{enumerate}
\end{Ex}
Fix a line bundle $L$ on $X$.
We remind the reader \cite{SchHit} that a \it Hitchin pair of type
$(d,r,L)$ \rm is a triple $(E,\psi,\eps)$ where $E$ is a vector bundle
of degree $d$ and rank $r$, $\psi\colon E\lra E\otimes L$ is a twisted
endomorphism, and $\eps$ is a complex number.
Two Hitchin pairs $(E_1,\psi_1,\eps_1)$ and $(E_2,\psi_2,\eps_2)$
are called \it equivalent\rm, if there exist an isomorphism
$h\colon E_1\lra E_2$ and a non zero complex number $\la$
with $\la \psi_1=(h\otimes\id_L)^{-1}\circ\psi_2\circ h$ and
$\la \eps_1=\eps_2$.
We fix a point $x_0$ and choose $n$ large enough, so that $M_0:=L(nx_0)$
has a non trivial
global section. Fix such a global section $\sigma_0\colon
\Oh _X\lra M_0$ and an embedding $\iota\colon L\subset M_0$.
To every Hitchin pair $(E,\psi,\eps)$ of type $(d,r,L)$, we can assign
the $\rho$-pair $(E,\phi,\sigma)$ with $\phi:=(\id_E\otimes\iota)\circ\psi$
and $\sigma:=\eps\cdot \sigma_0$.
Note that this assignment is compatible with the equivalence
relations.
By Lemma~\ref{Hit}, for $\delta\ge\delta_\infty$, the $\rho$-pair
$(E,\phi,\sigma)$ is $\delta$-(semi)stable if and only if
$(E,\psi,\eps)$ is a (semi)stable Hitchin pair in the sense
of \cite{SchHit}.
Again, the above assignment carries over to families, so that
the general construction also yields a construction of
the moduli space of
semistable Hitchin pairs on curves, constructed in
\cite{SchHit} and \cite{Ha}. This space is a compactification
of the ``classical" Hitchin space \cite{Hi}, \cite{Fa}, \cite{Ni}.
\par
As we have seen, the semistability concept for Hitchin pairs is
parameter dependent in nature, though it might be difficult to describe
for low values of $\delta$. To illustrate that we get new semistable
objects for small values of $\delta$, let us look
at an
\begin{Ex}
i)
Let $x_0\in X$ be a point, and set $\Oh (1):=\Oh _X(x_0)$.
Define $E:=\Oh \oplus\Oh (1)$, and $\psi\colon E\lra E\otimes\Oh (1)
=\Oh (1)\oplus \Oh (2)$ as the homomorphism whose restriction to $\Oh $ is zero
and, moreover, the induced homomorphisms $\Oh (1)\lra \Oh (1)$ and
$\Oh (1)\lra\Oh (2)$ are the identity and zero, respectively.
First, consider the Hitchin pair $(E,\psi,1)$.
Then, the third condition in~\ref{pari} is void and the second condition
is satisfied. Indeed, a $\psi$-invariant subbundle $E^\p$ of $E$ of rank one
cannot be contained in $\Oh (1)$ whence $\deg E^\p\le 0< 1/2$.
Any other line subbundle $E^\p$ has degree at most one, and $E^\p:=\Oh (1)$
is a subbundle of degree exactly one.
The first condition then reads
$1(\le)1/2 + \delta$. In other words, $(E,\psi,1)$ is $\delta$-stable
for $\delta>1/2$, properly $(1/2)$-semistable, and not semistable
for $\delta<1/2$.
Finally, we claim that $(E,\psi,0)$ is properly $(1/2)$-semistable
(although $\psi$ is nilpotent).
For this, we only have to check the condition for $E^\p=\Oh $,
i.e., $0\le 1/2-1/2$, and this is clearly satisfied.
\par
ii) To see the r\^ole of $\delta$ in the whole theory,
let us look at Hitchin pairs of type $(1,2,{\omega_X})$.
Let $\delta_{\infty}$ be as in Lemma~\ref{Hit}.
For $\delta\ge\delta_\infty$, denote by ${\cal H}it_{{\omega_X}}$ the moduli
space of stable (in the usual sense) Hitchin pairs of type $(1,2,{\omega_X})$.
Let $\delta_0,...,\delta_m\in (0,\delta_\infty)$ be the critical values.
For $0< \delta<\delta_0$, the moduli space of $\delta$-stable Hitchin
pairs of type $(1,2,{\omega_X})$ equals $\Pe(\Oh _{\n}\oplus T_{\n})$, the
compactified cotangent bundle of $\n$, the moduli space of stable
rank two bundles of degree one.
Furthermore, let
${\cal M}^i_{\omega_X}$ be the moduli space of $\delta$-stable Hitchin
pairs of type $(1,2,{\omega_X})$ where $\delta\in (\delta_i,\delta_{i+1})$, $i=0,...,
m-1$,
and $\widetilde{\cal M}^i_{\omega_X}$ the moduli space
of $\delta_i$-semistable Hitchin pairs of type $(1,2,{\omega_X})$, $i=0,...,m$.
Between those spaces, we have morphisms
$$
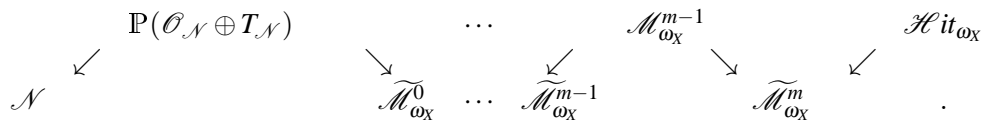
\begin{array}{ccccccccc}
          &         &\Pe(\Oh _\n\oplus T_\n)&                 &\cdots
          &         &{\cal M}_{\omega_X}^{m-1}\qquad\qquad& &{\cal H}it_{\omega_X}
\\
          &\swarrow &                      & \searrow &
         &\swarrow                 &\searrow&\swarrow
\\
\n        &         &                     &\qquad \widetilde{\cal M}^0_{\omega_X}
&\cdots&\widetilde{\cal M}^{m-1}_{\omega_X}& \qquad\qquad\widetilde{\cal M}^m_{\omega_X}&&.
\end{array}
$$
As in \cite{Th}, this is the factorization of the birational
correspondence $\Pe(\Oh _\n\oplus T_\n)\dasharrow {\cal H}it_{\omega_X}$ into
flips and is thus related to the factorization into blow ups and downs
(cf.\ \cite{Ha}).
\begin{Rem}[A.\ Teleman]
It might seem odd that we also obtain new semistability concepts for the
classical Higgs bundles $(E,\phi)$
where the semistability concept is known to be
parameter independent. In gauge theory, the reason is that, for studying
Higgs bundles, one fixes a flat
metric of infinite volume on the fibre $F=\End(\C^r)$ whereas we use
a metric of bounded volume induced by the embedding $\End(\C^r)\subset
\Pe(\End(\C^r)\oplus \C)$ which yields a different moment map. If we let the
parameter $\delta$ tend to infinity, we approximate the flat metric and
therefore recover the parameter independent semistability concept.
\end{Rem}
The related moduli problems of framed and oriented framed Hitchin
pairs discussed in \cite{Stup} and \cite{SchHit2} can also be dealt
with in our context. We leave this to the interested reader.
\end{Ex}
\subsection{Conic bundles}
Consider the representation $\rho\colon \GL(r)\lra \GL(S^2\C^r)$ and fix
a line bundle $M_0$ on $X$.
A $\rho$-pair of type $(d,r,M_0)$ is thus a pair $(E,\phi)$
consisting of a vector bundle $E$ of rank $r$ and degree $d$ and
a non-zero homomorphism $\rho\colon S^2 E\lra M_0$.
For $r\le 3$, these objects have been studied in \cite{GS}.
We apply Theorem~\ref{addniii} to analyze the notion of semistability,
using slightly different notation.
\par
To simplify the stability concept, we have to understand
the weights occurring for the action of $\SL(r)$ on $\Pe(S^2\C^r)$.
For this, let $[l]\in\Pe(S^2\C^r)$ be a point represented by the
linear form $l\colon S^2\C^r\lra \C$.
Set $I:=\{\, (i_1,i_2)\,|\, i_1,i_2\in\{\, 1,...,r\,\}, i_1\le i_2\,\}$.
For a basis $\ul{w}=(w_1,...,w_r)$ and $(i_1,i_2)\in I$,
we set $l(\ul{w})_{i_1i_2}:=l(w_{i_1}\otimes w_{i_2})$, so that the elements $l(\ul{w})_{i_1i_2}$, $(i_1,i_2)\in I$,
form a basis for $S^2\C^r$. We define a partial ordering on $I$,
by defining $(i_1,i_2)\preceq (j_1,j_2)$, if $i_1\le j_1$ and $i_2\le j_2$.
Furthermore, we define
$$
I(\ul{w},l):=\bigl\{\, (i_1,i_2)\in I\,|\, l(\ul{w})_{i_1i_2}\neq 0,
\hbox{ and $(i_1,i_2)$ is minimal w.r.t.\ ``$\preceq$"}\,\bigr\}.
$$
If $\# I(\ul{w},l)=1$, then one has the additivity property
(\ref{addi}) for all weight vectors $\ul{\gamma}_1$ and $\ul{\gamma}_2$.
In the other case, the cone of all weight vectors $(\gamma_1,...,\gamma_r)$
with $\gamma_1\le\cdots\le\gamma_r$ and $\sum\gamma_i=0$
becomes decomposed into subcones $C_{i_1i_2}(\ul{w},l)$,
$(i_1,i_2)\in I(\ul{w},l)$, where
$$
C_{i_1i_2}(\ul{w},l)
\q:=\q
\bigl\{\, (\gamma_1,...,\gamma_r)\,|\,
\gamma_{i_1}+\gamma_{i_2}\le \gamma_{i_1^\p}+\gamma_{i_2^\p}
\hbox{ for all } (i_1^\p,i_2^\p)\in I(\ul{w},l)\,\bigr\}.
$$
Then, (\ref{addi}) is still satisfied, if there is such a subcone
containing both $\ul{\gamma}_1$ and $\ul{\gamma}_2$. If one
chooses generators for these subcones, it therefore becomes
sufficient to compute the number
$\mu_{\rho}(\la(\ul{w},\ul{\gamma}),[l])$ for weight vectors
$\ul{\gamma}$ which are either of the form $\gamma^{(i)}$ or
belong to a set of generators for a cone $C_{i_1i_2}(\ul{w},l)$.
To see how this simplifies the concept of
$\delta$-(semi)stability, let us look at the cases $r=3$ and
$r=4$. 
\par
In the case $r=3$, one has $\#I(\ul{w}, l)=1$
unless $l(\ul{w})_{11}=0=l(\ul{w})_{12}$ and both
$l(\ul{w})_{22}$ and $l(\ul{w})_{13}$ are non-zero.
One checks that $C_{13}(\ul{w},l)$ is generated by $\gamma^{(1)}$ and
$\gamma^{(1)}+\gamma^{(2)}$ and that $C_{22}(\ul{w},l)$ is generated by
$\gamma^{(2)}$ and $\gamma^{(1)}+\gamma^{(2)}$.
To transfer this to our moduli problem,
let $E$ be a vector bundle of rank $3$ and $\tau\colon S^2E\lra M_0$
a non-zero homomorphism. Following \cite{GS}, given subbundles $F_1$
and $F_2$, we write $F_1\cdot F_2$ for the subbundle of $S^2 E$ generated
by local sections of the form $f_1\otimes f_2$ where $f_i$ is
a local section of $F_i$, $i=1,2$.
For any non-trivial proper subbundle $E^\p$ of $E$, one sets
\begin{itemize}
\item $c_\tau(E^\p)\q:=\q 2$,\qquad if $\tau_{|E^\p\cdot E^\p}\neq 0$,
\item $c_\tau(E^\p)\q:=\q 1$,\qquad if $\tau_{|E^\p\cdot E^\p}= 0$ and
                                 $\tau_{|E^\p\cdot E}\neq 0$, and
\item $c_\tau(E^\p)\q:=\q 0$,\qquad if $\tau_{|E^\p\cdot E}= 0$.
\end{itemize}
One checks
\begin{equation}
\label{wichti}
\mu_\rho\bigl(E^\p, \tau\bigr)\q=\q c_\tau(E^\p)\rk E-2\rk E^\p.
\end{equation}
Finally, call a filtration $E^\bullet:0\subset E_1\subset E_2\subset E$
with $\rk E_i=i$, $i=1,2$, \it critical\rm, if $\tau_{|E_1\cdot E_2}=0$,
and $\tau_{|E_1\cdot E}$ and $\tau_{|E_2\cdot E_2}$ are both non-zero.
Then
$$
\mu_\rho\bigl(E^\bullet, (1,1);\tau\bigr)\q=\q0.
$$
Putting everything together, we find
\begin{Lem}
A $\rho$-pair $(E,\tau)$ of type $(d,3,M_0)$ is $\delta$-(semi)stable
if and only if it satisfies the following two conditions
\begin{enumerate}
\item[\rm 1.]
For every non-zero proper subbundle $E^\p$ one has
$$
\mu(E^\p)-\delta {c_\tau(E^\p)\over \rk E^\p}\q (\le)\q
\mu(E)-\delta {2\over 3}.
$$
\item[\rm 2.]
For every critical filtration $0\subset E_1\subset E_2\subset E$
$$
\deg E_1+\deg E_2\q(\le)\q \deg E.
$$
\end{enumerate}
\end{Lem}
This is the stability condition formulated by G\'omez and Sols
\cite{GS}.
Next, we look at the case $r=4$. Set
$$
\nu(\ul{w},l):=\min\bigl\{\, i_1+i_2\,|\, l(\ul{w})_{i_1i_2}\neq 0, (i_1,i_2)
\in I\bigr\}.
$$
Suppose we are given a
linear form $l\colon S^2\C^4\lra \C$.
Then, for a basis $\ul{w}=(w_1,...,w_4)$, we have $\#I(\ul{w},l)=1$
except for the following cases
\begin{enumerate}
\item $\nu(\ul{w},l)=4$, $l(\ul{w})_{22}\neq 0$ and $l(\ul{w})_{13}\neq 0$,
\item $\nu(\ul{w},l)=4$, $l(\ul{w})_{22}\neq 0$, $l(\ul{w})_{13}= 0$, and $l(\ul{w})_{14}\neq 0$,
\item $\nu(\ul{w},l)=5$, $l(\ul{w})_{14}\neq 0$ and $l(\ul{w})_{23}\neq 0$,
\item $\nu(\ul{w},l)=5$, $l(\ul{w})_{14}\neq 0$, $l(\ul{w})_{23}= 0$, and $l(\ul{w})_{33}\neq 0$,
\item $\nu(\ul{w},l)=6$, $l(\ul{w})_{24}\neq 0$ and $l(\ul{w})_{33}\neq 0$.
\end{enumerate}
Straightforward computations show
\begin{Lem}
{\rm i)} In case {\rm 1.},
$C_{13}(\ul{w},l)$ is generated by $\gamma^{(1)}$, $\gamma^{(3)}$,
and $\gamma^{(1)}+\gamma^{(2)}$, and $C_{22}(\ul{w},l)$ by  $\gamma^{(2)}$,
$\gamma^{(3)}$, and $\gamma^{(1)}+\gamma^{(2)}$.
\par
{\rm ii)} In case {\rm 2.},
$C_{14}(\ul{w},l)$ is generated by $\gamma^{(3)}$, $\gamma^{(1)}+\gamma^{(3)}$,
and $\gamma^{(2)}+\gamma^{(3)}$, and $C_{22}(\ul{w},l)$ by  $\gamma^{(1)}$, $\gamma^{(2)}$,
$\gamma^{(1)}+\gamma^{(3)}$, and $\gamma^{(2)}+\gamma^{(3)}$.
\par
{\rm iii)} In case {\rm 3.},
$C_{14}(\ul{w},l)$ is generated by $\gamma^{(1)}$, $\gamma^{(2)}$,
and $\gamma^{(1)}+\gamma^{(3)}$, and $C_{23}(\ul{w},l)$ by $\gamma^{(2)}$,
$\gamma^{(3)}$, and $\gamma^{(1)}+\gamma^{(3)}$.
\par
{\rm iv)} In case {\rm 2.},
$C_{14}(\ul{w},l)$ is generated by $\gamma^{(2)}$, $\gamma^{(3)}$,
$\gamma^{(1)}+\gamma^{(2)}$, and $\gamma^{(1)}+\gamma^{(3)}$, 
and $C_{33}(\ul{w},l)$ by  $\gamma^{(1)}$, 
$\gamma^{(1)}+\gamma^{(2)}$, and $\gamma^{(1)}+\gamma^{(3)}$.
\par
{\rm v)} In case {\rm 5.},
$C_{24}(\ul{w},l)$ is generated by $\gamma^{(1)}$, $\gamma^{(2)}$,
and $\gamma^{(2)}+\gamma^{(3)}$, and $C_{33}(\ul{w},l)$ by $\gamma^{(1)}$,
$\gamma^{(3)}$, and $\gamma^{(2)}+\gamma^{(3)}$.
\end{Lem}
Now, let $(E,\tau)$ be a $\rho$-pair of type $(d,4,M_0)$.
For any non-zero, proper subbundle $E^\p$ of $E$, we define $c_\tau(E^\p)$
as before. One checks that (\ref{wichti}) remains valid.
Call a filtration $0\subset E_1\subset E_2\subset E_3\subset E$
with $\rk E_i=i$ \it critical of type \rm (I), (II), (III), (IV), (V), if
\begin{enumerate}
\item[(I)] $\tau_{|E_1\cdot E_2}=0$, and $\tau_{|{E_1\cdot E_3}}$
and $\tau_{|E_2\cdot E_2}$ are both non-zero;
\item[(II)] $\tau_{|E_1\cdot E_3}=0$, and $\tau_{|{E_1\cdot E}}$
and $\tau_{|E_2\cdot E_2}$ are both non-zero;
\item[(III)] $\tau_{|E_1\cdot E_3}=0$, $\tau_{|E_2\cdot E_2}=0$,
           and both $\tau_{|E_1\cdot E}$ and $\tau_{|E_2\cdot E_3}$ are
           non-zero;
\item[(IV)] $\tau_{|E_2\cdot E_3}=0$,
           and both $\tau_{|E_1\cdot E}$ and $\tau_{|E_3\cdot E_3}$ are
           non-zero;
\item[(V)] $\tau_{|E_1\cdot E}=0$, $\tau_{|E_2\cdot E_3}=0$,
           and both $\tau_{|E_2\cdot E}$ and $\tau_{|E_3\cdot E_3}$ are
           non-zero.
\end{enumerate}
respectively. In these cases, one has
\begin{itemize}
\item $\mu_{\rho}\bigl( 0\subset E_1\subset E_2\subset E, (1,1);\tau\bigr)
\q=\q -2$\qquad for type (I), (IV)
\item  $\mu_{\rho}\bigl( 0\subset E_1\subset E_3\subset E, (1,1);\tau\bigr)
\q=\q \phantom{-}0$\qquad for type (II), (III), (IV)
\item $\mu_{\rho}\bigl( 0\subset E_2\subset E_3\subset E, (1,1);\tau\bigr)
\q=\q \phantom{-}2$\qquad for type (II), (V).
\end{itemize}
Gathering all information,  we find
\begin{Lem}
The $\rho$-pair $(E,\tau)$ of type $(d,4,M_0)$
is $\delta$-(semi)stable
if and only if it satisfies the following two conditions
\begin{enumerate}
\item[\rm 1.]
For every non-zero proper subbundle $E^\p$ one has
$$
\mu(E^\p)-\delta {c_\tau(E^\p)\over \rk E^\p}\q (\le)\q
\mu(E)-{\delta \over 2}.
$$
\item[\rm 2.]
For every critical filtration $0\subset E_1\subset E_2\subset E_3\subset E$
$$
\begin{array}{lrcll}
\bullet & 4\deg E_1+4\deg E_2 &(\le) & 3\deg E-2,&
\hbox{\rm if it is of type (I), (IV)}
\\
\bullet & \deg E_1+\deg E_3 & (\le)&
\deg E, & \hbox{\rm if it is of type (II), (III), (IV)}
\\
\bullet &4\deg E_2+ 4\deg E_3&(\le)& 5\deg E+2, & \hbox{\rm if it
is of type \rm (II), (V)}.
\end{array}
$$
\end{enumerate}
\end{Lem}


\begin{thebibliography}{ABCD}

\bibitem{At} M.F.~Atiyah, \it Vector bundles on an elliptic curve\rm,
Proc.\ London Math.\ Soc.\ {\bf 7} (1957), 414-52.

\bibitem{Ba} D.~Banfield, \it Stable pairs and principal bundles\rm,  
Q.\ J.\ Math.\ \bf 51 \rm (2000), 417-36.

\bibitem{Bi} D.~Birkes, \it Orbits of linear algebraic groups\rm,
Ann.\ Math.\ {\bf 93} (1971), 459-75.

\bibitem{Br} S.B.~Bradlow, \it Special metrics and stability for holomorphic
bundles with global sections\rm, J.~Differential Geom.\ {\bf 33}
(1991), 169-213.

\bibitem{Bretal} S.B.~Bradlow, G.D.~Daskalopoulos,
O.~Garc\'\i a-Prada, R.~Wentworth, \it Stable augmented bundles over Riemann
surfaces \rm in \it Vector Bundles in Algebraic Geometry \rm (Durham 1993),
CUP 1995.

\bibitem{BG5} S.B.~Bradlow, O.~Garc\'\i a-Prada, \it Higher cohomology triples and
holomorphic extensions\rm, Comm.\ Anal.\ Geom.\ {\bf 3} (1995),  421-63.

\bibitem{HT2} S.B.~Bradlow, O.~Garc\'\i a-Prada, \it Stable triples,
equivariant bundles and dimensional reduction, \rm Math.~Ann.~{\bf 304}
(1996), 225-52.

\bibitem{ABCD} S.B.~Bradlow, O.~Garc\'\i a-Prada, I.~Mundet i Riera, \it 
Relative Hitchin--Kobayashi correspondences for principal pairs\rm,
Quart J.\ Math.\ \bf  54 \rm  (2003),  171-208.

\bibitem{UW} G.D.~Daskalopoulos, K.~Uhlenbeck, R.~Wentworth, \it
Moduli of extensions of holomorphic bundles on K\"ahler manifolds\rm,
Comm.\ Anal.\ Geom.\ {\bf 3} (1995), 479-522.

\bibitem{DH} I.V.~Dolgachev, Y.~Hu, \it Variation of geometric invariant
theory quotients\rm, Publ.\ I.H.E.S.\ {\bf 87} (1998), 5-56.

\bibitem{Fa} G.~Faltings, \it Stable $G$-bundles and projective
connections\rm, J.~Algebraic Geom.\ {\bf 2} (1993), 507-68.

\bibitem{FH} W.~Fulton, J.~Harris, \it Representation Theory --- A First Course\rm,
Springer 1991.

\bibitem{HT1} O.~Garc\'\i a-Prada, \it Dimensional reduction of stable
bundles,
vortices and stable pairs, \rm Int.\ J.~Math.~{\bf 5} (1994), 1-52.

\bibitem{Gies} D.~Gieseker, \it On the moduli of vector bundles on an
algebraic surface\rm, Ann.\ Math.\ {\bf 106} (1977), 45-60.

\bibitem{GS} T.L.~G\'omez, I.~Sols, \it Stability of conic bundles (with an
appendix by Mundet i Riera)\rm, Internat.\ J.\ Math.\ {\bf 11} (2000),
1027-55.

\bibitem{GS2} T.L.~G\'omez, I.~Sols, \it Stable tensor fields and moduli space of principal G-sheaves for classical groups\rm,
\tt math.AG/0103150\rm.

\bibitem{GS3} T.L.~G\'omez, I.~Sols, \it Moduli space of principal sheaves over projective
varieties\rm, \tt math.AG/0206277\rm, 40pp.

\bibitem{Gr} A.~Grothendieck, \it Sur la classification des fibr\'es
holomorphes sur la sph\`ere de Riemann\rm, Amer.\ J.\ Math.\ {\bf 79} (1957),
121-38.

\bibitem{Ha} T.~Hausel, \it Compactification of moduli of Higgs bundles\rm,
J.\ reine angew.\ Math.\ {\bf 503} (1998), 169-92.

\bibitem{Hi} N.J.~Hitchin, \it Self-duality equations on a compact Riemann
surface\rm, Proc.\ London Math.\ Soc.\ (3) {\bf 55} (1987), 59-126.

\bibitem{HL2} D.~Huybrechts, M.~Lehn, \it Stable pairs on curves and
surfaces\rm, J.\ Algebraic Geom.\ {\bf 4} (1995), 67-104.

\bibitem{HLF2} D.~Huybrechts, M.~Lehn, \it Framed modules and their
moduli\rm, Int.\ J.\ Math.\ {\bf 6} (1995), 297-324.

\bibitem{HL} D.~Huybrechts, M.~Lehn, \it The Geometry of Moduli Spaces
of Sheaves\rm, Vieweg 1997.

\bibitem{LP} J.~Le Potier, \it Lectures on Vector Bundles\rm, CUP 1997.

\bibitem{Lu} M.\ L\"ubke, \it The analytic moduli space
of framed vector bundles\rm,
J.\ reine angew.\ Math.\ {\bf 441} (1993), 45-59.

\bibitem{LT} M.~L\"ubke, A.~Teleman, \it The Kobayashi-Hitchin
Correspondence\rm,
World Scientific 1995.

\bibitem{Ma} M.~Maruyama, \it Moduli of stable sheaves \rm I, II,
J.\ Math.\ Kyoto Univ.\ {\bf 17} (1977), 91-126, {\bf 18} (1978), 557-614.

\bibitem{Mu2} D.\ Mumford, \it Projective invariants of projective structures
and applications \rm in \it Proc.\ Int.\ Congress Math.\ \rm (Stockholm 1962),
526-30.

\bibitem{GIT} D.~Mumford, \it Geometric Invariant Theory\rm,
Springer 1965.

\bibitem{MR} I.~Mundet i Riera, \it A Hitchin-Kobayashi correspondence
for K\"ahler fibrations\rm,  J.\ reine angew.\ Math.  \bf 528 \rm
(2000), 41-80.

\bibitem{NS1} M.S.~Narasimhan, C.S.~Seshadri, \it
Vector bundles on curves \rm in \it Algebraic Geometry \rm (Bombay 1968),
OUP 1969.

\bibitem{NS2} M.S.~Narasimhan, C.S.~Seshadri, \it
Moduli of vector bundles on a compact Riemann surface\rm, Ann.\ Math.\
{\bf 89} (1969), 14-51.

\bibitem{Ne} P.E.~Newstead, \it Introduction to Moduli Problems and
Orbit Spaces\rm, Springer 1978.

\bibitem{Ni} N.~Nitsure, \it Moduli space of semistable pairs on a curve\rm,
Proc.\ London Math. Soc. (3) {\bf 62} (1991), 275-300.

\bibitem{OST} Ch.~Okonek, A.~Schmitt, A.~Teleman, \it Master spaces
for stable pairs\rm, Topology {\bf 38} (1999), 117-39.

\bibitem{SchHit} A.~Schmitt, \it Projective moduli for Hitchin pairs\rm,
Int.\ J.\ Math.\ {\bf 9} (1998), 107-18; Erratum {\bf 11} (2000), 589.

\bibitem{SchHit2} A.~Schmitt, \it Framed Hitchin pairs\rm , Rev.\ roumaine
math.\ pures appl.\ \bf 45 \rm (2000), 681-711.

\bibitem{Sch4} A.~Schmitt, \it A useful semistability criterion\rm ,
Proc.\ American Math.\ Soc.\ {\bf 129} (2001), 1923-6.

\bibitem{Sch5} A.~Schmitt, \it Moduli problems of sheaves associated with
oriented trees\rm ,
Algebras and Representation Theory \bf 6 \rm (2003), 1-32.

\bibitem{SchPrin} A.~Schmitt, \it Singular principal bundles over higher-dimensional manifolds
and their moduli spaces\rm,  Int.\ Math.\ Res.\ Not.\ \bf
2002:23\rm, 1183-209.

\bibitem{SchQuiv} A.~Schmitt, \it 
Moduli for decorated tuples of sheaves and representation spaces for quivers\rm, \tt math.AG/0401173\rm.

\bibitem{Sesh} C.S.~Seshadri, \it Spaces of unitary vector bundles
on a compact Riemann surface\rm, Ann.\ Math.\ {\bf 85} (1967), 303-36.

\bibitem{Si} C.~Simpson, \it Moduli of representations of the fundamental
group of a smooth manifold \rm I, Publ.\ I.H.E.S. {\bf 79} (1994),
47-129.

\bibitem{Stup} M.-S.~Stupariu, \it The Kobayashi-Hitchin Correspondence
for Vortex-Type Equations Coupled with Higgs Fields\rm, PhD thesis, Z\"urich 1998.

\bibitem{Th} M.~Thaddeus, \it Stable pairs, linear systems and the Verlinde
formula\rm, Invent.\ Math.\ {\bf 117} (1994), 317-53.

\bibitem{Th2} M.~Thaddeus, \it Geometric Invariant Theory and Flips\rm,
J.\ Amer.\ Math.\ Soc.\ {\bf 9} (1996), 691-723.

\bibitem{Yo} K.~Yokogawa, \it Moduli of stable pairs\rm, J.\ Math.\
Kyoto Univ.\
{\bf 31} (1991), 311-27.

\end{thebibliography}
\end{document}